%% file: fredholm.tex
   \pdfoutput=1 

\documentclass{article}
\usepackage{excludeonly}

\usepackage{savesym}
\usepackage{scrextend}  
\usepackage{amsthm,amsmath}
\usepackage{amssymb,latexsym,graphicx}
\usepackage{amscd}
 \usepackage[all,cmtip]{xy}
\usepackage{tikz-cd}

\usepackage{accents}
\usepackage{cite}
\usepackage{mathtools}
\usepackage{stmaryrd} 

\usepackage{calligra}
\DeclareMathAlphabet{\mathcalligra}{T1}{calligra}{m}{n}
\DeclareMathAlphabet{\mathpzc}{OT1}{pzc}{m}{it}
\usepackage{hycolor}
\usepackage{xcolor}
\usepackage[
                       colorlinks=true,
                       linkcolor=black, 
                       citecolor=black, 
                       urlcolor=blue,
                     ]{hyperref}
\usepackage{soul} 

\usepackage{makeidx}

\usepackage[intoc]{nomentbl}
\makenomenclature

\usepackage{stackengine} 
\usepackage{booktabs} 

\input{joa-environments-english.tex}

\input{joa-latex-shortcuts.tex}

\begin{document}
\sloppy

\author{\quad Urs Frauenfelder \quad \qquad\qquad
             Joa Weber\footnote{
  Email: urs.frauenfelder@math.uni-augsburg.de
  \hfill
  joa@unicamp.br
  }
        %
        %
    \\
    Universit\"at Augsburg \qquad\qquad
    UNICAMP
}

\title{On the spectral flow theorem
         of Robbin-Salamon for finite intervals}

\date{\today}

\maketitle 
%


%
%

%





\begin{abstract}
In this article we consider operators of the form
$\p_s\xi+A(s)\xi$ where $s$ lies in an interval $[-T,T]$
and $s\mapsto A(s)$ is continuous.
Without boundary conditions these operators are not Fredholm.
However, using interpolation theory one can define suitable boundary
conditions for these operators so that they become Fredholm.
We show that in this case the Fredholm index is given by the spectral
flow of the operator path $A$.
\end{abstract}

\tableofcontents

\newpage 
\section{Introduction}

\subsection{Main results}

\begin{definition}\label{def:HS-pair}
A pair $H=(H_0,H_1)$ is called a \textbf{Hilbert space pair}
if $H_0$ and $H_1$ are both infinite dimensional Hilbert spaces such that
$H_1\subset H_0$ is a dense subset and the inclusion map
$\iota\colon H_1\to H_0$ is a compact linear map.
\\
Both Hilbert spaces in a Hilbert space pair
are separable by~\cite[Cor.\,A.5]{Frauenfelder:2024c}.
\end{definition}

Let $(H_0,H_1)$ be a Hilbert space pair.
An operator $A\in\Ll(H_1,H_0)$ is called \textbf{symmetric} if
\begin{equation}\label{eq:H_0-symmetric-INTRO}
   \INNER{A x}{y}_0
   =\INNER{x}{A y}_0 ,\quad
   \forall x,y\in H_1 .
\end{equation}

Note that while the notion of symmetric depends on the
inner product on $H_0$ it only depends on the inner product on $H_1$
up to equivalence. Namely, we call two inner products
$\INNER{\cdot}{\cdot}$ and $\INNER{\cdot}{\cdot}^\prime$ on a Hilbert
space $H$ \textbf{equivalent} if there exists a constant $c$ such that
$$
   \frac{1}{c}\norm{x}\le\norm{x}^\prime\le c\norm{x}
   ,\qquad
   \norm{x}:=\sqrt{\INNER{x}{x}}
   ,\quad
   \norm{x}^\prime:=\sqrt{\INNER{x}{x}^\prime} ,
$$
for every $x\in H$. One calls $\norm{\cdot}$ and $\norm{\cdot}^\prime$
the \textbf{induced norms}.

\medskip
The condition of being symmetric is kind of asymmetric.
While it depends on the $H_0$-inner product,
it only depends on the $H_1$-inner product up to equivalence
of norms.
A more symmetric notion which only depends on the equivalence classes
of the $H_1$- as well as the $H_0$-inner product is the following notion.

\begin{definition}
An element $A\in\Ll(H_1,H_0)$ is called \textbf{symmetrizable} if
there exists an inner product $\INNER{\cdot}{\cdot}$ on $H_0$
equivalent to the given inner product $\INNER{\cdot}{\cdot}_0$
such that $A$ is symmetric with respect to the new inner product
$\INNER{\cdot}{\cdot}$.
\end{definition}

We abbreviate by $\Ff=\Ff(H_1,H_0)\subset \Ll(H_1,H_0)$ the set of symmetrizable
Fredholm operators of index zero from $H_1$ to $H_0$.
We refer to the elements of $\Ff$ as \textbf{Hessians}.
We endow the set $\Ff$ with the subset topology inherited from
$\Ll(H_1,H_0)$.
We define $\Ff^*:=\{\A\in \Ff\mid \exists \A^{-1}\in\Ll(H_0,H_1)\}$.
To indicate invertibility visually we shall use for the elements of $\Ff^*$
the font $\A$.

\medskip
Taking adjoints gives rise to a bijection
(see Lemma~\ref{le:adjoint-map} for details)
\begin{equation}\label{eq:adjoint-map}
   *\colon \Ff(H_1,H_0)\to\Ff(H_0^*,H_1^*),\quad A\mapsto A^*
\end{equation}
which has the property $**=\Id_{\Ff(H_1,H_0)}$
and maps invertibles to invertibles.

\begin{remark}\label{rem:symmetrizable}
Note that~(\ref{eq:adjoint-map})
would not be true if one would replace symmetrizable by symmetric.
In fact, the adjoint of a symmetric operator $A\colon H_1\to H_0$ does
not need to be symmetric.
This is due to the asymmetric property of the symmetry condition
mentioned above.
Indeed the symmetry of $A\colon H_1\to H_0$ depends on the inner product on
$H_0$, while the symmetry of $A^*\colon H_0^*\to H_1^*$ depends on the
inner product on $H_1$ which can be used to identify $H_1$ with $H_1^*$.

In the following we will consider paths of Hessians.
Although in many applications one has paths of Hessians
which are symmetric for a fixed inner product on $H_0$,
and not for a time-dependent one as in the symmetrizable case,
the advantage of relaxing the symmetry condition to the
symmetrizability condition is that it gives a \emph{uniform} way to treat
paths of Hessians and the path of its adjoints. 
\end{remark}

\medskip
Let $I$ be an interval of the form
\begin{equation}\label{eq:I} 
   \R
   ,\quad
   I_-=\R_-=(-\infty,0]
   ,\quad I_+=\R_+=[0,\infty)
   ,\quad I_T=[-T,T] .
\end{equation}
Relevant \textbf{path spaces} are defined by
\begin{equation}\label{eq:P_1-I}
\begin{split}
   P_0(I)=P_0(I;H_0)
  &:=L^2 (I,H_0) ,
\\
   P_1(I)=P_1(I;H_1,H_0)
   &:= L^2 (I,H_1)\cap W^{1,2}(I,H_0) ,
\end{split}
\end{equation}
and these are Hilbert spaces with inner products
\[
   \INNER{v}{w}_{P_0}
   :=\int_I \INNER{v(s)}{w(s)}_0 \, ds
\]
and
\begin{equation}\label{eq:P_1-inner}
   \INNER{v}{w}_{P_1}
   :=\int_I \INNER{v^\prime(s)}{w^\prime(s)}_0 \, ds
   +\int_I \INNER{v(s)}{w(s)}_1 \, ds .
\end{equation}

\begin{definition}\label{def:Aa_I^*}
Denote the space of continuous paths of Hessians by
\begin{equation*}
   \Aa_I:=\{\text{$A\colon I\to \Ff$ continuous} \} .
\end{equation*}
The Hessian path spaces are defined by
\begin{equation*}
\begin{split}
   \Aa_{I_T}^*
   &:=\{A\in \Aa_{I_T}\mid \text{$\A_{-T}:=A(-T)$ and $\A_T:=A(T)$ are invertible} \}
\\
   \Aa_{I_+}^*
   &:=\{A\in \Aa_{I_+}\mid \text{
   $\A^+:=\lim_{s\to\infty} A(s)$ exists,
   $\A^+$ and $A(0)$ invertible} \}
\\
   \Aa_{I_-}^*
   &:=\{A\in \Aa_{I_-}\mid \text{
   $\A^-:=\lim_{s\to-\infty} A(s)$ exists,
   $\A^-$ and $A(0)$ invertible} \}
\\
   \Aa_\R^*
   &:=\{A\in \Aa_{I_-}\mid \text{
   $\A^\pm:=\lim_{s\to\pm\infty} A(s)$ exist
   and are invertible} \} .
\end{split}
\end{equation*}
\end{definition}

For $A\in\Aa^*_I$, where $I$ is one of the four interval types,
we define the bounded linear operator
\begin{equation}\label{eq:D_A}
   D_A\colon P_1(I)\to P_0(I)
   ,\quad
   \xi\mapsto \p_s\xi +A\xi .
\end{equation}

\begin{definition}[Projections]\label{def:projections}
Assume that $\A\in \Ff^*$.
Let $H_{1/2}=H_{1/2}(\A)$ be the interpolation space between the domain and
the co-domain of $\A$, namely between $H_1$ and $H_0$ in the case at hand.
We denote by
$$
   \pi_+^\A\in \Ll(H_\frac12)
   ,\qquad
   \pi_-^\A=\Id-\pi_+^\A\in \Ll(H_\frac12) ,
$$
the projection to the
positive eigenspaces of $\A$ along the negative ones,
respectively to the negative eigenspaces along the positive ones.
The images
$$
   H^\pm_\frac12(\A):=\pi_\pm^\A (H_\frac12)
$$
are complementary closed subspaces of $H_{1/2}$, as explained
in Section~\ref{sec:int-spaces-1/2}.
\end{definition}

\begin{definition}
For each of the four interval types $I$ we introduce
\textbf{augmented operators} as follows.
For $A\in \Aa_{I_T}^*$ we abbreviate $\A_{\pm T}:= A(\pm T)$ and define
\begin{equation}\label{eq:augmented-T}
\begin{split}
   \mathfrak{D}_A\colon P_1(I_T)
   &\to P_0(I_T)\times H^+_{\frac12}(\A_{-T}) \times H^-_{\frac12}(\A_T)
   =: \Ww(I_T;\A_{-T},\A_T)
   \\
   \xi
   &\mapsto \left(D_A\xi,\pi_+^{\A_{-T}}\xi_{-T},\pi_-^{\A_{T}}\xi_T\right).
\end{split}
\end{equation}
For $A\in \Aa_{I_+}^*$ we abbreviate $\A_0:= A(0)$
and define
\begin{equation*}
\begin{split}
   \mathfrak{D}_A\colon P_1(I_+)
   &\to P_0(I_+)\times H^+_{\frac12}(\A_0)
   =: \Ww(I_+;\A_0)
   \\
   \xi
   &\mapsto \left(D_A\xi,\pi_+^{\A_0}\xi_0\right).
\end{split}
\end{equation*}
For $A\in \Aa_{I_-}^*$ we abbreviate $\A_0:= A(0)$ and define
\begin{equation*}
\begin{split}
   \mathfrak{D}_A\colon P_1(I_-)
   &\to P_0(I_-)\times H^-_{\frac12}(\A_0)
   =: \Ww(I_-;\A_0)
   \\
   \xi
   &\mapsto \left(D_A\xi,\pi_-^{\A_0}\xi_0\right).
\end{split}
\end{equation*}
For $A\in \Aa_\R^*$ we define
\begin{equation*}
\begin{split}
   \mathfrak{D}_A\colon P_1(\R)
   &\to P_0(\R)
   \\
   \xi
   &\mapsto D_A\xi .
\end{split}
\end{equation*}
\end{definition}

The main result of this article is the following theorem
in which $I$ is any of the four interval types.
The proof uses~\cite[Thm.\,D]{Frauenfelder:2024c};
see Theorem~\ref{thm:proj}.

\begin{theoremABC}\label{thm:A}
For $A\in \Aa^*_I$ the augmented operator
$\mathfrak{D}_A$ is Fredholm and
$\INDEX{\mathfrak{D}_A}=\varsigma(A)$ where
$\varsigma(A)$ is the spectral flow of the path $A$ of Hessians.
\end{theoremABC}

\begin{remark}
In the case $I=\R$ Theorem~\ref{thm:A}
is the classical spectral flow theorem of Robbin and
Salamon~\cite{robbin:1995a}.
Strictly speaking, they proved the Fredholm property
under an additional assumption on the path $A$, namely they required
the existence of a weak derivative.
It was later shown by Rabier~\cite{Rabier:2004a} that such a weak
derivative is not needed to obtain the Fredholm property.

Although the special case $I=\R$ was known before our proof,
even in this case it differs rather from the proofs of Robbin-Salamon
and Rabier. While the Robbin-Salamon proof requires
the rather involved infinite dimensional transversality
theory~\cite{Abraham:1967a} to perturb the path of Hessians
to make it transverse in order to achieve only simple crossings
of the eigenvalues at zero,
our proof on concatenating finite intervals does not
require these techniques.
Instead we use elements $\lambda$ in the resolvent set
to shift the operators $D_A$ to operators $\mathfrak{D}_A^\lambda$,
defined by~(\ref{eq:frak-D-lambda}), for which the issue of non-simple
crossings of eigenvalues at zero can be avoided.
\end{remark}

Given $A\in \Aa^*_{I_T}(H_1,H_0)$, then
$-A^*\in\Aa^*_{I_T}(H_0^*,H_1^*)$ by~(\ref{eq:adjoint-map}).
We define the \textbf{adjoint} of $\mathfrak{D}_A$
in~(\ref{eq:augmented-T}) to be the augmented operator
associated to $-A^*$, i.e.
\begin{equation}\label{eq:augmented-T*}
\begin{split}
   \mathfrak{D}_A^*:=\mathfrak{D}_{-A^*}
   \colon P_1(I_T;H_0^*,H_1^*)
   &\to \Ww(I_T;-\A_{-T}^*,-\A_T^*) .
\end{split}
\end{equation}
Since the spectrum of an operator and its adjoint coincide, see
Lemma~\ref{le:spectrum-adjoint},
we have for the spectral flow
$$
   \varsigma(A)=\varsigma(A^*)=-\varsigma(-A^*)
$$ %
and therefore an immediate consequence of Theorem~\ref{thm:A} is the
index formula
$$
   \INDEX{\mathfrak{D}_A}=-\INDEX{\mathfrak{D}_A^*} .
$$

\smallskip\noindent
{\bf Acknowledgements.}
We would like to thank Fatih Kandemir for bringing to our attention
Rabier's article.
\\
UF~acknowledges support by DFG grant
FR~2637/4-1 and by Imecc Unicamp.

\subsection{Motivation and general perspective}

This article is part of our endeavor to find a general approach to Floer
homology as outlined in the section ``Motivation and general perspective''
in~\cite{Frauenfelder:2024c}.
With this goal in mind we therefore provide in the present article
a comprehensive study of these operators which play an important
role in a uniform approach to Floer homology.

Operators of the form $\p_s+A(s)$ for finite and half-infinite intervals
appear in the Hardy-approach to Lagrangian Floer gluing
of Tatjana Sim\v{c}evi\'{c}.


\newpage
\boldmath
\section{Preliminaries}
\unboldmath

\begin{notation}
The Kronecker symbol $\delta_{ij}$
is $1$ if $i=j$ and zero otherwise.
An operator is a bounded linear map.
\end{notation}

\boldmath
\subsection{$H$-self-adjoint operators}
\label{sec:H-sadj-operators-NEW}
\unboldmath

In Section~\ref{sec:H-sadj-operators-NEW}
let $H=(H_0,H_1)$ be a Hilbert space pair.
Let $\gf{h}\colon\N\to(0,\infty)$ be the growth function of $H$
and $H_\R=(H_r)_{r\in\R}$ the associated Hilbert $\R$-scale.

\boldmath
\subsubsection{$H$-self-adjointness}
\label{sec:H-self-adjointness}
\unboldmath

\begin{definition}
A bounded linear map $A\colon H_1\to H_0$ is called
\textbf{\boldmath$H$-self-adjoint} or, more precisely,
a self-adjoint Hilbert space pair operator,
if it is, firstly, \textbf{symmetric} in the sense that
\begin{equation}\label{eq:H_0-symmetric}
   \INNER{Ax}{y}_0=\INNER{x}{Ay}_0
   ,\qquad \forall x,y\in H_1
\end{equation}
and, secondly, a Fredholm operator of index zero.
\end{definition}

The requirement Fredholm of index zero guarantees non-emptiness of the
resolvent set $\R\setminus\spec A\not= \emptyset$, as we discuss right
below. Non-emptiness will be used over and over again in
Section~\ref{sec:semi-Fredholm} for perturbation arguments, see
e.g. Step~4 in the proof of Theorem~\ref{thm:Rabier}.

\begin{remark}[Why Fredholm of index zero is important]
\label{rem:Fredholm-index-zero}
As opposed to an operator acting on a Hilbert space, say $H_0\to H_0$,
the Fredholm requirement arises from domain $H_0$ and co-domain $H_1$
being different in the case at hand.

Suppose $A\colon H_1\to H_0$ is bounded,
but not Fredholm of index zero.
Then all reals lie in the spectrum
$$
   \R= \spec A:=\{\lambda\in\R\mid
  \text{$A-\lambda\iota\colon H_1\to H_0$ is not bijective}\}.
$$
To see equality suppose by contradiction that there is a real
$\lambda$ such that the bounded linear map
$A-\lambda\iota\colon H_1\to H_0$ is bijective and so, by the
open mapping theorem, admits a bounded inverse
$R_\lambda(A):=(A-\lambda\iota)^{-1}\colon H_0\to H_1$ called the
\textbf{resolvent of \boldmath$A$ at $\lambda$}$\,\notin\spec A$.
Thus $A-\lambda\iota$ is an isomorphism,
in particular $A-\lambda\iota$ is Fredholm of index zero.
But since $\iota$ is compact, so is $A$. Contradiction.
\end{remark}

\begin{remark}[Spectrum of $H$-self-adjoint operators is real and discrete]
\label{rem:compact-op}
In a Hilbert space pair both Hilbert spaces
are separable by~\cite[Cor.\,A.5]{Frauenfelder:2024c}.
If interpreted as an unbounded operator on $H_0$, then an
$H$-self-adjoint operator $A$ is a self-adjoint operator
$A\colon H_0\supset H_1\to H_0$ with dense domain $H_1$.
This, and what follows, is detailed in Appendix~\ref{sec:Hs.adj.ops}.

The spectrum of $A$ consists of infinitely many discrete
real eigenvalues~$a_\ell$, of finite multiplicity each,
which accumulate either at $+\infty$,
or at $-\infty$, or at both.
By Theorem~\ref{thm:ONB-basis},
there is a countable \textbf{orthonormal basis}
$\Vv(A)$ of $H_0$, see Definition~\ref{def:basis}, composed
of eigenvectors $v_\ell\in H_1$ of $A$.
The set of non-eigenvalues $\Rr(A):=\R\setminus \spec{A}$ is called
\textbf{resolvent set} of $A$. It is dense in~$\R$.

\smallskip\noindent
\textbf{a)~Invertible case.}
For invertible $H$-self-adjoint operators we use the boldface
letter $\A\colon H_1\to H_0$.
Accounting for multiplicities we enumerate the eigenvalues of $\A$
in increasing order and write them as a list with finite repetitions
\begin{equation}\label{eq:A^0-spectrum}
   \dots\le a_{-2}\le a_{-1}<0< a_1\le a_2\le\dots 
   ,\qquad
   \Ss(\A)=(a_\ell)_{\ell\in\Lambda} ,
\end{equation}
where the \textbf{eigenvalue index set} $\Lambda\subset \Z^*$,
counting multiplicities, is of the form
%
\begin{equation}\label{fig:Lambda-A}
\begin{tabular}{llll}
      
  & \text{Morse}
  & \text{co-Morse}
  & \text{Floer}
\\
\midrule
  $\Lambda$
  & \small $-\Lambda_-\cup\N$
  & \small $-\N\cup\Lambda_+$
  & \small $-\N\cup\N=:\Z^*$
\\
\midrule
  $\Lambda_-$
  & \small $\{\mu_-,\dots,2,1\}$ or $\emptyset$
  & \small 
  & \small $\N$
\\
  $\Lambda_+$
  & \small 
  & \small $\{1,2,\dots,\mu_+\}$ or $\emptyset$
  & \small $\N$
\end{tabular}
\end{equation}
%
The number of elements $\# \Lambda_-$ ($\# \Lambda_+$) is the
\textbf{Morse} (\textbf{co-Morse}) \textbf{index} of~$\A$.
Using the same index set $\Lambda$ we write the orthonormal
basis of $H_0$ in the form
\begin{equation}\label{eq:ONB-A}
   \Vv(\A)=\{v_\ell\}_{\ell\in\Lambda}\subset H_1
   ,\qquad
   \A v_\ell=a_\ell v_\ell ,
\end{equation}
where the eigenvalues accumulate on the set $\{-\infty,+\infty\}$;
see Theorem~\ref{thm:ONB-basis}.

\smallskip\noindent
\textbf{b)~Non-invertible case.}
In this case the only difference is that $A\colon H_1\to H_0$
has a nontrivial, but finite dimensional kernel for which one
chooses an ONB $\Vv(\ker{A})$. In the notation $A=0\oplus\A$
of Appendix~\ref{sec:Hs.adj.ops}, where $\A$ is invertible as in~a),
the eigenvalue list of $A$ and the corresponding ONB of $H_0$ are the unions
\begin{equation}\label{eq:Vv(A)-general}
   \Vv(A)\stackrel{\text{(\ref{eq:ONB-A})}}{=}\Vv(\A)\cup\Vv(\ker{A})
   ,\qquad
   \Ss(A)\stackrel{\text{(\ref{eq:A^0-spectrum})}}{=}\Ss(\A)\cup\{0\} .
\end{equation}
\end{remark}

\begin{definition}[$H$-self-adjoint operators come in three types]
We distinguish three \textbf{types} of $H$-self-adjoint
operators $A\colon H_1\to H_0$; cf. Remark~\ref{rem:compact-op}.
\begin{enumerate}
\item
  \textbf{Morse.} Finitely many negative, infinitely many positive
  eigenvalues.
\item
  \textbf{Co-Morse.} Finitely many positive, infinitely many negative
  eigenvalues.
\item
  \textbf{Floer.} Infinitely many negative and positive 
  eigenvalues.
\end{enumerate}
\end{definition}

\boldmath
\subsubsection{Banach adjoint}
\unboldmath

\begin{definition}\label{def:Banach-adjoint}
Let $A\in\Ll(H_1,H_0)$.
For a dual space element $\eta\in H_0^*:=\Ll(H_0,\R)$, the
image $A^*\eta\in H_1^*:=\Ll(H_1,\R)$ under the Banach \textbf{adjoint}
$$
   A^*\colon H_0^*\to H_1^*
$$
is characterized by
$$
   (A^*\eta)\xi=\eta(A\xi)
   ,\quad \forall \xi\in H_1 .
$$
The inner products on the dual spaces
are defined via the musical isomorphisms
$$
   \flat_0\colon H_0\to H_0^*,\quad \xi\mapsto \INNER{\xi}{\cdot}_0
   ,\qquad
   \flat_1\colon H_1\to H_1^*,\quad \xi\mapsto \INNER{\xi}{\cdot}_1 ,
$$
by
$$
   \INNER{\cdot}{\cdot}_0^*
   :=\bigl\langle{\flat_0}^{-1}\cdot , {\flat_0}^{-1}\cdot\bigr\rangle_0
   ,\qquad
   \INNER{\cdot}{\cdot}_1^*
   :=\bigl\langle{\flat_1}^{-1}\cdot , {\flat_1}^{-1}\cdot\bigr\rangle_1 .
$$
\end{definition}

\begin{lemma}\label{le:spectrum-adjoint}
Let $B\in\Ll(H_1,H_0)$.
Then $\spec B=\spec B^*$.
\end{lemma}

\begin{proof}
We first show that $B$ is invertible iff $B^*$ is invertible.
Assume $B\colon H_1\to H_0$ is invertible. This means that there is
$C\in\Ll(H_0,H_1)$ such that $BC=\Id_{H_0}$ and $CB=\Id_{H_1}$.
Applying $*$ to these equations we get $C^*B^*=\Id_{H_0^*}$ and $B^*C^*=\Id_{H_1^*}$.
Hence we have shown that invertibility of $B$ implies invertibility of
$B^*$. Hence invertibility of $B^*$ implies invertibility of $B^{**}$,
but $B^{**}=B$. Therefore invertibility of $B$ and $B^*$ are equivalent.

Observe that $\lambda\in\spec B$ iff $B-\lambda\iota\colon H_1\to H_0$
is not invertible.
As we have just seen this is equivalent that
$B^*-\lambda\iota^*\colon H_0^*\to H_1^*$ is not invertible.
This shows that the spectrum of $B$ coincides with the spectrum of $B^*$.
\end{proof}

\begin{lemma}\label{le:adjoint-map}
Taking adjoints gives rise to a bijection
$$
   *\colon \Ff(H_1,H_0)\to\Ff(H_0^*,H_1^*),\quad A\mapsto A^*
$$
which has the property $**=\Id_{\Ff(H_1,H_0)}$ and maps invertibles to invertibles.
\end{lemma}

\begin{proof}
That $*$ maps invertibles to invertibles holds by Lemma~\ref{le:spectrum-adjoint}.
  By~\cite[\S 16 Thm.\,4]{Muller:2007a} an operator $A\colon H_1\to H_0$
  is Fredholm iff $A^*\colon H_0^*\to H_1^*$ is Fredholm. In our
  case $\INDEX{A^*}=-\INDEX{A}=0$.

We first discuss the case when $A$ is invertible:
After replacing the inner product on $H_1$ and $H_0$ by equivalent
ones, we can assume without loss of generality, that $A$ is a
symmetric isometry.
Such inner products are referred to as $A$-adapted and the existence
is discussed around~(\ref{eq:1A}).
Using the $A$-adapted inner products we can naturally identify $H_0^*$
with $H_0$ and $H_1^*$ with $H_{-1}$ and $A^*$ becomes a symmetric
isometry $H_0\to H_{-1}$ as explained by Lemma~\ref{le:A*-A}.
This shows that $A^*$ is in $\Ff(H_0^*,H_1^*)$ and is invertible as well.

It remains to discuss the case when $A$ is not invertible:
Choose $\lambda$ in the resolvent set $\Rr(A)$.
Then $\A_\lambda:=A-\lambda\iota$ is an invertible element in $\Ff(H_1,H_0)$.
By the discussion before $\A_\lambda^*\in \Ff(H_0^*,H_1^*)$.
Using that $\A_\lambda^*=A^*-\lambda\iota^*$ we conclude that
$\A^*$ is in $\Ff(H_0^*,H_1^*)$ as well.

Since $H_0$ and $H_1$ are Hilbert spaces, they are in particular
reflexive so that we have the canonical isomorphisms
$H_0=H_0^{**}$ and $H_1=H_1^{**}$ which
does not depend on the choice of any inner product,
so that $A^{**}$ naturally becomes $A$.
\end{proof}

\boldmath
\subsection{The interpolation classes $H_{\frac12}^\pm(\A)$}
\label{sec:int-spaces-1/2}
\unboldmath

For a pair of Hilbert spaces $H=(H_0,H_1)$ let $H_{1/2}$ be the
$\R$-scale interpolation space of $H_0$ and $H_1$ as defined
in~(\ref{eq:H_r}) for $r=1/2$.
The construction of the interpolation space $H_{1/2}$
uses the $0$-inner product on $H_0$ and the 
$1$-inner product on $H_1$ to get a $\frac12$-inner product.
A useful formula, in terms of a pair growth function and a
scale basis, is~(\ref{eq:isometric}).

A consequence of the Stein-Weiss interpolation theorem is that if
we replace the inner products by equivalent ones, say a
$0^\prime$- and a $1^\prime$-inner product, then we obtain
on $H_\frac12$ as well an equivalent $\frac12^\prime$-inner product.
To see this abbreviate by $H_0^\prime$ the vector space
$H_0$ endowed with the $0^\prime$-inner product
and analogously for $H_1^\prime$.
Interpret the identity map as a map $\Id\colon H_0\to H_0^\prime$.
The identity map restricts to a map $H_1\to H_1^\prime$.
Since the $0$- and $0^\prime$-inner products are equivalent
there exists a constant $c_0$ such that $\norm{\Id}_{\Ll(H_0,H_0^\prime)}\le c_0$.
For the same reason there exists a 
constant $c_1$ such that $\norm{\Id}_{\Ll(H_1,H_1^\prime)}\le c_1$.
It follows from the Stein-Weiss interpolation theorem, see
e.g.~\cite[5.4.1 p.115, $U=V=\N$, $p=2$, $\theta=\frac12$]{Bergh:1976a}
or~\cite[App.\,B]{Frauenfelder:2024c},
that the identity map maps $H_\frac12$ to $H_\frac12^\prime$ and satisfies
$\norm{\Id}_{\Ll(H_\frac12,H_\frac12^\prime)}\le \sqrt{c_0c_1}$.
Interchanging the roles of $H_0$ and $H_0^\prime$
shows that the restriction of the identity to $H_\frac12$
actually is an isomorphism between $H_\frac12$ and~$H_\frac12^\prime$.

\begin{definition}\label{def:A-adapted}
Assume that $\A\in\Ff^*(H_1,H_0)$ is a symmetrizable invertible
bounded linear map from $H_1$ to $H_0$.
We say that equivalent inner products $1^\prime$ on $H_1$ and
$0^\prime$ on $H_0$ are \textbf{\boldmath$\A$-adapted} if $\A$ is an
isometry with respect to the inner products $1^\prime$ and $0^\prime$
and symmetric with respect to the inner product $0^\prime$.
\end{definition}

\textbf{Existence.} Note that $\A$-adapted inner products always exist:
Indeed since $\A\colon H_1\to H_0$ is symmetrizable there exists an
inner product $0^\prime$ on $H_0$ such that $\A$ is symmetric with
respect to the inner product $0^\prime$.
Now define the $1^\prime$ inner product on $H_1$
as the pull-back of the $0^\prime$-inner product on $H_0$,
i.e.
\begin{equation}\label{eq:1A}
   \INNER{\xi}{\eta}_{1^\prime}
   =\INNER{\A\xi}{\A\eta}_{0^\prime}
\end{equation}
for all $\xi,\eta\in H_1$.
Since $\A$ is invertible the $1^\prime$-inner product is equivalent
to the $1$-inner product on $H_1$.
By construction of the $1^\prime$-inner product
$\A$ becomes an isometry with respect to the
$1^\prime$- and $0^\prime$-inner products.

\medskip
\textbf{Spectral decomposition.}
An operator $\A\in\Ff^*$ gives rise to a decomposition 
of interpolation space into two closed subspaces
\begin{equation}\label{eq:H_1/2}
   H^\pm_\frac12 (\A)
   :=\pi_\pm^\A(H_\frac12)
   ,\quad
   H_\frac12=H^-_\frac12 (\A)\oplus H^+_\frac12 (\A)
   ,\quad
   (\pi^\A_\pm)^2=\pi^\A_\pm\in\Ll(H_\frac12) ,
\end{equation}
corresponding to the negative and the positive eigenspaces of $\A$.

\smallskip
\textbf{Case 1} (Symmetric isometry)\textbf{.}
We first explain the construction of the spaces $H_\frac12^\pm (\A)$ in
the special case where $\A\colon H_1\to H_0$ is a symmetric isometry.
In this case choose an orthonormal basis $\Vv(\A)=\{v_\ell\}_{\ell\in\Lambda}$
of $H_0$ as in~(\ref{eq:ONB-A}), in particular consisting of eigenvectors, namely
$\A v_\ell=a_\ell v_\ell$.
The basis is orthogonal in $H_1$, indeed
$\INNER{v_\ell}{v_k}_1=\INNER{\A v_\ell}{\A v_k}_0=a_\ell a_k \delta_{\ell k}$.
So the $H_1$-lengths are given by 
\begin{equation}\label{eq:1-length-v_ell}
   \norm{v_\ell}_1
   =\abs{a_\ell} .
\end{equation}
This basis is also orthogonal in $H_{\frac12}$
and the $H_\frac12$-lengths are given~by\footnote{
  By~(\ref{eq:inner_r}) for $\xi=\eta=v_\ell$ (the growth
  operator is $Tv_\ell=a_\ell^{-2}v_\ell$ since $\A$ is an isometry).
  }
\begin{equation}\label{eq:1/2-length-v_ell}
   \norm{v_\ell}_\frac12=\abs{a_\ell}^\frac12 .
\end{equation}
Now consider the projections $\pi^\A_\pm\colon H_\frac12\to H_\frac12$
to the positive/negative eigenspaces of $\A$, defined by
$$
   \pi^\A_+ v_\ell
   =\begin{cases}
      v_\ell&\text{, $\ell>0$,}\\
      0&\text{, $\ell<0$,}\\
   \end{cases}
   \qquad\quad
   \pi^\A_- v_\ell
   =\begin{cases}
      0&\text{, $\ell>0$,}\\
      v_\ell&\text{, $\ell<0$,}\\
   \end{cases}
$$
for every eigenvector $v_\ell\in\Vv(\A)$.
The definition shows that $(\pi^\A_\pm)^2=\pi^\A_\pm$, so image and
fixed point set coincide.
But the latter is a closed subspace by continuity.

\smallskip
\textbf{Case 2} (Symmetrizable invertible)\textbf{.}
Now consider the case where the bounded linear map
$\A\colon H_1\to H_0$ is still invertible, but only symmetrizable.
In this case we can replace the $0$- and $1$-inner products on $H_0$
and $H_1$ respectively, by $\A$-adapted ones, say $0^\prime$ and $1^\prime$,
as explained after Definition~\ref{def:A-adapted}.
The space $H_\frac12$ gets endowed with an equivalent
$\frac12^\prime$-inner product, too.
For the new inner products $0^\prime$, $1^\prime$, and $\frac12^\prime$
we obtain projections $\pi^\A_\pm$ as above.
By equivalence of the new and the original inner product on $H_\frac12$
the projections
\begin{equation}\label{eq:pi}
   \pi^\A_\pm\colon H_\frac12 \to H_\frac12
\end{equation}
are still continuous with respect to the original inner product, 
although in general they will not be orthogonal any more.
The images of $\pi_-^\A$ and $\pi_+^\A$ are 
complementary closed subspaces and we get~(\ref{eq:H_1/2}).

\newpage
\section{Spectral flow}\label{sec:spectral-flow}

Inspired by Hofer, Wysocki, and Zehnder~\cite[p.\,216]{Hofer:1998a}
we define the spectral flow as follows.
For the spaces $\Aa_I^*$ see~Definition~\ref{def:Aa_I^*}.

\subsubsection*{Finite intervals}

Given $A\in\Aa_{I_T}^*$, we consider the invertible operator
$\A_{-T}:=A(-T)\colon H_1\to H_0$ and we write its spectrum,
similar to~(\ref{eq:A^0-spectrum}) but now denoting the
{\color{orange} non-eigenvalue~$0$} by $a_0^{-T}$, in the form
\begin{equation}\label{eq:A_-T-spectrum}
   \dots\le a^{-T}_{-2}\le a^{-T}_{-1}<{\color{orange}\underbrace{a_0^{-T}}_{:=0}}
   <a^{-T}_1\le a^{-T}_2\le\dots
   ,\quad
   \Ss(\A_{-T})=(a^{-T}_\ell)_{\ell\in\Lambda} .
\end{equation}
Extend these eigenvalues and ${\color{orange}a_0^{-T}:=0}$
to continuous functions
$a_\ell \colon[-T,T]\to\R$ for $\ell\in\Lambda\cup\{0\}$ satisfying,
for any $s\in[-T,T]$,
the following conditions under inclusion of the zero function
${\color{orange}[-T,T]\ni s\mapsto 0}$, namely
\begin{itemize}\setlength\itemsep{0ex} 
\item[\rm (i)]
   $\dots\le a_{-2}(s)\le a_{-1}(s)
   \le a_0(s)
   \le a_1(s)\le a_2(s)\le\dots$
\item[\rm (ii)]
  $\Ss(A(s))\cup \{0\}=(a_\ell(s))_{\ell\in\Lambda\cup\{0\}}$.
\end{itemize}
Since eigenvalues depend continuously on the operator,
the functions $a_\ell$ exist and are uniquely determined
by these two conditions.

\begin{definition}[Spectral flow - finite interval]\label{def:SF-finite}
Let $A\in\Aa_{I_T}^*$ be a path of Hessians.
The \textbf{spectral flow} $\varsigma(A)\in\Z$ is defined by
\begin{equation}\label{eq:SF-finite}
   \varsigma(A)
   :=-i\text{ if $a_i(T)=0$.}
\end{equation}
This is the net count of eigenvalues that change from negative to positive.
\end{definition}

\begin{figure}[h]
  \centering
  \includegraphics
                             {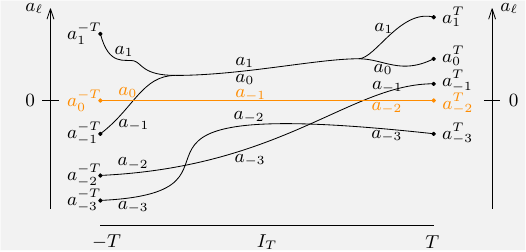}
  \caption{Spectral flow $\varsigma(A)=2$ along $[-T,T]$}
   \label{fig:fig-specflow-I_T-0}
\end{figure}

\begin{lemma}\label{le:SF-axioms}
The map $\Aa_{I_T}^*\to\Z$, $A\mapsto \varsigma(A)$, has
the following properties.
\begin{labeling}{\texttt{(Normalization)}}\setlength\itemsep{0ex} 
\item[\rm {\texttt{(Homotopy)}}]
  $\varsigma$ is constant on the connected components of $\Aa_{I_T}^*$.
\item[\rm{\texttt{(Constant)}}]
  If $A$ is constant, then $\varsigma(A)=0$.
\item[\rm{\texttt{(Direct Sum)}}]
  $\varsigma(A_1\oplus A_2)=\varsigma(A_1)+\varsigma(A_2)$.  
\item[\rm{\texttt{(Normalization)}}]
  For $H_1=H_0=\R$ and $A(s)=\arctan(s)$, it holds $\varsigma(A)=1$.
\item[\rm{\texttt{(Catenation)}}]
  If $A=A_\ell\# A_r$, then $\varsigma(A)=\varsigma(A_\ell)+\varsigma(A_r)$.
\end{labeling}
\end{lemma}

The first four properties guarantee uniqueness and the first three imply
\texttt{(Catenation)}, see~\cite[\S\,4]{robbin:1995a}.
Thus $\varsigma$ is the spectral flow as defined in~\cite{robbin:1995a}.

\begin{proof}
\texttt{(Homotopy)} Assume that $A_0$ and $A_1$ lie in the same
connected component of $\Aa_{I_T}^*$. This means that there exists
a homotopy $\{A_r\}_{r\in[0,1]}\subset \Aa_{I_T}^*$ between them.
Consider the map $[0,1]\to\Z$, $r\mapsto \varsigma(A_r)$.
By continuous dependence of the eigenvalues this map is continuous
and since its image is discrete, the map is constant.
In particular $\varsigma(A_0)=\varsigma(A_1)$ and therefore the
homotopy property is proved.
\\
\texttt{(Constant)} In this case $a_\ell(s)\equiv a_\ell(-T)$ $\forall s,\ell$,
in particular we have $a_0(T)=a_0(-T)=0$, hence $\varsigma(A)=0$.
\\
\texttt{(Direct Sum)}
The net number of eigenvalues of the direct sum $A_1\oplus A_2$
crossing zero corresponds to the sum of the net number of eigenvalues
of $A_1$ crossing zero and $A_2$ crossing zero. Therefore the direct
sum property holds.
\\
\texttt{(Normalization)} Initially, since $a_0(-T)=0$ and $\arctan(-T)<0$,
we have $\arctan(-T)=a_{-1}(-T)$.
At the end, since $\arctan(T)>0$, we have
$0=a_{-1}(T)$. Thus $\varsigma(A)=-(-1)=1$.
\end{proof}

\subsubsection*{Half-infinite forward interval}

Given $A\in\Aa_{I_+}^*$, we consider the invertible operator
$\A_0:=A(0)\colon H_1\to H_0$ and we write its spectrum,
similar to~(\ref{eq:A^0-spectrum}) but now denoting the
{\color{orange} non-eigenvalue~$0$} by $a_0^0$, in the form
\begin{equation}\label{eq:A_0-spectrum}
   \dots\le a^0_{-2}\le a^0_{-1}
   <{\color{orange}\underbrace{a_0^{0}}_{:=0}}
   < a^0_1\le a^0_2\le\dots 
   ,\quad
   \Ss(\A_0)=(a^0_\ell)_{\ell\in\Lambda} .
\end{equation}
Extend these eigenvalues and ${\color{orange}a_0^{0}:=0}$
to continuous functions
$a_\ell \colon[0,\infty)\to\R$ for $\ell\in\Lambda\cup\{0\}$ satisfying,
for any $s\in[-T,T]$,
the following conditions under inclusion of the zero function
${\color{orange}[-T,T]\ni s\mapsto 0}$, namely 
\begin{itemize}\setlength\itemsep{0ex} 
\item[\rm (i)]
   $\dots\le a_{-2}(s)\le a_{-1}(s)\le a_{0}(s)\le a_1(s)\le a_2(s)\le\dots$
\item[\rm (ii)]
  $\Ss(A(s))\cup\{0\}=(a_\ell(s))_{\ell\in\Lambda\cup\{0\}}$.
\end{itemize}
Since the limit $\lim_{s\to\infty} A(s)=:\A^+$ exists so do the limits
$\lim_{s\to\infty} a_\ell(s)=: a_\ell(\infty)$ for every $\ell\in\Lambda\cup\{0\}$
and they satisfy
\begin{itemize}\setlength\itemsep{0ex} 
\item[\rm (iii)]
   $\dots\le a_{-2}(\infty)\le a_{-1}(\infty)\le
   a_{0}(\infty)\le a_1(\infty)\le a_2(\infty)\le\dots$
\item[\rm (iv)]
  $\Ss(\A^+)=(a_\ell(\infty))_{\ell\in\Lambda\cup\{0\}}$.
\end{itemize}

\begin{figure}[h]
  \centering
  \includegraphics
                             {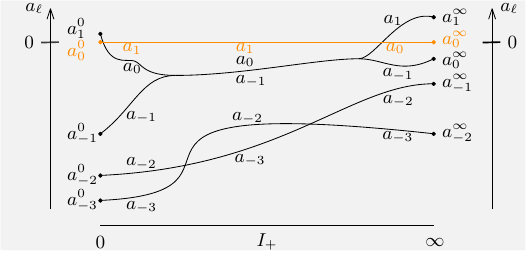}
  \caption{Spectral flow $\varsigma(A)=0$ along $\R_+$}
   \label{fig:fig-specflow-I_+-0}
\end{figure} 

\begin{definition}[Spectral flow - half-infinite forward interval]
\label{def:SF-hif}
The spectral flow of $A\in\Aa_{I_+}$ is defined as in Definition~\ref{def:SF-finite}
just by replacing $T$ by $\infty$.
\end{definition}

\subsubsection*{Half-infinite backward interval}

\begin{definition}[Spectral flow - half-infinite backward interval]
\label{def:SF-hib}
Given a backward path $A\in\Aa_{I_-}^*$, we define a forward path
$\tilde A\in\Aa_{I_+}^*$ by $\tilde A(s):=A(-s)$.
Then we define the spectral flow of the backward path
as the spectral flow of the negative forward path, in symbols
$\varsigma(A):=\varsigma(-\tilde A)$.
\end{definition}

\subsubsection*{Real line}

Similarly as in the case of half infinite intervals
the spectral flow can be defined along the whole real line.

Note that since the asymptotics are invertible,
no eigenvalues will cross zero any more
whenever $\abs{s}\ge T$ for some sufficiently large $T>0$.
Therefore, alternatively, one could also define the spectral flow
of $A\in \Aa_\R^*$ as the spectral flow of $A$ restricted to the
finite interval $[-T,T]$.

\newpage 
\section{Fredholm operators}\label{sec:semi-Fredholm}

Throughout $(H_0,H_1)$ is a Hilbert space pair.
Let $I\subset\R$ be connected,
then\footnote{
  In the notation of Def.\,2.10 and by Le.\,2.15
  in~\cite{neumeister:2020a}, cf.~\cite[Le.\,7.1]{Roubicek:2013a}, we have that
  \begin{equation}\label{eq:Oliver}
  \begin{split}
    P_1(\R_+;H_1,H_0)&=W^{1,(2,2)}(\R_+;H_1,H_0)
    \subset C^0(\R_+,H_0) ,
  \\
    P_1(\R_+;H_0^*,H_1^*)&=W^{1,(2,2)}(\R_+;H_0^*,H_1^*)
    \subset C^0(\R_+,H_1^*) .
  \end{split}
  \end{equation}
  }
\begin{equation}\label{eq:P_1-continuous}
\begin{split}
   P_1(I)=P_1(I;H_1,H_0)
   := L^2 (I,H_1)\cap W^{1,2}(I,H_0)
   &\subset C^0(I,H_0) ,
\\
   P_1^*(I)=P_1(I;H_0^*,H_1^*)
   := L^2 (I,H_0^*)\cap W^{1,2}(I,H_1^*)
   &\subset C^0(I,H_1^*) .
\end{split}
\end{equation}

\subsection{Real line}
\label{sec:Rabier}

\begin{definition}[Hessian path space $\Aa_\R^*$]
\label{def:Pp}
Let $\Aa_\R^*$ be the space of
continuous maps $A\colon(-\infty,\infty)\to\Ff(H_1,H_0)$
such that
both asymptotic limits exist and are invertible and symmetrizable, in
symbols
$$
   \A^\pm:=\lim_{s\to\pm\infty} A(s)
   \in \Ff^*(H_1,H_0) .
$$
\end{definition}

\boldmath
\subsubsection{Rabier's semi-Fredholm estimate for $D_A$}
\label{sec:Rabier-semi-Fred}
\unboldmath

The following theorem is due to Rabier~\cite{Rabier:2004a}.
In the case where $s\mapsto A(s)$ has a derivative 
the theorem is due to Robbin and Salamon~\cite[Lemma~3.9]{robbin:1995a}.
For the readers convenience we give a detailed explanation
of Rabier's ingenious argument of how to \emph{overcome the need of a derivative}.
In fact, Rabier proved his theorem even more general for some Banach
and not just Hilbert spaces which however requires additional arguments.

\begin{theorem}[Rabier]\label{thm:Rabier}
Given $A\in \Aa_\R^*$, there are constants $c, T>0$ such that
\[
   \Norm{\xi}_{P_1(\R)}
   \le c\left(\Norm{\xi}_{P_0([-T,T])}+\Norm{D_A\xi}_{P_0(\R)}\right)
\]
for every $\xi\in P_1(\R)$ where $D_A\colon P_1(\R)\to P_0(\R)$ is
defined by~(\ref{eq:D_A}) for $I=\R_+$.
\end{theorem}

\begin{remark}[Idea of proof]
One relates the operator $D_A$ associated to a path of Hessians
$s\mapsto A(s)$ to finitely many invertible operators
$D_{\A^{\lambda(\sigma_j)}}$ associated to constant-in-$s$ invertible paths
$\A^{\lambda(\sigma_j)}:=A(\sigma_j)-\lambda(\sigma_j)\iota$,
as illustrated by Figure~\ref{fig:fig-Rabier-idea}.
We indicate invertibility of Hessians by using the font $\A$.
Step~1: One shows that invertibility of a constant path $\A$ implies
invertibility of $D_\A$.
\begin{figure}[h]
  \centering
  \includegraphics
                             {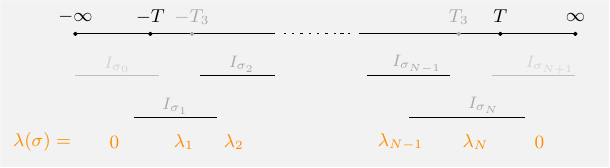}
  \caption{Approximate $D_A$ via finitely many invertible
                $D_{\A^{\lambda(\sigma_j)}}$,\,(\ref{eq:interv44}),\,(\ref{eq:jhu678-2})}
   \label{fig:fig-Rabier-idea}
\end{figure} 

\noindent
\textsc{Asymptotic ends.}
Step~2: The asymptotic limits $\A^\pm$ are invertible by assumption.
Hence so is, by continuity of the path $s\mapsto A(s)$,
each member of the path outside a sufficiently large compact
interval $[-T_2,T_2]$.
Step~3: One derives the estimate for $D_A$ outside a larger interval $[-T_3,T_3]$.

\noindent
\textsc{Compact center.}
Step~4: Since at each time $\sigma\in[-T_3,T_3]$ the operator
$A(\sigma)\colon H_1\to H_0$ is $H$-self-adjoint,
there exist non-eigenvalues $\mu_\sigma$, see Remark~\ref{rem:compact-op}.
Pick one, then the shifted operator $\A^{\mu_\sigma}:=A(\sigma)-\mu_\sigma\iota$
is invertible.
Step~5: But invertibility is an open property, thus
$\A^{\mu_\sigma}(\tau):=A(\tau)-\mu_\sigma\iota$ is still invertible
in a sufficiently narrow interval about $\sigma$, in symbols
$\forall \tau \in I_\sigma:=(\sigma-\eps_\sigma,\sigma+\eps_\sigma)$.
The compact interval $[-T_3,T_3]$ is covered by finitely many
intervals $I_{\sigma_1}, I_{\sigma_2},\dots, I_{\sigma_N}$,
see Figure~\ref{fig:fig-Rabier-idea}.
Define $\lambda(\sigma_j):=\mu_{\sigma_j}$ for $j=1,\dots,N$.
Step~6. If $A(\tau)$ is sufficiently close to $A(\sigma)$,
then one derives the desired estimate in a small neighborhood.

\noindent
\textsc{Putting things together.}
Step~7: One chooses $T>T_3$ and a suitable partition of unity for $\R$
to put the obtained estimates near $\pm \infty$ and such in the
compact center together. The closeness condition in Step~6
is achieved by subdividing $[-T,T]$ in sufficiently small
subintervals using that continuity of $s\mapsto A(s)$ along the
compact $[-T,T]$ is uniform.
\end{remark}

\begin{proof}
The proof is in seven steps.
Let $A\in \Aa_\R^*$, notation
$
   \A^\pm:=\lim_{s\to\pm\infty} A(s)
$.

\smallskip
\noindent
\textbf{Step~1} (Constant invertible path $\A$)\textbf{.}
Let $A(s)\equiv \A \in \Ff^*(H_1,H_0)$ be constant in time.
Then the following is true.
There exists a constant $C_1$ such that the following injectivity
estimate holds
\begin{equation}\label{eq:inject-Rabier}
   \Norm{\xi}_{P_1(\R)}
   \le C_1 \Norm{D_{\A }\xi}_{P_0(\R)}
\end{equation}
for every $\xi\in P_1(\R)$
and $D_{\A } \colon P_1(\R)\to P_0(\R)$ is an isomorphism with inverse
bounded by
\begin{equation}\label{eq:inject-Rabier-inv-1}
   \norm{(D_{\A })^{-1}}_{\Ll(P_0(\R),P_1(\R))}\le C_1 .
\end{equation}

\begin{proof}[1a -- Proof of the injectivity estimate~(\ref{eq:inject-Rabier})]
We can assume without loss of generality
that $\A\colon H_1\to H_0$ is symmetric. Indeed replacing the inner product
on $H_0$ by an equivalent one leads to equivalent norms
on $P_1(\R)$ and $P_0(\R)$ and therefore~(\ref{eq:inject-Rabier})
continues to hold after adapting the constant.

By definition~(\ref{eq:P_1-I}) of the space $P_1(\R)$ we get
\begin{equation}\label{eq:jhghgbhj77}
\begin{split}
   \Norm{\xi}_{P_1(\R)}^2
   &=\int_{-\infty}^\infty \Norm{\A^{-1} \A \xi(s)}_{H_1}^2+\Norm{\p_s\xi(s)}_{H_0}^2\, ds
\\
   &\le (1+\norm{\A^{-1}}_{\Ll(H_0,H_1)}^2)
   \left(\Norm{\A \xi}_{P_0(\R)}^2+\Norm{\p_s\xi}_{P_0(\R)}^2\right).
\end{split}
\end{equation}
On the other hand, by partial integration and symmetry of $\A $,
the mixed term is zero and we get
\begin{equation*}
\begin{split}
   \Norm{D_{\A }\xi}_{P_0(\R)}^2
   &=\int_{-\infty}^\infty\left(
   \Norm{\p_s\xi(s)}_{H_0}^2
   +2\INNER{\p_s\xi(s)}{\A \xi(s)}_0
   +\Norm{\A \xi(s)}_{H_0}^2
   \right) ds
\\
   &= \Norm{\A \xi}_{P_0(\R)}^2+\Norm{\p_s\xi}_{P_0(\R)}^2 .
\end{split}
\end{equation*}
This identity, together with~(\ref{eq:jhghgbhj77}), proves the
injectivity estimate in Step~1.
\end{proof}

\begin{proof}[1b -- Proof of an injectivity estimate for $D_{-\A^*}$]
There is a constant $C_1^*$ with
\begin{equation*}\label{eq:inject-Rabier-*}
   \Norm{\eta}_{P_1(\R;H_0^*,H_1^*)}
   \le C_1^* \Norm{D_{-\A^* }\eta}_{P_0(\R.H_1^*)}
\end{equation*}
for every $\eta\in P_1(\R;H_0^*,H_1^*)$.

\smallskip
\noindent
To see this note that the assumption $\A\in \Ff(H_1,H_0)$ implies
$\A^*\in \Ff(H_0^*,H_1^*)$, by Lemma~\ref{le:adjoint-map}.
Moreover, since $\A$ is invertible, the adjoint $\A^*$ is invertible as well.
Therefore Step 1b follows from Step~1a.
\end{proof}

\begin{proof}[1c -- Proof of surjectivity]
To prove surjectivity we first show that the image of $D_{\A }$ is closed
in $P_0(\R)$. For this we use the obtained injectivity estimate.
Suppose that $\eta_\nu$ is a sequence in the image of $D_{\A }$ which
converges to some $\eta\in P_0(\R)$.
We need to show that $\eta\in\im D_{\A }$.
Since $\eta_\nu\in\im D_{\A }$ there exists $\xi_\nu\in P_1(\R)$ such that
$\eta_\nu=D_{\A }\xi_\nu$. Since the sequence $\eta_\nu$ converges
it is a Cauchy sequence in $P_0(\R)$.
By the injectivity estimate the sequence $\xi_\nu$ is as well a Cauchy
sequence.
Since $P_1(\R)$ is complete the Cauchy sequence $\xi_\nu$ has a limit
$\xi\in P_1(\R)$. It follows that $\eta=D_{\A }\xi$ and therefore lies in
the image of $D_{\A }$. So $\im D_{\A }$ is closed.
\\
Hence to show that $D_{\A }$ is an isomorphism it suffices to check that
the orthogonal complement of $\im D_{\A }$ is trivial.
To see this pick $\eta\in(\im D_{\A })^\perp\subset P_0(\R)$.
This means that $\INNER{\eta}{D_{\A }\xi}_{P_0(\R)}=0$ for every $\xi\in P_1(\R)$,
hence
\begin{equation}\label{eq:Rabier-1c}
\begin{split}
   0
   &=\int_{-\infty}^\infty \INNER{\eta(s)}{\p_s\xi(s)+\A \xi(s)}_0 ds\\
   &=\int_{-\infty}^\infty \left({\color{red}\flat}\eta(s)\right) \p_s\xi(s)\, ds
   +\int_{-\infty}^\infty
   \underbrace{\left(\A^*{\color{red}\flat}\eta(s)\right)}_{\A^*\colon H_0^*\to H_1^*}
   \underbrace{\xi(s)}_{\in H_1} \, ds
\end{split}
\end{equation}
for every $\xi\in P_1(\R)$ and where ${\color{red}\flat}\colon H_{0}\to H_0^*$
is the insertion isometry~(\ref{eq:flat}).
Thus ${\color{red}\flat}\eta\in L^2(\R,H_0^*)$ has a weak derivative in $H_1^*$
satisfying $\p_s {\color{red}\flat}\eta=\A^*{\color{red}\flat}\eta$
where $\A^*\colon H_0^*\to H_1^*$ is the
adjoint of $\A \colon H_1\to H_0$.
In particular, ${\color{red}\flat}\eta$ lies in the kernel of the operator
$D_{-\A^*}\colon P_1(\R;H_0^*,H_1^*)\to P_0(\R;H_1^*)$.
But $D_{-\A^*}$ is injective by part 1b of the proof,
thus ${\color{red}\flat}\eta=0$, hence $\eta=0$.
This shows that $D_{\A } \colon P_1(\R;H_1,H_0)\to P_0(\R;H_0)$ is an
isomorphism.
Hence~(\ref{eq:inject-Rabier-inv-1}) follows from~(\ref{eq:inject-Rabier}).
This concludes the proof of Step~1.
\end{proof}

From now on we abbreviate
$$
   A_s:=A(s)\colon H_1\to H_0
   ,\qquad
   \text{$\A_s:=A(s)$ indicates invertibility.}
$$
We enumerate the constants by the step where they appear,
e.g. constant $T_2$ arises in Step~2.

\smallskip
\noindent
\textbf{Step~2} (Invertible near $\A^\pm$)\textbf{.}
There are constants $T_2,C_2>0$ such that 
for any fixed time $\sigma\in (-\infty,-T_2]\cup[T_2,\infty)$
the operators $\A_\sigma$ and $D_{\A_\sigma}$ are invertible~and
$$
   \Norm{(D_{\A_\sigma})^{-1}}_{\Ll(P_0(\R),P_1(\R))}\le C_2 .
$$

\begin{proof}
To prove Step~2 we show that the map
$$
   T\colon \Ff\to\Ll(P_1,P_0)
   ,\quad A\mapsto D_A
$$
is continuous.
To see this, given $A,\tilde A\in\Ff$, we calculate
\begin{equation}\label{eq:D(s)-cont}
\begin{split}
   \norm{D_A-D_{\tilde A}}_{\Ll(P_1,P_0)}:
   &=\sup_{\norm{\xi}_{P_1}=1}
   \norm{(A-\tilde A)\xi}_{P_0=L^2(\R,H_0)}
\\
   &=\sup_{\norm{\xi}_{P_1}=1}
   \left(
      \int_\R \norm{(A-\tilde A)\xi(s)}_{H_0}^2 \, ds
   \right)^{\frac12}
\\
   &\le \norm{A-\tilde A}_{\Ll(H_1,H_0)}
   \sup_{\norm{\xi}_{P_1}=1}
   \underbrace{
   \left(
      \int_\R \norm{\xi(s)}_{H_1}^2 \, ds
   \right)^{\frac12}
   }_{=\norm{\xi}_{L^2(\R,H_1)}\le \norm{\xi}_{P_1}=1}
\\
   &\le \norm{A-\tilde A}_{\Ll(H_1,H_0)} .
\end{split}
\end{equation}
Step~2 follows now from Step~1 (invertibility of asymptotic operators
$D_{\A^\pm}$) and with the help of Lemma~\ref{le:invertibility-open}
since the path $\R\ni s\mapsto A(s)$ converges at the ends to $D_{\A^\pm}$.
This proves Step~2.
\end{proof}

\smallskip
\noindent
\textbf{Step~3} (Asymptotic estimate)\textbf{.}
Let $T_2>0$ be the constant of Step~2.
There exists $T_3\ge T_2$ such that the following is true.
Suppose $\beta\in C^\infty(\R,\R)$ satisfies
\[
   \text{either $\supp\beta\subset(-\infty,-T_3)$}
   \qquad
   \text{or $\supp\beta\subset (T_3,\infty)$.}
\]
Then for every $\xi\in P_1(\R)$ we have the estimate
\[
   \Norm{\beta\xi}_{P_1(\R)}
   \le 2C_2\left(
   \Norm{\beta D_A\xi}_{P_0(\R)}+\Norm{\beta^\prime\xi}_{P_0(\R)}
   \right)
\]
where $C_2>0$ is the constant from Step~2.

\begin{proof}
Since $\lim_{s\to\pm\infty} A(s)=\A^\pm$
there exists $T_3\ge T_2$ such that
\begin{equation*}
\begin{split}
   \Norm{\A_\sigma-\A^+}_{\Ll(H_1,H_0)}\le\tfrac{1}{4C_2}\quad
   &\forall \sigma\ge T_3
\\
   \Norm{\A_\sigma-\A^-}_{\Ll(H_1,H_0)}\le\tfrac{1}{4C_2}\quad
   &\forall \sigma\le -T_3 
\end{split}
\end{equation*}
where Step~2 provides the constant $C_2>0$ and invertibility of
$\A_\sigma:=A(\sigma)$.

In the following we only discuss the case where
$\supp \beta\subset (T_3,\infty)$, the other case is analogous.
Assume that $\sigma\in \supp\beta\subset (T_3,\infty)$.
We calculate
\begin{equation}\label{eq:jhu678}
\begin{split}
   D_{\A_\sigma}\beta\xi
   &=\p_s(\beta\xi)+\A_\sigma \beta\xi\\
   &=\beta^\prime\xi+\beta\left(\p_s\xi+(\A_\sigma-A+A)\xi\right)\\
   &=\beta^\prime\xi+\beta D_A\xi+(\A_\sigma-A)\beta\xi .
\end{split}
\end{equation}
By Step~2 the operator $D_{\A_\sigma}$ is invertible,
we multiply with $(D_{\A_\sigma})^{-1}$ to get
\begin{equation*}
\begin{split}
   \beta\xi
   &=(D_{\A_\sigma})^{-1}\bigl(\beta^\prime\xi+\beta D_A\xi+(\A_\sigma-A)\beta\xi\bigr) .
\end{split}
\end{equation*}
Taking norms we estimate
\begin{equation*}
\begin{split}
   \Norm{\beta\xi}_{P_1(\R)}
   &\le\Norm{ (D_{\A_\sigma})^{-1}}_{\Ll(P_0(\R),P_1(\R))}
   \Norm{\beta^\prime\xi+\beta D_A\xi+(\A_\sigma-A)\beta\xi}_{P_0(\R)}
   \\
   &\le C_2
   \left(
   \Norm{\beta D_A\xi}_{P_0(\R)}
   +\Norm{\beta^\prime\xi}_{P_0(\R)}
   +\Norm{(\A_\sigma-A)\beta\xi}_{P_0(\R)}
   \right) .
\end{split}
\end{equation*}
It remains to estimate the difference
\begin{equation*}
\begin{split}
   &\Norm{(\A_\sigma-A)\beta\xi}_{P_0(\R)}^2\\
   &=\int_{-\infty}^\infty \Norm{\left(\A_\sigma-A(s)\right)\beta(s)\xi(s)}_{H_0}^2\, ds
\\
   &=\int_{T_3}^\infty
   \Norm{\A_\sigma-\A^++\A^+-A(s)}_{\Ll(H_1,H_0)}^2
   \Norm{\beta(s)\xi(s)}_{H_1}^2 ds
\\
   &\le\int_{T_3}^\infty
   2\left(
   \Norm{\A_\sigma-\A^+}_{\Ll(H_1,H_0)}^2
   +\Norm{\A^+-A(s)}_{\Ll(H_1,H_0)}^2
   \right)
   \Norm{\beta(s)\xi(s)}_{H_1}^2 ds\\
   &\le\frac{1}{4C_2^2}\int_{T_3}^\infty \Norm{\beta(s)\xi(s)}_{H_1}^2 ds\\
   &\le\frac{1}{4C_2^2}\Norm{\beta\xi}_{P_1(\R)}^2 .
\end{split}
\end{equation*}
The last two estimates together show that
\begin{equation*}
\begin{split}
   \Norm{\beta\xi}_{P_1(\R)}
   &\le 2C_2
   \left(
   \Norm{\beta D_A\xi}_{P_0(\R)}
   +\Norm{\beta^\prime\xi}_{P_0(\R)}
   \right) .
\end{split}
\end{equation*}
This proves Step~3.
\end{proof}

\smallskip
\noindent
\textbf{Step~4} (Invertibility perturbation)\textbf{.}
For any $\sigma\in\R$ there is $\mu_\sigma\in\R$ such that
the \textbf{shifted operator}
$$
   \A^{\mu_\sigma}:=A_\sigma-\mu_\sigma \iota\colon H_1\to H_0
$$
is invertible where $ \iota\colon H_1\INTO H_0$ is inclusion.

\begin{proof}
Since $A_\sigma$ is symmetric
and inclusion $\iota\colon H_1\INTO H_0$ is compact,
the spectrum of $A_\sigma$, as unbounded operator on $H_0$
with dense domain $H_1$, is a discrete unbounded subset of $\R$
with no (finite) accumulation point, see
Remarks~\ref{rem:Fredholm-index-zero} and~\ref{rem:compact-op}.
Pick $\mu_\sigma$ in the complement of the spectrum of $A_\sigma$,
that is $\mu_\sigma\in\Rr(A)$.
\end{proof}

\smallskip
\noindent
\textbf{Step~5} (Invertibilizing $A$ along $[-T_3,T_3]$ by
finitely many shifts $\lambda_1,\dots,\lambda_N$)\textbf{.}
Let $T_3>0$ be from Step~3.
There is a finite set 
$\Lambda^\prime=\{\lambda_1,\dots,\lambda_N\}\subset \R$ and a
constant $C_5>0$ such the following~holds. Fix any $\sigma\in[-T_3,T_3]$.
Then there exists an element 
$\lambda(\sigma)\in \{\lambda_1,\dots,\lambda_N\}$, such that the
operator $D_{\A^{\lambda(\sigma)}}=\p_s+\A^{\lambda(\sigma)}$ is invertible
and there is the estimate
$$
   \norm{(D_{\A^{\lambda(\sigma)}})^{-1}}
      _{\Ll(P_0(\R),P_1(\R))}
   \le C_5 .
$$

\begin{proof}
Invertibility is an open property.
Given any $\sigma\in\R$ there exists, by Step~4,
a real number $\mu_\sigma$, pick $\mu_\sigma\in\Rr(A_\sigma)$, such that
$\A^{\mu_\sigma}=A_\sigma-\mu_\sigma\iota$ is invertible.
Hence $D_{\A^{\mu_\sigma}}$ is invertible by Step~1.
For $\tau$ near $\sigma$ we vary $\A^{\mu_\sigma}$ in the form
$$
   \A^{\mu_\sigma}(\tau)
    :=A_\tau-\mu_\sigma\iota
   ,\qquad
   \A^{\mu_\sigma}(\sigma)=\A^{\mu_\sigma} .
$$
By continuity~(\ref{eq:D(s)-cont}) of $s\mapsto D_{A(s)}$ there exists,
by Lemma~\ref{le:invertibility-open},
a constant $\eps_\sigma>0$ with the following significance.
At any time $\tau\in I_\sigma:=(\sigma-\eps_\sigma,\sigma+\eps_\sigma)$
the operator $D_{\A^{\mu_\sigma}(\tau)}$
is invertible and the inverse is bounded by
\begin{equation}\label{eq:interv33}
       \norm{\left(D_{\A^{\mu_\sigma}(\tau)}\right)^{-1}}_{\Ll(P_0,P_1)}
       \le 2
       \norm{\left(D_{\A^{\mu_\sigma}}\right)^{-1}}_{\Ll(P_0,P_1)} .
\end{equation}
Since $[-T_3,T_3]$ is compact there exist $N\in\N$
and $\sigma_1,\dots,\sigma_N\in[-T_3,T_3]$ with
\begin{equation}\label{eq:interv44}
   [-T_3,T_3]\subset\bigcup_{i=1}^N
   I_{\sigma_i}
   ,\qquad
   I_{\sigma_i}:=(\sigma_i-\eps_{\sigma_i}, \sigma_i+\eps_{\sigma_i}) .
\end{equation}
Now define
$$
   \Lambda^\prime:=\{\lambda_i:=\mu_{\sigma_i}\mid i=1,\dots,N\}
   ,\quad
   C_5
   :=2
   \max_{i=1,\dots,N} \norm{\left(D_{\A^{\mu_{\sigma_i}}}\right)^{-1}}_{\Ll(P_0,P_1)}
   .
$$
Suppose now that $\sigma\in[-T_3,T_3]$, then by
the finite covering property~(\ref{eq:interv44})
there exists $i\in\{1,\dots,N\}$ such that $\sigma\in I_{\sigma_i}$.
%
We choose such $i$ and define
$$
   \lambda(\sigma):=\mu_{\sigma_i} .
$$
For this choice $\A^{\lambda(\sigma)}=A_\sigma-\mu_{\sigma_i}\iota
=:\A^{\mu_{\sigma_i}}(\sigma)$ and there is the estimate
$$
   \norm{(D_{\A^{\lambda(\sigma)}})^{-1}}_{\Ll(P_0,P_1)}
   =\norm{(D_{\A^{\mu_{\sigma_i}}(\sigma)})^{-1}}_{\Ll}
   \stackrel{\text{(\ref{eq:interv33})}}{\le} 2
   \norm{\left(D_{\A^{\mu_{\sigma_i}}}\right)^{-1}}_{\Ll}
   \le C_5
$$
where $\norm{\cdot}_\Ll=\norm{\cdot}_{\Ll(P_0,P_1)}$.
This proves Step~5.
\end{proof}

\smallskip
\noindent
\textbf{Step~6} (Small intervals)\textbf{.}
Let $C_6:=\max\{C_2,C_5\}$.
Let $\lambda^*:=\max\abs{\Lambda^\prime}$
be the maximal absolute value of the elements
of the finite set $\Lambda^\prime\subset\R$ in Step~5.
Then for every $\beta\in C^\infty(\R,\R)$ with the property
$$
   \sup_{\sigma,\tau\in\supp\beta}\norm{A_\sigma-A_\tau}_{\Ll(H_1,H_0)}\le\frac{1}{2C_6}
$$
it holds that
\[
   \Norm{\beta\xi}_{P_1(\R)}
   \le 2C_6\left(
   \Norm{\beta D_A\xi}_{P_0(\R)}+\Norm{\beta^\prime\xi}_{P_0(\R)}
   +\lambda^*\Norm{\beta\xi}_{P_0(\R)}
   \right)
\]
for every $\xi\in P_1(\R)$.

\begin{proof}
By Step~2 and Step~5 there exists a map
$\lambda\colon\R\to\R$ satisfying 
$\lambda(\sigma)\in\Lambda^\prime$ if $\sigma\in[-T_3,T_3]$ and
$\lambda(\sigma)=0$ if $\sigma\in(-\infty,-T_3)\cup(T_3,\infty)$ and
such that
$$
   \norm{(D_{\A^{\lambda(\sigma)}})^{-1}}_{\Ll(P_0(\R),P_1(\R))}
   \le C_5 \le C_6
   ,\qquad
   \A^{\lambda(\sigma)}
   :=A_\sigma-\lambda(\sigma)\iota .
$$
Now the proof of Step~6 proceeds similarly as the proof of Step~3.
Suppose that $\sigma$ lies in the support of $\beta$.
Computing as in~(\ref{eq:jhu678}) we get
\begin{equation}\label{eq:jhu678-2}
\begin{split}
   D_{\A^{\lambda(\sigma)}}\beta\xi
   &=\beta^\prime\xi+\beta D_A\xi+(A_\sigma-A)\beta\xi -\lambda(\sigma)\beta\xi.
\end{split}
\end{equation}
By construction of $\lambda$ the operator
$D_{\A^{\lambda(\sigma)}}$ is invertible so that we can write
\begin{equation*}
\begin{split}
   \beta\xi
   &=(D_{\A^{\lambda(\sigma)}})^{-1}
   \bigl(\beta^\prime\xi+\beta D_A\xi+(A_\sigma-A)\beta\xi
   -\lambda(\sigma)\beta\xi\bigr) .
\end{split}
\end{equation*}
Taking norms we estimate
\begin{equation*}
\begin{split}
   &\Norm{\beta\xi}_{P_1(\R)}\\
   &\le C_6
   \left(
   \Norm{\beta^\prime\xi}_{P_0(\R)}
   +\Norm{\beta D_A\xi}_{P_0(\R)}
   +\Norm{(A_\sigma-A)\beta\xi}_{P_0(\R)}
   +\Abs{\lambda(\sigma)}\Norm{\beta\xi}_{P_0(\R)}
   \right) .
\end{split}
\end{equation*}
Now we use the hypothesis
$\sup_{\sigma,\tau\in\supp\beta}\Norm{A_\sigma-A_\tau}_{\Ll(H_1,H_0)}\le\frac{1}{2C_6}$
to estimate
\begin{equation*}
\begin{split}
   \Norm{(A_\sigma-A)\beta\xi}_{P_0(\R)}^2
   &=\int_{-\infty}^\infty \Norm{\left(A_\sigma-A(s)\right)\beta(s)\xi(s)}_{H_0}^2\, ds
\\
   &=\int_{\supp\beta}
   \Norm{A_\sigma-A(s)}_{\Ll(H_1,H_0)}^2
   \Norm{\beta(s)\xi(s)}_{H_1}^2 ds
\\
   &\le\frac{1}{4C_6^2}\int_{\supp\beta} \Norm{\beta(s)\xi(s)}_{H_1}^2 ds\\
   &\le\frac{1}{4C_6^2}\Norm{\beta\xi}_{P_1(\R)}^2 .
\end{split}
\end{equation*}
The last two estimates together imply Step~6.
\end{proof}

\smallskip
\noindent
\textbf{Step~7} (Partition of unity)\textbf{.}
We prove Theorem~\ref{thm:Rabier}.

\begin{proof}
Choose $T>T_3$
and a finite partition of unity $\{\beta_j\}_{j=0}^{M+1}$ for $\R$,
where each $\beta_j\colon [0,1]\to\R$ is smooth, with the properties
$$
   \supp\beta_0\subset (-\infty,-T_3)
   ,\qquad
   \supp\beta_{M+1}\subset (T_3,\infty) ,
$$
and for $j=1,\dots,M$ it holds
$$
   \sup_{\sigma,\tau\in\supp\beta_j}\norm{A_\sigma-A_\tau}_{\Ll(H_1,H_0)}
   \le\frac{1}{2C_6}
   ,\qquad
  \supp\beta_j\subset (-T,T) .
$$
That such a partition exists follows from the continuity
of $s\mapsto A(s)$ and the fact that on the compact set $[-T_3,T_3]$
continuity becomes uniform continuity.
Let $\xi\in P_1(\R)$. Then by Step~3 we have
\begin{equation*}
\begin{split}
   \Norm{\beta_0\xi}_{P_1(\R)}
   &\le 2C_6\left(
   \Norm{\beta_0 D_A\xi}_{P_0(\R)}+\Norm{\beta_0^\prime\xi}_{P_0(\R)}
   \right) ,
\\
   \Norm{\beta_{M+1}\xi}_{P_1(\R)}
   &\le 2C_6\left(
   \Norm{\beta_{M+1} D_A\xi}_{P_0(\R)}+\Norm{\beta_{M+1}^\prime\xi}_{P_0(\R)}
   \right) .
\end{split}
\end{equation*}
By Step~6 we have for each $j=1,\dots,M$ an estimate
\[
   \Norm{\beta_j\xi}_{P_1(\R)}
   \le 2C_6\left(
   \Norm{\beta_j D_A\xi}_{P_0(\R)}+\Norm{\beta_j^\prime\xi}_{P_0(\R)}
   +\lambda^*\Norm{\beta_j\xi}_{P_0(\R)}
   \right) .
\]
Let $B:=\max\{\norm{\beta_0^\prime}_\infty,
\norm{\beta_1^\prime}_\infty,\dots, \norm{\beta_{M+1}^\prime}_\infty\}$.
Put the estimates together to~get
\begin{equation*}
\begin{split}
   \Norm{\xi}_{P_1(\R)}
   &\le\sum_{j=0}^{M+1} \Norm{\beta_j\xi}_{P_1(\R)}\\
   &\le 2C_6 \sum_{j=0}^{M+1}\left(
   \Norm{\beta_j D_A\xi}_{P_0(\R)}+\Norm{\beta_j^\prime\xi}_{P_0([-T,T])}
  \right)
   +2C_6 \lambda^* \sum_{j=1}^{M}\Norm{\beta_j\xi}_{P_0([-T,T])}\\
   &\le 2C_6 (M+2) \Norm{D_A\xi}_{P_0(\R)}
   +2C_6 \left(B (M+2)+\lambda^* M\right) \Norm{\xi}_{P_0([-T,T])}
\end{split}
\end{equation*}
where in the second inequality we replaced the $P_0(\R)$ norm by the
$P_0([-T,T])$ norm due to the supports of the $\beta_j$'s and their
derivatives.\footnote{
   $\beta_0\equiv 1$ on $(-\infty,-T]$
   and
   $\beta_{M+1}\equiv 1$ on $[T,\infty)$, so the derivatives vanish.
   }
Setting
$$
   c:=\max\{2C_6 (M+2), 2C_6 \left(B (M+2)+\lambda^* M\right)\}
$$
proves Step~7.
\end{proof}
\noindent
The proof of Theorem~\ref{thm:Rabier} is complete.
\end{proof}

\boldmath
\subsubsection[Semi-Fredholm estimate for the adjoint $D_A^*:=D_{-A^*}$]{Semi-Fredholm estimate for the adjoint $D_A^*$}
\label{sec:Rabier-D-A*}
\unboldmath

Let $A\in  \Aa_\R^*$.
We call the following operator the
\textbf{adjoint of \boldmath$D_A$}, namely
$$
   D_A^*:=D_{-A^*}\colon P_1(\R;H_0^*,H_1^*)\to P_0(\R;H_1^*)
   ,\quad
   \xi\mapsto \p_s\xi-A(s)^*\xi .
$$

\begin{corollary}\label{cor:Rabier-*}
For $A\in  \Aa_\R^*$ there are constants $c$ and $T>0$ with
\[
   \Norm{\xi}_{P_1(\R;H_0^*,H_1^*)}
   \le c\left(\Norm{\xi}_{P_0([-T,T];H_1^*)}+\Norm{D_A^*\xi}_{P_0(\R;H_1^*)}\right)
\]
for every $\xi\in P_1(\R;H_0^*,H_1^*)$.
\end{corollary}

\begin{proof}
Theorem~\ref{thm:Rabier} and Lemma~\ref{le:adjoint-map};
see also Remark~\ref{rem:symmetrizable}.
\end{proof}

\boldmath
\subsubsection{Fredholm property of $D_A$}
\unboldmath

An immediate consequence of  Theorem~\ref{thm:Rabier} is that the
operator $D_A\colon P_1(\R)\to P_0(\R)$ is \textbf{semi-Fredholm},
i.e. the kernel of $D_A$ is of finite dimension and the range is closed.
Indeed the restriction map $P_1(\R)\to P_0([-T,T])$ is compact,
see e.g.~\cite[Lemma~3.8]{robbin:1995a}.
Hence the semi-Fredholm property follows
from~\cite[Lemma~A.1.1]{mcduff:2004a}.

To see that $D_A$ is actually a Fredholm operator we need to examine
its cokernel. For that purpose let
$\eta\in\left(\im D_A\right)^\perp\subset P_0(\R)=L^2(\R,H_0)$, that is
\[
   \INNER{\eta}{D_A\xi}_{P_0(\R)}
   =0
   ,\quad
   \forall \xi\in P_1(\R)=L^2 (\R,H_1)\cap W^{1,2}(\R,H_0) .
\]
To put
it differently
\begin{equation}\label{eq:5ru6tyubhj}
\begin{split}
   0
   &=\int_{-\infty}^\infty\INNER{\eta(s)}{\p_s\xi(s)+A(s)\xi(s)}_0 \, ds\\
   &=\int_{-\infty}^\infty\left({\color{red}\flat}\eta(s)\right) \p_s\xi(s) \, ds
   +\int_{-\infty}^\infty
   \underbrace{(A(s)^*{\color{red}\flat}\eta(s))}_{A(s)^*\colon H_0^*\to H_1^*}\xi(s) \, ds
\end{split}
\end{equation}
for every $\xi\in P_1(\R)$ and where ${\color{red}\flat}\colon H_{0}\to H_0^*$
is the insertion isometry~(\ref{eq:flat}).
Interpreting the $\xi$'s as test functions
this means that ${\color{red}\flat}\eta\in L^2(\R,H_0^*)$ has a weak
derivative in $H_1^*$, namely $\p_s {\color{red}\flat}\eta=A^*{\color{red}\flat}\eta$.
Hence ${\color{red}\flat}\eta$ lies in
$L^2(\R,H_0^*)\cap W^{1,2}(\R,H_1^*)=P_1(\R;H_0^*,H_1^*)$ and satisfies
$D_A^*{\color{red}\flat}\eta=\p_s {\color{red}\flat}\eta-A(s)^*{\color{red}\flat}\eta=0$.
Observe that $D_A^*$ is a map
\begin{equation}\label{eq:D_-A*}
   D_A^*=\p_s-A(s)^*\colon 
   \underbrace{L^2(\R,H_0^*) \cap W^{1,2}(\R,H_1^*)}
      _{P_1(\R;H_0^*,H_1^*)\stackrel{\text{(\ref{eq:P_1-continuous})}}{\subset} C^0(\R,H_1^*)}
   \to \underbrace{L^2(\R,H_1^*)}_{P_0(\R;H_1^*)} .
\end{equation}
We proved ${\color{red}\flat}\left(\im D_A\right)^\perp\subset \ker D_A^*$.
Vice versa, fix ${\color{red}\flat}\eta\in \ker D_A^*$ Then
$(D_A^*{\color{red}\flat}\eta)(s)=0_{H_1^*}$ for every $s\in\R$. Pick $\xi\in
P_1(\R)$ and integrate $(D_A^*{\color{red}\flat}\eta)(s)\xi(s)=0$ over $s\in\R$
to get back to $\INNER{\eta}{D_A\xi}_{P_0(\R)}=0$. We proved

\begin{lemma}\label{le:im-perp}
Given $A\in \Aa_\R^*$,
consider $D_A\colon P_1(\R;H_1,H_0)\to P_0(\R;H_0)$
and $D_A^*\colon P_1(\R;H_0^*,H_1^*)\to P_0(\R;H_1^*)$,
then
$$
   \left(\im D_A\right)^\perp
   \stackrel{{\color{red}\flat}}{\simeq}\ker D_A^*
   \subset L^2(\R,H_0^*)\cap W^{1,2}(\R,H_1^*).
$$
\end{lemma}

\begin{corollary}\label{cor:Fredholm-D_A}
$D_A=\mathfrak{D}_A\colon P_1(\R;H_1,H_0)\to P_0(\R;H_0)$ is Fredholm
$\forall A\in \Aa_\R^*$.
\end{corollary}

\begin{proof}
By Theorem~\ref{thm:Rabier} the operator $D_A$ is semi-Fredholm
(finite dimensional kernel and closed image). Since $D_A$ has closed
image, it follows that $\coker D_A=\left(\im D_A\right)^\perp$, 
but $\left(\im D_A\right)^\perp\simeq\ker D_A^*$ by 
Lemma~\ref{le:im-perp}.
By Corollary~\ref{cor:Rabier-*} the operator~(\ref{eq:D_-A*})
is semi-Fredholm as well.
This proves Corollary~\ref{cor:Fredholm-D_A}.
\end{proof}

\newpage
\subsection{Finite interval}\label{sec:finite-interval}

Pick a Hessian path $A\in \Aa_{I_T}^*
:=\{A\in \Aa_{I_T}\mid \text{$A(-T)$ and $A(T)$ are invertible} \}$
along the finite interval $I_T=[-T,T]$.
Then $A\colon [-T,T]\to \Ff=\Ff(H_1,H_0)$
takes values in the symmetrizable Fredholm operators of index zero;
cf. Remarks~\ref{rem:symmetrizable} and~\ref{rem:compact-op}.
In order to eventually get to Fredholm operators,
it is not enough that the Hessians at the interval ends are invertible, notation
$$
   \A_{-T}:=A(-T)
   ,\qquad
   \A_T:=A(T) .
$$
In addition, one must impose boundary conditions
formulated in terms of the spectral projections
$\pi_+^{\A_{-T}}$ sitting at time $-T$ and
$\pi_-^{\A_T}$ at time $T$; see~(\ref{eq:H_1/2}).

\boldmath
\subsubsection{Estimate for $D_A$}\label{sec:I_T-D_A}
\unboldmath

In~this~section we study the linear operator  $\p_s+A$ as a map
\begin{equation*} 
   D_A\colon P_1(I_T)\to P_0(I_T)
   ,\quad
   \xi\mapsto \p_s\xi+A(s)\xi
\end{equation*}
and the augmented operator
$\mathfrak{D}_{A}\colon P_1(I_T)\to \Ww(I_T;\A_{- T},\A_T)$
in~(\ref{eq:augmented-T}).
Define Hilbert spaces $P_0(I_T )=P_0(I_T;H_0)$ and
$P_1(I_T )=P_1(I_T;H_1,H_0)$ by~(\ref{eq:P_1-I}).
These operators are \emph{not} Fredholm: although $D_A$ has closed
image and finite dimensional co-kernel, the kernel is infinite
dimensional in the Floer and Morse case; see Figure~\ref{fig:D_A-finite-interval}.

\begin{theorem}\label{thm:finite-interval}
For $A\in\Aa_{I_T}^*$ there exists a constant $c>0$ such that
\[
   \norm{\xi}_{P_1(I_T)}
   \le c\left(\norm{\xi}_{P_0(I_T)}+\norm{D_A\xi}_{P_0(I_T)}
   +\norm{\pi_+^{\A_{-T}}\xi(-T)}_{\frac12}+\norm{\pi_-^{\A_T}\xi(T)}_{\frac12}\right)
\]
for every $\xi\in P_1(I_T)$.
\hfill 
{\color{gray}\small we abbreviate
$\norm{\cdot}_{\frac12}:=\norm{\cdot}_{H_{\frac12}}$
}
\end{theorem}

Orthogonal projections are bounded by $1$,
hence the theorem reduces to

\begin{corollary}\label{cor:finite-interval}
For $A\in\Aa_{I_T}^*$ there exists a constant $c>0$ such that
\[
   \Norm{\xi}_{P_1(I_T)}
   \le c\left(\Norm{\xi}_{P_0(I_T)}+\Norm{D_A\xi}_{P_0(I_T)}
   +\Norm{\xi(-T)}_{\frac12}+\Norm{\xi(T)}_{\frac12}\right)
\]
for every $\xi\in P_1(I_T)$.
\end{corollary}

\begin{proof}[Proof of Theorem~\ref{thm:finite-interval}]
We prove the theorem in five steps.
It is often convenient to abbreviate $\xi_s:=\xi(s)$ and $A_s:=A(s)$.
We enumerate the constants by the step where they appear,
e.g. constant $C_1$ arises in Step~1.

\smallskip
\noindent
\textbf{Step~1} (Constant invertible case)\textbf{.}
Let $A(s)\equiv\A {\color{gray}\; =\A_T=\A_{-T}}$
be constant in time and invertible.
Then there is a constant $C_1>0$ such that
\begin{equation}\label{eq:Step-1-Sim}
   \Norm{\xi}_{P_1(I_T)}
   \le C_1 \left(\Norm{D_{\A}\xi}_{P_0(I_T)}
   +\norm{\pi_+^\A\xi_{-T}}_{\frac12}
   +\norm{\pi_-^\A\xi_{T}}_{\frac12}
   \right)
\end{equation}
for every $\xi\in P_1(I_T)$. 
Moreover, for the constant path $\A$ the augmented operator
\begin{equation}\label{eq:augmented-T-AA}
\begin{split}
   \mathfrak{D}_\A\colon P_1(I_T)
   &\to P_0(I_T)\times H^+_{\frac12}(\A) \times H^-_{\frac12}(\A)
   =: \Ww(I_T;\A,\A)
   \\
   \xi
   &\mapsto
   \left(D_A\xi,\pi_+^{\A}\xi_{-T},\pi_-^{\A}\xi_T\right)
\end{split}
\end{equation}
is bijective.

\begin{proof}
Step~1 was proved in~\cite[Thm.\,3.1.6~i)]{Simcevic:2014a}.
We follow her proof.
By changing the constant $C_1$ if necessary, we can assume without
loss of generality, as explained in Section~\ref{sec:int-spaces-1/2}, that
\begin{equation}\label{eq:hyp-I_T}
   \text{$\A\colon H_1\to H_0$ is a symmetric isometry.}
\end{equation}
In this case~(\ref{eq:Step-1-Sim}) actually holds with unit constant $C_1=1$.
We abbreviate $\pi_\pm:=\pi_\pm^\A$.
Since $D_{\A}\xi=\p_s\xi+\A\xi$,
integration by parts yields
\begin{equation}\label{eq:gh654-NEW}
\begin{split}
   &\Norm{D_{\A}\xi}_{P_0(I_T)}^2\\
   &=\int_{-T}^T\left(
   \Norm{\p_s\xi_s}_{H_0}^2
   +2\INNER{\p_s\xi_s}{\A\xi_s}_0
   +\Norm{\A\xi_s}_{H_0}^2
   \right) ds
\\
   &= \Norm{\A\xi}_{P_0(I_T)}^2+\Norm{\p_s\xi}_{P_0(I_T)}^2
   +\INNER{\xi_T}{\A\xi_T}_0-\INNER{\xi_{-T}}{\A\xi_{-T}}_0 
\end{split}
\end{equation}
for every $\xi\in P_1(I_T)$.
To see this note that
\begin{equation*}
\begin{split}
   \int_{-T}^T \INNER{\p_s\xi_s}{\A\xi_s}_0 \, ds
   &=\int_{-T}^T \left(
   \p_s\INNER{\xi_s}{\A\xi_s}_0-\INNER{\xi_s}{\A\p_s\xi_s}_0
   \right) ds \\
   &=\int_{-T}^T \left(
   \p_s\INNER{\xi_s}{\A\xi_s}_0-\INNER{\A\xi_s}{\p_s\xi_s}_0
   \right) ds
\end{split}
\end{equation*}
where in the last step we used $H_0$-symmetry of $\A$.
Thus
\begin{equation*}
\begin{split}
   2\int_{-T}^T \INNER{\p_s\xi_s}{\A\xi_s}_0 \, ds
   &=\int_{-T}^T \p_s\INNER{\xi_s}{\A\xi_s}_0 ds \\
   &=\INNER{\xi_T}{\A\xi_T}_0-\INNER{\xi_{-T}}{\A\xi_{-T}}_0 .
\end{split}
\end{equation*}
This proves~(\ref{eq:gh654-NEW}). The opposite signs
will be crucial soon.
As $\A$ is an isometry we get
\begin{equation*}
\begin{split}
   -\INNER{\pi_-\xi_T}{\A \pi_-\xi_T}_0
   &\stackrel{\text{(\ref{eq:1/2-length-v_ell})}}{=}
   \norm{\pi_-\xi_T}_{\frac12}^2
\\
   \INNER{\pi_+\xi_{-T}}{\A \pi_+\xi_{-T}}_0
   &\stackrel{\text{(\ref{eq:1/2-length-v_ell})}}{=}
   \norm{\pi_+\xi_{-T}}_{\frac12}^2 .
\end{split}
\end{equation*}
So we get
{\small
\begin{equation}\label{eq:vghvhj95478}
\begin{split}
   \INNER{\xi_{T}}{\A\xi_{T}}_0
   =\overbrace{\INNER{\pi_+\xi_{T}}{\A \pi_+\xi_{T}}_0}^{\ge 0}
   +\INNER{\pi_-\xi_{T}}{\A \pi_-\xi_{T}}_0
   &\ge-\norm{\pi_-\xi_{T}}_{\frac12}^2
\\
   \INNER{\xi_{-T}}{\A\xi_{-T}}_0
   =\INNER{\pi_+\xi_{-T}}{\A \pi_+\xi_{-T}}_0
   +\underbrace{\INNER{\pi_-\xi_{-T}}{\A \pi_-\xi_{-T}}_0}_{\le 0}
   &\le \norm{\pi_+\xi_{-T}}_{\frac12}^2
\end{split}
\end{equation}
}where the identities use that the mixed terms vanish by
$H_0$-orthogonality; see Case~1 in Section~\ref{sec:int-spaces-1/2}.
Since $\A$ is an isometry we get identity one in the following
\begin{equation*}
\begin{split}
   \Norm{\xi}_{P_1(I_T)}^2
   &\stackrel{\text{\,(\ref{eq:P_1-inner})\,}}{=}
   \Norm{\A\xi}_{P_0(I_T)}^2+\Norm{\p_s\xi}_{P_0(I_T)}^2
   \\
   &\stackrel{\text{(\ref{eq:gh654-NEW})}}{=}
   \Norm{D_{\A}\xi}_{P_0(I_T)}^2
   +\INNER{\xi_{-T}}{\A\xi_{-T}}_0
   -\INNER{\xi_T}{\A\xi_T}_0
   \\
   &\stackrel{{\color{white}\text{(4.16)}}}{\le}
   \Norm{D_{\A}\xi}_{P_0(I_T)}^2
   +\Norm{\pi_+\xi_{-T}}_{H_{\frac12}}^2
   +\Norm{\pi_-\xi_T}_{H_{\frac12}}^2 .
\end{split}
\end{equation*}
The last inequality is by the previous displayed estimate.
This proves~(\ref{eq:Step-1-Sim}).
In particular, this implies that the augmented operator
$\mathfrak{D}_{\A}$ is injective.

\medskip\noindent
\textsc{Claim 1.} $\mathfrak{D}_{\A}$ is surjective.

\smallskip\noindent
To see this consider the orthonormal basis
$\Vv(\A)=(v_\ell)_{\ell\in \Lambda}\subset H_1$
of $H_0$ from~(\ref{eq:ONB-A}) enumerated by the
ordering~(\ref{eq:A^0-spectrum}) of the eigenvalues $a_\ell$ of $\A$.
Pick $\zeta=(\eta,x,y)\in \Ww(I_T;\A,\A)$.
Given $s\in[-T,T]$, with respect to the common basis $\Vv(\A)$
of $H_0$ and $H^-_\frac12(\A)\oplus H^+_\frac12(\A)$ we write
$$
   \eta(s)=\sum_{\ell\in\Lambda} \eta_\ell(s) v_\ell \in H_0
   ,\qquad
   x=\sum_{\nu\in\Lambda_+} x_\nu v_\nu
   ,\quad
   y=\sum_{\nu\in\Lambda_-} y_{-\nu} v_{-\nu} .
$$
We are looking for a map $\xi\in P_1(I_T)$ of the form
$s\mapsto \xi(s)=\sum_{\ell\in\Lambda} \xi_\ell(s) v_\ell$
which, for each $\ell\in\Lambda$, satisfies a linear inhomogeneous ODE of
the form
\begin{equation}\label{eq:inhom-ODE}
   \p_s\xi_\ell(s)+a_\ell\xi_\ell(s)=\eta_\ell(s)
\end{equation}
with the mixed boundary condition
\begin{equation}\label{eq:mixed-bdy}
   \xi_\nu(-T)=x_\nu
   ,\quad 
   \forall \nu\in\Lambda_+
   ,\qquad
   \xi_{-\nu}(T)=y_{-\nu}
   ,\quad 
   \forall \nu\in\Lambda_-.
\end{equation}
We make the variation of constant Ansatz $\xi_\ell(s)=c_\ell(s) e^{-a_\ell s}$.
Apply $\frac{d}{ds}$ to both sides and use~(\ref{eq:inhom-ODE}) to get
\[
   \p_s c_\ell(s)=\eta_\ell(s) e^{a_\ell s} .
\]
\textit{Positive eigenvalue $a_\nu$.}
Then we get $x_\nu=\xi_\nu(-T)=c_\nu(-T) e^{a_\nu T}$, so
$c_\nu(-T)=x_\nu e^{-a_\nu T}$.
Integrate $\p_s c_\nu(s)$ from $-T$ to $s$ to get
\[
   c_\nu(s)
   =x_\nu e^{-a_\nu T}+\int_{-T}^s \eta_\nu(t) e^{a_\nu t}\, dt
\]
and therefore
\begin{equation}\label{eq:xi12}
   \xi_\nu(s)
   =\underbrace{x_\nu e^{-a_\nu (T+s)}}_{=:\xi_\nu^1(s)}
   +\underbrace{\int_{-T}^s \eta_\nu(t) e^{a_\nu (t-s)}\, dt}_{=:\xi_\nu^2(s)} .
\end{equation}
\textit{Negative eigenvalue $a_{-\nu}$.} 
Then we get $y_{-\nu}=\xi_{-\nu}(T)=c_{-\nu}(T) e^{-a_{-\nu} T}$, so
$c_{-\nu}(T)=y_{-\nu} e^{a_{-\nu} T}$.
Integrate $\p_s c_{-\nu}(s)$ from $s$ to $T$ to get
\[
   c_{-\nu}(s)
   =y_{-\nu} e^{a_{-\nu} T}-\int_s^T \eta_{-\nu}(t) e^{a_{-\nu} t}\, dt
\]
and therefore
\[
   \xi_{-\nu}(s)
   =y_{-\nu} e^{a_{-\nu} (T-s)}-\int_s^T \eta_{-\nu}(t) e^{a_{-\nu} (t-s)}\, dt .
\]
To finish the proof of Claim~1 it suffices to show

\medskip\noindent
\textsc{Claim 2.} $\xi$ lies in $P_1(I_T)= L^2 (I_T,H_1)\cap W^{1,2}(I_T,H_0)$.

\smallskip\noindent
To we see this consider the case of a positive eigenvalue $a_\nu$
and write $\xi_\nu(s)=\xi_\nu^1(s)+\xi_\nu^2(s)$, see~(\ref{eq:xi12}).
To estimate $\xi_\nu^2$ define a function $g_\nu\colon\R\to\R$ by
$$
   g_\nu(s)
   :=\begin{cases}
      e^{-a_\nu s}&\text{, $s\ge 0$,}\\
      0&\text{, $s<0$.}
   \end{cases}
$$
We estimate the $L^1(I_T):=L^1([-T,T],\R)$ norm of $g_\nu$ by
\begin{equation*}
\begin{split}
   \norm{g_\nu}_{L^1(I_T)}
   =\int_{-T}^T g(s) \, ds
   =\int_0^T e^{-a_\nu s} \, ds
   =\frac{1-e^{-a_\nu T}}{a_\nu}\le \frac{1}{a_\nu} .
\end{split}
\end{equation*}
Thus, writing $\xi_\nu^2(s)=(\eta_\nu * g_\nu)(s)$
and by Young's inequality, we obtain
\begin{equation*}
\begin{split}
   \norm{\xi_\nu^2}_{L^2(I_T)}
   =\norm{\eta_\nu * g_\nu}_{L^2(I_T)}
   \le \norm{\eta_\nu}_{L^2(I_T)}\norm{g_\nu}_{L^1(I_T)}
   \le \tfrac{1}{a_\nu}\norm{\eta_\nu}_{L^2(I_T)}.
\end{split}
\end{equation*}
We estimate $\xi_\nu^1$ as follows
\begin{equation*}
\begin{split}
   \norm{\xi_\nu^1}_{L^2(I_T)}^2
   =\int_{-T}^T x_\nu^2 e^{-2a_\nu (T+s)}\, ds
   =x_\nu^2\frac{1-e^{-4a_\nu T}}{2a_\nu}
   \le \frac{x_\nu^2}{2a_\nu} .
\end{split}
\end{equation*}
Now we use these estimates to obtain for the sum
$\xi_\nu=\xi_\nu^1+\xi_\nu^2$ that
\begin{equation}\label{eq:xi-nu-pos}
\begin{split}
   \norm{\xi_\nu}_{L^2(I_T)}^2
   \le 2 \norm{\xi_\nu^1}_{L^2(I_T)}^2+2 \norm{\xi_\nu^2}_{L^2(I_T)}^2
   \le \frac{x_\nu^2}{a_\nu}+\frac{2}{a_\nu^2}
   \norm{\eta_\nu}_{L^2(I_T)}^2 .
\end{split}
\end{equation}
In the case of a negative eigenvalue $a_{-\nu}$
we observe that
$$
   \xi_{-\nu}(-s)
   =y_{-\nu} e^{a_{-\nu}(T+s)}
   -\int_{-T}^s \eta_{-\nu}(-\tau) e^{-a_{-\nu}(\tau-s)} \, d\tau .
$$
Comparing this expression with~(\ref{eq:xi12})
shows that we get the analogous estimate
\begin{equation}\label{eq:xi--nu-pos}
\begin{split}
   \norm{\xi_{-\nu}}_{L^2(I_T)}^2
   \le \frac{y_{-\nu}^2}{-a_{-\nu}}+\frac{2}{a_{-\nu}^2}
   \norm{\eta_{-\nu}}_{L^2(I_T)}^2 .
\end{split}
\end{equation}
Now we show that $\xi\in L^2 (I_T,H_1)$, namely
\begin{equation*}
\begin{split}
   \int_{-T}^T \norm{\xi(s)}_{H_1}^2 \, ds
   &\stackrel{1}{=}\int_{-T}^T \sum_{\ell\in\Lambda} a_\ell^2 \xi_\ell(s)^2 \, ds\\
   &=\sum_{\ell\in\Lambda} a_\ell^2 \norm{\xi_\ell}_{L^2(I_T)}^2\\
   &\stackrel{3}{=}\sum_{\nu\in\Lambda_-} a_{-\nu}^2 \norm{\xi_{-\nu}}_{L^2(I_T)}^2
   +\sum_{\nu\in\Lambda_+} a_\nu^2 \norm{\xi_\nu}_{L^2(I_T)}^2\\
   &\stackrel{4}{\le} 2 \sum_{\nu\in\Lambda_-} \norm{\eta_{-\nu}}_{L^2(I_T)}^2
   +\sum_{\nu\in\Lambda_-} - a_{-\nu} y_{-\nu}^2
\\
   &\quad
   +2 \sum_{\nu\in\Lambda_+}\norm{\eta_{\nu}}_{L^2(I_T)}^2
   +\sum_{\nu\in\Lambda_+} a_{\nu} x_{\nu}^2\\
   &\stackrel{5}{=} 2 \sum_{\ell\in\Lambda} \norm{\eta_{\ell}}_{L^2(I_T)}^2
   +\norm{y}_{\frac12}^2 +\norm{x}_{\frac12}^2\\
   &=2\int_{-T}^T \norm{\eta(s)}_{H_0}^2 \, ds
   +\norm{y}_{\frac12}^2 +\norm{x}_{\frac12}^2\\
   &=2\norm{\eta}_{P_0(I_T)}^2
   +\norm{y}_{\frac12}^2 +\norm{x}_{\frac12}^2 .
\end{split}
\end{equation*}
Equality 1 uses that $\A$ is an isometry and the fact
that the basis $\Vv(\A)$ consists of eigenvectors of $\A$.
Equality 3 uses the decomposition~(\ref{fig:Lambda-A}) of $\Lambda$.
Inequality 4 uses~(\ref{eq:xi-nu-pos}) and~(\ref{eq:xi--nu-pos}).
To see equality 5 go backwards and write $x=\sum_{\nu\in\Lambda_+} x_\nu v_\nu$.
Then use that the basis is $1/2$-orthogonal
and that the $1/2$-length of $v_\nu$ is $\sqrt{a_\nu}$
by~(\ref{eq:1/2-length-v_ell}). Similarly for $y$.

It remains to show that $\p_s\xi\in L^2 (I_T,H_0)$.
To see this we consider the case of a positive eigenvalue $a_\nu$.
Use~(\ref{eq:inhom-ODE}) in 1 and~(\ref{eq:xi-nu-pos}) in 3 to obtain
\begin{equation*}
\begin{split}
   \norm{\p_s\xi_\nu}_{L^2(I_T)}^2
   &\stackrel{1}{=}\norm{\eta_\nu-a_\nu\xi_\nu}_{L^2(I_T)}^2\\
   &\le 2 \norm{\eta_\nu}_{L^2(I_T)}^2+2 a_\nu^2\norm{\xi_\nu}_{L^2(I_T)}^2\\
   &\stackrel{3}{\le} 6 \norm{\eta_\nu}_{L^2(I_T)}^2+2 a_\nu x_\nu^2.
\end{split}
\end{equation*}
Similarly, by using~(\ref{eq:xi--nu-pos}) instead
of~(\ref{eq:xi-nu-pos}) we obtain
\begin{equation*}
\begin{split}
   \norm{\p_s\xi_{-\nu}}_{L^2(I_T)}^2
   \le 6 \norm{\eta_{-\nu}}_{L^2(I_T)}^2-2 a_{-\nu} y_{-\nu}^2.
\end{split}
\end{equation*}
Similarly as above we obtain the estimate
\begin{equation*}
\begin{split}
   \int_{-T}^T \norm{\p_s\xi(s)}_{H_0}^2\, ds
   =\sum_{\ell\in\Lambda}\norm{\p_s\xi_\ell}_{L^2(I_T)}^2
   \le 6 \norm{\eta}_{P_0(I_T)}^2+2\norm{y}_{\frac12}^2 +2\norm{x}_{\frac12}^2 .
\end{split}
\end{equation*}
We have shown that $\xi\in L^2 (I_T,H_1)\subset L^2 (I_T,H_0)$ and
$\p_s\xi\in L^2 (I_T,H_0)$. Thus $\xi\in P_1(I_T;H_1,H_0)$ and
\begin{equation*}
\begin{split}
   \norm{\xi}_{P_1(I_T;H_1,H_0)}^2
   \le 10 \norm{\eta}_{P_0(I_T)}^2+4\norm{y}_{\frac12}^2 +4\norm{x}_{\frac12}^2 .
\end{split}
\end{equation*}
This concludes the proof of Claim~2, hence of Claim~1 and Step~1.
\end{proof}

From now on we abbreviate $A_\sigma:=A(\sigma)\in\Ll(H_1,H_0)$.

\smallskip
\noindent
\textbf{Step~2} (Small interior interval)\textbf{.}
There is a finite subset $\Lambda^\prime\subset\R$
and a constant $C_2>0$
such that for every $\beta\in C^\infty(I_T,\R)$ which vanishes
on the interval boundary, in symbols $\beta(-T)=\beta(T)=0$,
and has the property
$$
   \sup_{\sigma,\tau\in\supp\beta}\norm{A_\sigma-A_\tau}_{\Ll(H_1,H_0)}\le\frac{1}{C_2}
$$
it holds that
\[
   \Norm{\beta\xi}_{P_1(I_T)}
   \le C_2\left(
   \Norm{\beta D_A\xi}_{P_0(I_T)}+\Norm{\beta^\prime\xi}_{P_0(I_T)}
   +\Norm{\beta\xi}_{P_0(I_T)}
    \right)
\]
for every $\xi\in P_1(I_T)$.

\begin{proof}
This is Step~6 in the proof of the Rabier Theorem~\ref {thm:Rabier}.
\end{proof}

\smallskip
\noindent
\textbf{Step~3} (Small interval at right boundary)\textbf{.}
There exist constants $\eps_3>0$ and $C_3>0$ such that for every
$\beta\in C^\infty(I_T,\R)$ which vanishes on the left interval
boundary, in symbols $\beta(-T)=0$, and has the property
$$
   \sup_{\sigma,\tau\in\supp\beta}\norm{A_\sigma-A_\tau}_{\Ll(H_1,H_0)}\le\eps_3
$$
it holds that
\[
   \Norm{\beta\xi}_{P_1(I_T)}
   \le C_3\left(
   \Norm{\beta D_A\xi}_{P_0(I_T)}+\Norm{\beta^\prime\xi}_{P_0(I_T)}
   +\Norm{\beta\xi}_{P_0(I_T)}
   +\norm{\pi_-\beta_T\xi_T}_{\frac12}
    \right)
\]
for every $\xi\in P_1(I_T)$.

\begin{proof}
By continuity of the map $s\mapsto A(s)=:A_s$ there exists the limit
\[
   \lim_{\sigma\to T}
   A_\sigma
   =\A_T .
\]
By Step~1 the augmented operator associated to the constant path
$\A_T$, namely
\begin{equation*}
\begin{split}
   \mathfrak{D}_{\A_T}\colon P_1(I_T)
   &\to P_0(I_T)\times H^+_{\frac12}(\A_T) \times H^-_{\frac12}(\A_T)
   =: \Ww(I_T;\A_T,\A_T)
   \\
   \xi
   &\mapsto \left(D_{\A_T}\xi,\pi_+\xi_{-T},\pi_-\xi_T\right)
\end{split}
\end{equation*}
is bijective, in particular invertible. Together these two facts imply that
there is $\eps_3>0$ such that the following is true.
If $\norm{A_\sigma-\A_T}_{\Ll(H_1,H_0)}\le\eps_3$, then
$\mathfrak{D}_{A_\sigma}$ is still invertible with
inverse bound
\[
   \Norm{(\mathfrak{D}_{A_\sigma})^{-1}}_{\Ll(\Ww(I_T;\A_t,\A_T),P_1(I_T))}
   \le\frac{1}{2\eps_3} .
\]
This follows from the fact that the map
$$
   T\colon \Ff^*\to\Ll(P_1,P_0)
   ,\quad A\mapsto D_A
$$
is continuous, by~(\ref{eq:D(s)-cont}),
in view of the injectivity Lemma~\ref{le:invertibility-open}.

Now pick $\sigma\in\supp\beta$. Then, in particular, the previous
estimate for the inverse is valid.
By formula~(\ref{eq:jhu678-2}) in Step~6 in the proof of the Rabier
Theorem~\ref{thm:Rabier} we get a formula for the first component
$D_{A_\sigma}\beta\xi$, namely
\begin{equation*}
\begin{split}
   \mathfrak{D}_{A_\sigma}\beta\xi
   &=\left(D_{A_\sigma}\beta\xi\,,\,
   \pi_+\beta_{-T}\xi_{-T}\,,\,\pi_-\beta_T\xi_T \right)\\
   &=\left(\beta^\prime\xi+\beta D_A\xi
   +(A_\sigma-A)\beta\xi 
   \,,\,
    0\,,\, \pi_-\beta_T\xi_T \right) .
\end{split}
\end{equation*}
In the second equality we use the assumption $\beta_{-T}=0$.
Since the operator
$\mathfrak{D}_{A_\sigma}$ is invertible we
can write
\begin{equation*}
\begin{split}
   \beta\xi
   &=(\mathfrak{D}_{A_\sigma})^{-1}
   \left(\beta^\prime\xi+\beta D_A\xi+(A_\sigma-A)\beta\xi\,,\,
    0\,,\, \pi_-\beta_T\xi_T \right) .
\end{split}
\end{equation*}
Taking norms we estimate
\begin{equation*}
\begin{split}
   \Norm{\beta\xi}_{P_1(I_T)}
   &\le \frac{1}{\eps_3}
   \Bigl(
   \Norm{\beta^\prime\xi}_{P_0(I_T)}
   +\Norm{\beta D_A\xi}_{P_0(I_T)}
   +\Norm{(A_\sigma-A)\beta\xi}_{P_0(I_T)}\\
   &\quad
   +\norm{\pi_-\beta_T\xi_{T}}_{\frac12}
   \Bigr) .
\end{split}
\end{equation*}
Now we estimate
\begin{equation*}
\begin{split}
   &\Norm{(A_\sigma-A)\beta\xi}_{P_0(I_T)}^2
\\
   &=\int_{-T}^T \Norm{\left(A_\sigma-A_s\right)\beta_s\xi_s}_{H_0}^2\, ds
\\
   &\le \int_{\supp\beta} 
   \Norm{A_\sigma-A_s}_{\Ll(H_1,H_0)}^2
   \Norm{\beta_s\xi_s}_{H_1}^2 ds\\
   &\le \eps_3^2\Norm{\beta\xi}_{P_1(I_T)}^2 .
\end{split}
\end{equation*}
The last two estimates together imply Step~5
with $C_3:=\frac{2}{\eps_3}$.
\end{proof}

\smallskip
\noindent
\textbf{Step~4} (Small interval at left boundary)\textbf{.}
There exist constants $\eps_4>0$ and $C_4>0$ such that for every
$\beta\in C^\infty(I_T,\R)$ which vanishes on the right interval
boundary, in symbols $\beta(T)=0$, and has the property
$$
   \sup_{\sigma,\tau\in\supp\beta}\norm{A_\sigma-A_\tau}_{\Ll(H_1,H_0)}\le\eps_4
$$
it holds that
\begin{equation*}
\begin{split}
   &\Norm{\beta\xi}_{P_1(I_T)}\\
   &\le C_4\left(
   \Norm{\beta D_A\xi}_{P_0(I_T)}+\Norm{\beta^\prime\xi}_{P_0(I_T)}
   +\Norm{\beta\xi}_{P_0(I_T)}
   +\norm{\pi_+\beta_{-T}\xi_{-T}}_{\frac12}
    \right)
\end{split}
\end{equation*}
for every $\xi\in P_1(I_T)$.

\begin{proof}
Same argument as in Step~3.
\end{proof}

\smallskip
\noindent
\textbf{Step~5} (Partition of unity)\textbf{.}
We prove Theorem~\ref{thm:finite-interval}.

\begin{proof}
Let $\eps:=\min\{\eps_2,\eps_3,\eps_4\}$ and $C:=\max\{C_2,C_3,C_4\}$.
Choose a finite partition of unity $\{\beta_j\}_{j=0}^{M+1}$ for $I_T=[-T,T]$
with the properties
$$
   \beta_0(-T)=1,\quad\beta_0(T)=0,
   ,\qquad
   \beta_{M+1}(-T)=0,\quad\beta_{M+1}(T)=1,
$$
and
$$
   \sup_{\sigma,\tau\in\supp\beta_i}\norm{A_\sigma-A_\tau}_{\Ll(H_1,H_0)}
   \le\eps
   ,\qquad
  \supp\beta_j\subset (-T,T) ,
$$
for $i=0,1,\dots,M,M+1$ and $j=1,\dots,M$.
That such a partition exists follows from the continuity
of $s\mapsto A(s)$ and since on the compact set $[-T,T]$
continuity becomes uniform continuity.
Let $\xi\in P_1(I_T)$. By Steps~4 and~3 we get
\begin{equation*}
\begin{split}
   \Norm{\beta_0\xi}_{P_1(I_T)}
   &\le C\Bigl(
   \Norm{\beta_0
     D_A\xi}_{P_0(I_T)}+\Norm{\beta_0^\prime\xi}_{P_0(I_T)}\\
   &\qquad
   +\Norm{\beta_0\xi}_{P_0(I_T)}
   +\norm{\pi_+\xi_{-T}}_{\frac12}
    \Bigr)
\\
   \Norm{\beta_{M+1}\xi}_{P_1(I_T)}
   &\le C\Bigl(
   \Norm{\beta_{M+1}
     D_A\xi}_{P_0(I_T)}+\Norm{\beta_{M+1}^\prime\xi}_{P_0(I_T)}
\\
   &\qquad
   +\Norm{\beta_{M+1}\xi}_{P_0(I_T)}
   +\norm{\pi_-\xi_{T}}_{\frac12}
   \Bigr) .
\end{split}
\end{equation*}
By Step~2 we have
\[
   \Norm{\beta_j\xi}_{P_1(I_T)}
   \le C \left(
   \Norm{\beta_j D_A\xi}_{P_0(I_T)}+\Norm{\beta_j^\prime\xi}_{P_0(I_T)}
   +\Norm{\beta_j\xi}_{P_0(I_T)}
    \right)
\]
for $j=1,\dots,M$.
We abbreviate $B:=\max\{\norm{\beta_0^\prime}_\infty,
\norm{\beta_1^\prime}_\infty,\dots, \norm{\beta_{M+1}^\prime}_\infty\}$.
Putting these estimates together we obtain
\begin{equation*}
\begin{split}
   \Norm{\xi}_{P_1(I_T)}
   &\le\sum_{j=0}^{M+1} \Norm{\beta_j\xi}_{P_1(I_T)}\\
   &\le C \sum_{j=0}^{M+1}\left(
   \Norm{\beta_j D_A\xi}_{P_0(I_T)}+\Norm{\beta_j^\prime\xi}_{P_0(I_T)}
   +\Norm{\beta_j\xi}_{P_0(I_T)}\right)
\\
   &\quad
   +C\norm{\pi_+\xi_{-T}}_{\frac12}
   +C\norm{\pi_-\xi_{T}}_{\frac12}
   \\
   &\le C (M+2) \Norm{D_A\xi}_{P_0(I_T)}
   +C(B+1) (M+2) \Norm{\xi}_{P_0(I_T)}
\\
   &\quad
   +C\norm{\pi_+\xi_{-T}}_{\frac12}
   +C\norm{\pi_-\xi_{T}}_{\frac12} .
\end{split}
\end{equation*}
Setting
$
   c:=C(B+1)(M+2)
$
proves Step~5.
\end{proof}
\noindent
The proof of Theorem~\ref{thm:finite-interval} is complete.
\end{proof}

\boldmath
\subsubsection{Estimate for the adjoint $D_A^*$}
\unboldmath

Let $A\in\Aa_{I_T}^*$.
We call the following operator the
\textbf{adjoint of \boldmath$D_A$}, namely
$$
   D_A^*:=D_{-A^*}\colon P_1(I_T;H_0^*,H_1^*)\to P_0(I_T;H_1^*)
   ,\quad
   \eta\mapsto \p_s\eta-A(s)^*\eta .
$$

\begin{corollary}\label{cor:finite-interval-*}
For $A\in\Aa_{I_T}^*$ there exists a constant $c>0$ such that
\begin{equation*}
\begin{split}
   \norm{\eta}_{P_1(I_T;H_0^*,H_1^*)}
   &\le c\Bigl(\norm{\eta}_{P_0(I_T;H_1^*)}+\norm{D_A^*\eta}_{P_0(I_T;H_1^*)}
   \\
   &\qquad
   +\norm{\pi_+^{-\A^*_{-T}}\eta(-T)}_{\frac12}
   +\norm{\pi_-^{-\A^*_T}\eta(T)}_{\frac12}\Bigr)
\end{split}
\end{equation*}
for every $\eta\in P_1(I_T;H_0^*,H_1^*)$.
\end{corollary}

\begin{proof}
Theorem~\ref{thm:finite-interval} and Lemma~\ref{le:adjoint-map};
see also Remark~\ref{rem:symmetrizable}.
\end{proof}

\boldmath
\subsubsection{Fredholm under boundary conditions: $D_A^{+-}$}
\unboldmath

Let $A\in\Aa_{I_T}^*$.
For the spectral projections $\pi_+^{\A_{-T}}$ and
$\pi_-^{\A_T}$ see~(\ref{eq:H_1/2}).
To turn Theorem~\ref{thm:finite-interval}
to a semi-Fredholm estimate we restrict the domain
of the operator 
$$
   D_A\colon P_1(I_T;H_1,H_0)\to P_0(I_T;H_0)
   ,\quad
   \xi\mapsto \p_s\xi+A(s)\xi
$$
by imposing appropriate boundary conditions
that cut down the operator kernel to finite dimension.
To this end we define a subspace of the domain as follows
\begin{equation}\label{eq:P_1^+-}
\begin{split}
   &P_1^{+-}(I_T,\A_{\pm T};H_1,H_0)\\
   &:=\{\xi\in P_1(I_T;H_1,H_0)\mid \pi_+^{\A_{-T}}\xi_{-T}=0\;\wedge\;
   \pi_-^{\A_{T}}\xi_{T}=0\} .
\end{split}
\end{equation}
The associated restriction of $D_A$ we denote by
\begin{equation}\label{eq:D_A^+-}
\boxed{
   D_A^{+-}=\p_s+A\colon P_1^{+-}(I_T,\A_{\pm T};H_1,H_0) \to P_0(I_T;H_0) .
}
\end{equation}

\noindent
To $A\in\Aa_{I_T}^*$ we associate
the path $-A^*\in\Aa_{I_T}^*$,
namely $s\mapsto -A(s)^*\colon H_0^*\to H_1^*$.
Define a Hilbert space
$
   P_1(I_T;H_0^*,H_1^*)
   :=L^2(I_T,H_0^*)\cap W^{1,2}(I_T,H_1^*)
$,
analogous to~(\ref{eq:P_1-I}).
A closed linear subspace is defined by imposing boundary conditions
\begin{equation}\label{eq:P_1^*}
\begin{split}
   &P_1^{+-}(I_T,-\A^*_{\pm T};H_0^*,H_1^*)\\
   &:=\{\eta\in P_1(I_T;H_0^*,H_1^*)\mid \pi_+^{-\A_{-T}^*}\eta_{-T} =0\;\wedge\;
   \pi_-^{-\A_{T}^*}\eta_T=0\} .
\end{split}
\end{equation}
It includes into
$
   P_1^{+-}(I_T,-\A^*_{\pm T};H_0^*,H_1^*)
   \subset P_1(I_T;H_0^*,H_1^*)
   \stackrel{\text{(\ref {eq:P_1-continuous})}}{\subset}
   C^0(I_T,H_1^*)
$.
The restriction of the linear operator
\begin{equation*}
   D_{-A^*}\colon P_1(I_T;H_0^*,H_1^*)\to P_0 (I_T;H_1^*)
   ,\quad
   \eta\mapsto \p_s\eta-A(s)^*\eta
\end{equation*}
to the closed linear subspace~(\ref{eq:P_1^*}) is denoted by
\begin{equation}\label{eq:D_-A*^+-}
\begin{split}
\boxed{
   D_{-A^*}^{+-}=\p_s-A^*\colon P_1^{+-}(I_T,-\A^*_{\pm T};H_0^*,H_1^*)
   \to P_0 (I_T;H_1^*) .
}
\end{split}
\end{equation}

\begin{corollary}[Semi-Fredholm]\label{cor:semi-Fredholm}
For any $A\in\Aa_{I_T}^*$ the operators
$$
   D_A^{+-}\colon P_1^{+-}(I_T,\A_{\pm T};H_1,H_0)\to P_0(I_T;H_0)
$$
and
$$
   D_{-A^*}^{+-}\colon P_1^{+-}(I_T,-\A^*_{\pm T};H_0^*,H_1^*)
   \to P_0(I_T;H_1^*)
$$
are semi-Fredholm {\color{gray}(finite dimensional
kernel and closed image)}.
\end{corollary}

\begin{proof}
Theorem~\ref{thm:finite-interval} provides
for $D_A^{+-}$ the semi-Fredholm estimate\footnote{
  The inclusion map
  $P_1(I_T)\to P_0(I_T)$ is compact,
  see e.g.~\cite[Lemma~3.8]{robbin:1995a}.
  Hence the semi-Fredholm property follows
  from~\cite[Lemma~A.1.1]{mcduff:2004a}.
  }
\[
   \norm{\xi}_{P_1(I_T;H_1,H_0)}
   \le c\left(\norm{D_A^{+-}\xi}_{P_0(I_T;H_0)}+\norm{\xi}_{P_0(I_T;H_0)}\right)
\]
$\forall\xi\in P_1^{+-}(I_T,\A_{\pm T};H_1,H_0)$.
Corollary~\ref{cor:finite-interval-*} provides
the semi-Fredholm estimate
\[
   \norm{\eta}_{P_1(I_T;H_0^*,H_1^*)}
   \le c\left(\norm{D_{-A^*}^{+-}\eta}_{P_0(I_T;H_1^*)}+\norm{\eta}_{P_0(I_T;H_1^*)}\right)
\]
for the operator $D_{-A^*}^{+-}$
and every $\eta\in P_1^{+-}(I_T,-\A^*_{\pm T};H_0^*,H_1^*)$.
\end{proof}

\begin{theorem}[Fredholm]\label{thm:Fredholm-I_T-D_A+-}
For any Hessian path $A\in \Aa_{I_T}^*$ the operator
$D_A^{+-}\colon P_1^{+-}(I_T,\A_{\pm T};H_1,H_0)\to P_0(I_T;H_0)$
is Fredholm.
\end{theorem}

\begin{corollary}\label{cor:D_A-I_T-closed-image}
The operator $D_A\colon P_1(I_T)\to P_0(I_T)$
in~(\ref{eq:D_A}) has closed image of finite co-dimension
for any Hessian path $A\in\Aa_{I_T}^*$.
\end{corollary}

\begin{proof}
By Theorem~\ref{thm:Fredholm-I_T-D_A+-} the image of $D_A^{+-}$
is closed and of finite co-dimension.
Since $D_A^{+-}$ is a restriction of $D_A$ we have inclusion
$\im D_A^{+-}\subset \im D_A$.
So $\im D_A$ is of finite co-dimension.
Thus $\im D_A$ is closed by~\cite[Prop.\,11.5]{brezis:2011a}.
\end{proof}
\begin{figure}[h]
\begin{center}
\begin{tabular}{lll}
      
  & $D_A\colon P_1\to P_0$         
  & $D_A^{+-}\colon P_1^{+-}\to P_0$ 
\\
\midrule
  $\dim \ker$
  & \small ${\color{red}\infty}$
  & \small $k<\infty$ 
\\
  $\dim \coker$
  & \small $\le\ell$ \qquad\qquad ${\color{gray}\Leftarrow}$
  & \small $\ell<\infty$
\\
\midrule
     \small 
  & \small co-semi-Fredholm
  & \small Fredholm 
\\
\midrule
  $\mathrm{image}$
  & \small closed \qquad\quad ${\color{gray}\Leftarrow}$
  & \small closed
\\
  $\coker$
  & \small 
  & \small $\coker D_A^{+-}\simeq\ker D_{-A^*}^{+-}$
\\
  $\ker$
  & \small {\color{red} huge}
  & \small {\color{gray}$\ker D_A^{+-}\simeq\coker D_{-A^*}^{+-}$}
\end{tabular}
\caption{$D_A=\p_s+A(s)$ on $P_1(I_T)$ and its restriction $D_A^{+-}$ to $P_1^{+-}$}
\label{fig:D_A-finite-interval}
\end{center}
\end{figure}
\begin{proof}[Proof of Theorem~\ref{thm:Fredholm-I_T-D_A+-}]
Pick $A\in\Aa_{I_T}^*$, then $\A_{-T}:=A(-T)$ and $\A_T:=A(T)$ are
invertible.
By Corollary~\ref{cor:semi-Fredholm}
the operator  $D_A^{+-}$ (and also $D_{-A^*}^{+-}$)
has finite dimensional kernel and closed image.
It remains to show that $D_A^{+-}$ has finite dimensional co-kernel.
This is proved in the following Proposition~\ref{prop:coker-ker-I_T}.
The proof of Theorem~\ref{thm:Fredholm-I_T-D_A+-} is complete.
\end{proof}

To prove that $D_A^{+-}$ has finite dimensional co-kernel we show
how the annihilator of $D_A^{+-}$
can be identified with the kernel of the semi-Fredholm operator
$D_{-A^*}^{+-}$.
The \textbf{annihilator} of $D_A^{+-}$ consists of all linear
functionals on the co-domain which vanish along the image
$$
   \Ann{(D_A^{+-})}
   :=\{\zeta\in P_0(I_T;H_0)^*\mid \zeta(D_A\xi)=0
   \;\;\forall \xi\in P_1^{+-}(I_T,\A_{\pm T};H_1,H_0)\}.
$$
Since $\zeta(D_A\xi)=\INNER{\zeta}{D_A\xi}_0$
the annihilator identifies naturally with the orthogonal complement of
the image of $D_A^{+-}$, which is the cokernel, in symbols
$$
   \Ann{(D_A^{+-})}\simeq \coker{D_A^{+-}} .
$$
We have a natural map
$$
   \Kk\colon P_0(I_T;H_0^*)
   \to P_0(I_T;H_0)^*
$$
defined by
$$
   (\Kk\eta)\chi:=\int_{I_T}\eta(s)\chi(s)\, ds
$$
for every $\chi\in P_0(I_T;H_0) =L^2(I_T,H_0)$.
As shown by Kreuter~\cite[Thm.\,2.22]{Kreuter:2015a}
the map $\Kk$ is an isometry.

\begin{proposition}\label{prop:coker-ker-I_T}
The vector spaces 
$
   \Kk^{-1}\Ann{(D_A^{+-})}
   =
   \ker{D_{-A^*}^{+-}}
$
coincide as vector subspaces of $P_0(I_T;H_0^*)$.
\end{proposition}

\begin{proof}
Note that both are subspaces of $P_0(I_T;H_0^*) =L^2(I_T,H_0^*)$, indeed
$$
   \Kk^{-1}\Ann{(D_A^{+-})}
   \subset P_0(I_T;H_0^*)
   \supset P_1^*(I_T)
   \supset P_1^{*,+-}(I_T)
   \supset \ker{D_{-A^*}^{+-}} .
$$
For equality $\Kk^{-1}\Ann{(D_A^{+-})}=\ker{D_{-A^*}^{+-}}$ we
show inclusions I. $\subset$ and II. $\supset$.

\medskip\noindent
\textbf{I. The inclusion \boldmath''$\subset$''.}
Pick $\eta \in \Kk^{-1}\Ann{(D_A^{+-})}\subset L^2(I_T,H_0^*)$.
For $\xi\in P_1^{+-}(I_T)$ we calculate
\begin{equation}\label{eq:eta-I_T}
\begin{split}
   0
   &=(\Kk\eta)D_A\xi
\\
   &=\int_{I_T} \eta(D_A\xi) \, ds
\\
   &\stackrel{3}{=}\int_{I_T}\eta(\p_s\xi)\, ds+\int_{I_T}\eta(A\xi)\, ds
\\
   &=\int_{I_T}\eta(\p_s\xi)\, ds+\int_{I_T} (A^*\eta)\xi\, ds .
\end{split}
\end{equation}
This shows that $\eta$ admits a weak derivative in $H_1^*$, notation
$\p_s\eta$, with
$$
   \p_s\eta-A^*\eta=0 .
$$
We first show that this equation implies that $\p_s\eta\in L^2(I_T,H_1^*)$.
To see this we first note that since $A\in\Aa^*_{I_T}$ and $I_T$ is compact
there is a constant $c$ such that $\norm{A(s)}_{\Ll(H_1,H_0)}\le c$
for every $s\in I_T$.
The map $*\colon \Ll(H_1,H_0)\to \Ll(H_0^*,H_1^*)$, $A\mapsto A^*$,
is an isometry. Hence we also have
$\norm{A(s)^*}_{\Ll(H_0^*,H_1^*)}\le c$ for every $s\in I_T$.
Using this we obtain finiteness of
$$
   \norm{\p_s \eta}_{L^2(I_T,H_1^*)}^2
   =\norm{A^* \eta}_{L^2(I_T,H_1^*)}^2
   =\int_{I_T} \norm{A(s)^*\eta(s)}_{H_1^*}^2\, ds
   \le c^2 \norm{\eta}_{L^2(I_T,H_0^*)}^2 .
$$
Indeed, since $\eta\in L^2(I_T,H_0^*)$, it follows that
$\p_s \eta\in L^2(I_T,H_1^*)$.
To summarize, we proved that
$$
   \eta\in P_1(I_T;H_0^*,H_1^*)
   \;\wedge\;
   \eta\in\ker{D_{-A^*}} .
$$

\smallskip
It remains to show that $\eta$ satisfies the boundary
conditions~(\ref{eq:P_1^*}).
We check the boundary condition at $-T$,
namely
$$
   0=\pi_+^{-\A_{T}^*}\eta(-T) =\pi_-^{\A_T^*}\eta(-T) ,
$$
the boundary condition at $T$ one checks analogously.
By Section~\ref{sec:int-spaces-1/2}
the boundary condition does not depend on the choice of the
inner products on $H_1$ and $H_0$, therefore we can assume without
loss of generality that the inner products are $\A_{-T}$-adapted,
i.e. from now on
\begin{equation*}\label{eq:hyp-A_T}
   \text{$\A_{-T}\colon H_1\to H_0$ is a symmetric isometry.}
\end{equation*}
We pick an orthonormal basis $\Vv(\A_{-T})=\{v_\ell\}_{\ell\in\Lambda}\subset H_1$
of $H_0$, see~(\ref{eq:ONB-A}), which consists of eigenvectors of $\A_{-T}$,
more precisely $\A_{-T} v_\ell=a_\ell v_\ell$. 
From now on we identify $\A^*\colon H_0^*\to H_1^*$ isometrically with
$A\colon H_1\to H_0$ according to Lemma~\ref{le:A*-A}.

We have that
$$
   \eta
   \in L^2(I_T;H_0)\cap W^{1,2}(I_T;H_{-1})
   =P_1(I_T;H_{-1},H_0)
   \subset C^0(I_T,H_{-1}).
$$
Here the last inclusion follows from~\cite[Le.\,7.1]{Roubicek:2013a};
see also~(\ref{eq:Oliver}).

Since $\eta$ is continuous, it makes sense to consider $\eta$
pointwise at any time $s$ and use the orthogonal
basis $\Vv(\A_{-T})$ of $H_{-1}$ to write
$$
   \eta(s)=\sum_{\ell\in\Lambda} \eta_\ell(s) v_\ell
   \in H_{-1}
$$
where $\eta_\ell(s)\in\R$ depends continuously on $s$.
Moreover, the norm of $\eta$ is related to the norms of the
coefficients as follows
$$
   \norm{\eta}_{P_1(I_T;H_{-1},H_0)}^2
   =\sum_{\ell\in\Lambda}\left(\norm{\eta_\ell}_{L^2(I_T,\R)}^2
   +\frac{1}{a_\ell^2}\norm{\eta_\ell}_{W^{1,2}(I_T,\R)}^2\right) .
$$
Since $\pi_-^{\A_{-T}}\xi(-T)$ is arbitrary and $\pi_+^{\A_{-T}}\xi(-T)=0$
we prove that~(\ref{eq:eta-I_T}) implies

\medskip
\noindent
\textbf{Claim.}
$\eta_{-\nu}(-T)=0$ for every $\nu\in\Lambda_-=\Lambda_-(\A_{-T})$.

\begin{proof}[Proof of Claim]
We pick a smooth cut-off function $\beta\colon [-T,T]\to[0,1]$
such that $\beta(-T)=1$, that $\beta\equiv 0$ outside a small interval
$[-T,-T+\delta]$, and that $\norm{\beta^\prime}_{L^\infty}\le\frac{2}{\delta}$.
Pick $\nu\in\Lambda_-$ and consider the corresponding eigenvector
$v_{-\nu}$.
Let us define a map $\zeta\colon [-T,T]\to \R v_{-\nu}$ and conclude
the following two properties
\begin{equation}\label{eq:zeta}
\begin{split}
\boxed{
   \zeta
   :=\beta\eta_{-\nu}(-T)  v_{-\nu}
}\,
   ,\;\;
   \pi_+^{\A_{-T}}\zeta(-T)=0
   ,\;\;
   \zeta\in P_1^{+-}(I_T,\A_{\pm T};H_1,H_0) .
\end{split}
\end{equation}
Recall that $\A_{-T} v_{-\nu}=a_{-\nu} v_{-\nu}$ where $a_{-\nu}<0$.
Given $\eps>0$, pick a parameter
$0<\delta\le \min\{1,(8\abs{a_{-\nu}})^{-1}\}$ so small that
\begin{equation}\label{eq:777} 
\begin{split}
   &\norm{\eta}_{L^2(I_T,H_0)}^2
   \sup_{s\in[-T,-T+\delta]}\norm{A(s)-\A_{-T}}_{\Ll(H_1,H_0)}^2(a_{-\nu})^2 4\\
   &+\max\{16,2\Abs{a_{-\nu}},4(a_{-\nu})^2\}
   \sup_{s\in[-T,-T+\delta]} \Abs{\eta_{-\nu}(s)-\eta_{-\nu}(-T)}^2
   \\
   &\le\tfrac{1}{4} \eps^2.
\end{split}
\end{equation}
This will be used together with $ab\le\frac{a^2}{2}+\frac{b^2}{2}$
for terms 1-4 below. We calculate
\begin{equation*}
\begin{split}
   &\eta_{-\nu}(-T)^2\\
   &=-\int_{-T}^T \eta_{-\nu}(-T)^2 \beta^\prime(s) \, ds\\
   &=-\int_{-T}^T \INNER{\eta(-T)}{\p_s\zeta(s)}_{0} ds\\
   &\stackrel{3}{=}-\int_{-T}^T
   \INNER{\underline{\eta(s)}}{{\color{cyan} \p_s\zeta(s)}}_{0} ds
   +\int_{-T}^T \INNER{\underline{\eta(s)}-\eta(-T)}{\p_s\zeta(s)}_{0} ds
\\
   &\stackrel{4}{=}
   \int_{-T}^{-T+\delta}\INNER{\eta(s)}{({\color{cyan} A(s)}-\underline{\A_{-T}})\zeta(s)}_{0} \, ds
   +\int_{-T}^{-T+\delta}\INNER{\eta(s)-\eta(-T)}{\p_s\zeta(s)}_{0} \, ds\\
   &\quad
   +\int_{-T}^{-T+\delta}\INNER{\eta(s)}{\underline{\A_{-T}}\zeta(s)}_{0} \, ds
\\
   &\stackrel{5}{=}
   \int_{-T}^{-T+\delta}\INNER{\eta(s)}{(A(s)-\A_{-T})\zeta(s)}_{0} \, ds
   +\int_{-T}^{-T+\delta}\INNER{\eta(s)-\eta(-T)}{\p_s\zeta(s)}_{0} \, ds\\
   &\quad
   +\int_{-T}^{-T+\delta}\INNER{\eta(s)-\underline{\eta(-T)}}{\A_{-T}\zeta(s)}_{0} \, ds
   +\int_{-T}^{-T+\delta}\INNER{\underline{\eta(-T)}}{\A_{-T}\zeta(s)}_{0} \, ds
\end{split}
\end{equation*}
where in equalities 3, 4, and 5 we added \underline{zero},
equality 4 is by {\color{cyan}line 3} in~(\ref{eq:eta-I_T})
for $\xi:=\zeta$ and since
$\supp \zeta\subset \supp \beta\subset [-T,-T+\delta]$.
Now we discuss each of the four terms in the sum individually
using the estimate $ab\le \frac{a^2}{2}+\frac{b^2}{2}$.

\smallskip\noindent
\textbf{Term 1.}
By Cauchy-Schwarz and definition~(\ref{eq:zeta}) of $\zeta$ we obtain
\begin{equation*}
\begin{split}
   &\int_{-T}^{-T+\delta}\INNER{\eta(s)}{(A(s)-\A_{-T})\zeta(s)}_{0} ds\\
   &\le \int_{-T}^{-T+\delta}
   \norm{\eta(s)}_0\norm{A(s)-\A_{-T}}_{\Ll(H_1,H_0)}\norm{\zeta(s)}_1
   ds\\
   &\stackrel{2}{\le}\norm{\eta}_{L^2(I_T,H_0)}
   \sup_{s\in[-T,-T+\delta]}\norm{A(s)-\A_{-T}}_{\Ll(H_1,H_0)}
   \Abs{a_{-\nu}}2\cdot\tfrac12\Abs{\eta_{-\nu}(-T)}\\
   &\stackrel{3}{\le} \tfrac{\eps^2}{8}+\tfrac{1}{8} \eta_{-\nu}(-T)^2.
\end{split}
\end{equation*}
Inequality 2
uses that $\norm{v_{-\nu}}_1=\Abs{a_{-\nu}}$, by~(\ref{eq:1-length-v_ell}),
inequality 3 is by~(\ref{eq:777}).

\smallskip\noindent
\textbf{Term 2.}
By definition of $\zeta$ and $(v_\ell)$ being an orthonormal basis of $H_0$ we get
\begin{equation*}
\begin{split}
   &\int_{-T}^{-T+\delta}\INNER{\eta(s)-\eta(-T)}{\p_s\zeta(s)}_{0} ds\\
   &= \int_{-T}^{-T+\delta} 
   (\eta_{-\nu}(s)-\eta_{-\nu}(-T)) \beta^\prime(s) \eta_{-\nu}(-T) \INNER{v_{-\nu}}{v_{-\nu}}_0
   ds\\
   &\stackrel{2}{\le} \delta\sup_{s\in[-T,-T+\delta]}\Abs{\eta_{-\nu}(s)-\eta_{-\nu}(-T)}
   \tfrac{2}{\delta}2\cdot\tfrac12
   \Abs{\eta_{-\nu}(-T)}\\
   &\stackrel{3}{\le} \tfrac{\eps^2}{8}+\tfrac{1}{8} \eta_{-\nu}(-T)^2.
\end{split}
\end{equation*}
Inequality 2 pulls out the supremum,
uses $\norm{\beta^\prime}_{L^\infty}\le \tfrac{2}{\delta}$,
inequality 3 is by~(\ref{eq:777}).

\smallskip\noindent
\textbf{Term 3.}
By definition of $\zeta$ and $(v_\ell)$ being an orthonormal basis of $H_0$ we get
\begin{equation*}
\begin{split}
   &\int_{-T}^{-T+\delta}\INNER{\eta(s)-\eta(-T)}{\A_{-T}\zeta(s)}_{0} ds\\
   &\stackrel{1}{\le} \int_{-T}^{-T+\delta}
   \left(\underline{\eta_{-\nu}(s)-\eta_{-\nu}(-T)}\right) \beta(s)
   a_{-\nu}
   \left(\underline{\eta_{-\nu}(s){\color{gray}\,-\eta_{-\nu}(-T)}}
   {\color{gray}\,+\eta_{-\nu}(-T)}\right)
   ds\\
   &\stackrel{2}{\le} 
   \delta\sup_{s\in[-T,-T+\delta]}\Abs{\eta_{-\nu}(s)-\eta_{-\nu}(-T)}^2
   \Abs{a_{-\nu}}\\
   &\quad
   +\delta\sup_{s\in[-T,-T+\delta]}\Abs{\eta_{-\nu}(s)-\eta_{-\nu}(-T)}
   \Abs{a_{-\nu}} 2\cdot\tfrac12\Abs{\eta_{-\nu}(-T)}
\\
   &\stackrel{3}{\le} \tfrac{\eps^2}{8}
   +\tfrac{\eps^2}{8}+\tfrac{1}{8} \eta_{-\nu}(-T)^2.
\end{split}
\end{equation*}
Inequality~1 uses the eigenvalue $a_{-\nu}$ of $\A_{-T}$
and we {\color{gray}added zero}.
Inequality~2 pulls out the supremum, uses $\norm{\beta}_{L^\infty}\le 1$.
Inequality~3 uses $\delta<1$
and~(\ref{eq:777}).

\smallskip\noindent
\textbf{Term 4.}
By definition of $\zeta$ and $(v_\ell)$ being an orthonormal basis of $H_0$ we get
\begin{equation*}
\begin{split}
   \int_{-T}^{-T+\delta}\INNER{\eta(-T)}{\A_{-T}\zeta(s)}_{0} ds
   &\stackrel{1}{=}  \int_{-T}^{-T+\delta}
   \eta_{-\nu}(-T)^2\beta(s) a_{-\nu}
   ds\\
   &\le \delta \eta_{-\nu}(-T)^2 \Abs{a_{-\nu}} \\
   &\le \tfrac{1}{8} \eta_{-\nu}(-T)^2.
\end{split}
\end{equation*}
Equality~1 uses the eigenvalue $a_{-\nu}$ of $\A_{-T}$.
Then we use that $\norm{\beta}_{L^\infty}\le 1$ and the
final inequality exploits the choice of $\delta$.

\smallskip\noindent
The analysis of terms 1-4 shows
$
   \eta_{-\nu}(-T)^2
   \le 4\tfrac{\eps^2}{8} + 4 \tfrac{1}{8} \eta_{-\nu}(-T)^2
$.
Thus $\eta_{-\nu}(-T)^2\le\eps^2$ for every $\eps>0$. So $\eta_{-\nu}(-T)=0$.
This proves the claim.
\end{proof}

\medskip\noindent
\textbf{II. The inclusion 
\boldmath$\Kk \ker{D_{-A^*}^{+-}}\subset \Ann{(D_A^{+-})}$.}
Pick $\eta \in \ker{D_{-A^*}^{+-}}\subset P_1^{+-}(I_T,-\A^*_{\pm T};H_0^*,H_1^*)$.
For $\xi\in P_1^{+-}(I_T;H_1,H_0)$ we calculate
\begin{equation}\label{eq:eta-I_T-reverse}
\begin{split}
   &(\Kk\eta)D_A\xi
\\
   &=\int_{I_T} \eta(D_A\xi) \, ds
\\
   &=\int_{I_T}\eta(\p_s\xi)\, ds+\int_{I_T}\eta(A\xi)\, ds
\\
   &\stackrel{3}{=} -\int_{I_T}(\p_s\eta) \xi\, ds+\int_{I_T}(A^*\eta)\xi\, ds
\\
   &=(D_{-A^*}\eta)\xi
\\
   &=0 .
\end{split}
\end{equation}
Equation 1 is by definition of $\Kk$.
Equation 3 is integration by parts together with the fact that
$\eta$ and $\xi$ satisfy mutually orthogonal boundary conditions
at $-T$ as well as at $T$.
This proves that $\Kk \eta \in \Ann{(D_A^{+-})}$.

\smallskip
This concludes the proof of
Proposition~\ref{prop:coker-ker-I_T}.
\end{proof}

\boldmath
\subsubsection{Theorem~\ref{thm:A} -- Fredholm property}
\unboldmath

In order to prove Theorem~\ref{thm:A}
in the finite interval case, namely that $\mathfrak{D}_A$
is Fredholm and $\INDEX{\mathfrak{D}_A}=\varsigma(A)$,
we need Theorem~\ref{thm:proj}
(\hspace{-.1pt}\cite[Thm.\,D]{Frauenfelder:2024c}).

\begin{corollary}[to Theorem~\ref{thm:Fredholm-I_T-D_A+-}, Fredholm]
\label{cor:Fredholm-I_T}
For any $A\in \Aa_{I_T}^*$ the operator
\begin{equation*}
\begin{split}
   \mathfrak{D}_A=\mathfrak{D}_A^{I_T}\colon P_1(I_T)
   &\to P_0(I_T)\times H^+_{\frac12}(\A_{-T}) \times H^-_{\frac12}(\A_T)
   =: \Ww(I_T;\A_{-T},\A_T)
   \\
   \xi
   &\mapsto
   \left(D_A\xi,\pi_+^{\A_{-T}}\xi_{-T},\pi_-^{\A_{T}}\xi_T\right)
\end{split}
\end{equation*}
is Fredholm, where $\A_{\pm T}:= A(\pm T)$,
and of the same index as $D_A^{+-}$.
More precisely, the kernels coincide and the co-kernels are of equal dimension.
\end{corollary}

\begin{proof}
By Theorem~\ref{thm:finite-interval}
the operator $\mathfrak{D}_A$ is semi-Fredholm.
So it has finite dimensional kernel and closed image.
We shall show that kernel and image of $\mathfrak{D}_A$ are equal,
respectively isomorphic, to those of the Fredholm operator $D_A^{+-}$
from Theorem~\ref{thm:Fredholm-I_T-D_A+-}.

\medskip\noindent
\textbf{Step~1.} $\ker{\mathfrak{D}_A}=\ker{D_A^{+-}}$.

\begin{proof}
Clearly $\xi\in P_1(I_T)$ and $(0,0,0)=
\mathfrak{D}_A\xi:=(D_A\xi,\pi_+^{\A_{-T}}\xi_{-T},\pi_-^{\A_{T}}\xi_T)$
is equivalent to $D_A\xi=0$ and $\xi\in P_1^{+-}$; see~(\ref{eq:P_1^+-}).
\end{proof}

\medskip\noindent
\textbf{Step~2.} $D_A^{+-}$ is surjective iff $\mathfrak{D}_A$ is surjective.

\begin{proof}
``$\Rightarrow$"
Given $(\eta,x,y)\in \Ww(I_T;\A_{-T},\A_T)$, we need to find a
$\xi\in P_1$ such that $\mathfrak{D}_A\xi=(\eta,x,y)$.
To this end, pick $\xi_1\in P_1(I_T)$ which satisfies the
given boundary conditions $\xi_1(-T)=x$ and $\xi_1(T)=y$;
see Corollary~\ref{cor:Ev-onto}.
Since $D_A^{+-}$ is surjective, there exists $\xi_0\in P_1^{+-}$
such that
$
   D_A\xi_0=\eta-D_A\xi_1 \in P_0(I_T)
$.
We define $\xi:=\xi_0+\xi_1\in P_1(I_T)$.
Then $D_A\xi=\eta$ and
$\pi_+^{\A_{-T}}\xi(-T)=\pi_+^{\A_{-T}}\xi_0(-T)+\pi_+^{\A_{-T}}\xi_1(-T)
=0+\pi_+^{\A_{-T}}x=x$
and similarly $\pi_-^{\A_T}\xi(T)=y$.
Hence $\mathfrak{D}_A\xi=(\eta,x,y)$. 

``$\Leftarrow$"
Pick $\eta\in P_0(I_T)$ and consider $(\eta,0,0)\in \Ww(I_T;\A_{-T},\A_T)$.
Since $\mathfrak{D}_A$ is surjective there exists $\xi\in P_1(I_T)$ such that
$(\eta,0,0)=\mathfrak{D}_A\xi
=(D_A\xi,\pi_+^{\A_{-T}}\xi_{-T},\pi_-^{\A_{T}}\xi_T)$.
Since $\pi_+^{\A_{-T}}\xi_{-T}=0$ and $\pi_-^{\A_{T}}\xi_T=0$
we conclude that $\xi\in P_A^{+-}$ so that
$D_A^{+-}\xi=\eta$. This shows that $D_A^{+-}$ is surjective and
concludes the proof of Step~2.
\end{proof}

\medskip\noindent
\textbf{Step~3.} $\dim\coker{\mathfrak{D}_A}\le\dim\coker{D_A^{+-}}$
{\color{gray}\small$<\infty$ since $D_A^{+-}$ is Fredholm}

\begin{proof}
Suppose $n:=\dim\coker{D_A^{+-}}\ge 1$.
Let $\Bb=\{\beta_1,\dots,\beta_n\}$ be a basis of the orthogonal complement
$(\im{D_A^{+-}})^\perp$. Define an image filling operator by
$$
   \widetilde D_A^{+-}\colon P_1^{+-}\times\R^n\to P_0
   ,\quad
   (\xi,a)\mapsto D_A\xi+\sum_{i=1}^n a_i\beta_i
$$
and define a candidate to be image filling by
$$
   \widetilde{\mathfrak{D}}_A\colon P_1\times\R^n\to \Ww
   ,\quad
   (\xi,a)\mapsto \left(D_A\xi+\sum_{i=1}^n a_i\beta_i,
   \pi_+^{\A_{-T}}\xi(-T), \pi_-^{\A_T}\xi(T)\right) .
$$
We now show that $\widetilde{\mathfrak{D}}_A$ is surjective as well.
Pick $(\eta,x,y)\in \Ww(I_T;\A_{-T},\A_T)$.
We need to find $(\xi,a)\in  P_1\times\R^n$ such that
$\widetilde{\mathfrak{D}}_A(\xi,a)=(\eta,x,y)$.
To this end, pick $\xi_1\in P_1(I_T)$ which satisfies the
given boundary conditions $\xi_1(-T)=x$ and $\xi_1(T)=y$;
see Corollary~\ref{cor:Ev-onto}.
Since $\widetilde{D}_A^{+-}$ is surjective, there exists
$(\xi_0,a)\in P_1^{+-}\times\R^n$ such that
$$
   \widetilde{D}_A^{+-} (\xi,a)
   :=D_A\xi_0+\sum_{i=1}^n a_i\beta_i
   =\eta-D_A\xi_1 \in P_0(I_T) .
$$
This identity applied to $\xi:=\xi_0+\xi_1\in P_1(I_T)$ yields
\begin{equation*}
\begin{split}
   \widetilde{\mathfrak{D}}_A(\xi,a)
   =\left(D_A\xi+\sum_{i=1}^n a_i\beta_i,
   \pi_+^{\A_{-T}}\xi(-T), \pi_-^{\A_T}\xi(T)\right) 
   =\left(\eta,x,y\right) .
\end{split}
\end{equation*}
This proves that $\widetilde{\mathfrak{D}}_A$ is surjective.
Hence $\dim\coker{\mathfrak{D}_A}\le n{\color{gray}\,=\dim\coker{D_A^{+-}}}$.
This proves Step~3.
\end{proof}

\medskip\noindent
\textbf{Step~4.} $\dim\coker{D_A^{+-}}\le\dim\coker{\mathfrak{D}_A}$
{\color{gray}\small$<\infty$ by Step~3}

\begin{proof}
We choose a finite basis $\Bb=\{\beta_1,\dots,\beta_n\}\subset P_0$
of the orthogonal complement of the image of $D_A^{+-}$.
Suppose by contradiction that $n>\dim\coker{\mathfrak{D}_A}$.
Then the span of $\{(\beta_1,0,0),\dots, (\beta_n,0,0)\}\in \Ww$
has non-trivial intersection with $\im{\mathfrak{D}_A}\subset \Ww$.
Otherwise, the span would form a complement of $\im{\mathfrak{D}_A}$
and so $\dim\coker{\mathfrak{D}_A}\ge n$. Contradiction.
Hence some non-zero element in the span lies in the image of $\mathfrak{D}_A$,
in symbols $\exists a\in\R^n\setminus \{0\}$ such that
$$
   \left(\sum_{j=1}^n a_j \beta_j,0,0\right)
   =\sum_{j=1}^n a_j(\beta_j,0,0)
   \in \im{\mathfrak{D}_A} .
$$
Thus there exists $\xi\in P_1$ such that
$$
   \left(\sum_{j=1}^n a_j \beta_j,0,0\right)
   =\mathfrak{D}_A\xi
   =\left(D_A\xi,\pi_+^{\A_{-T}}\xi_{-T},\pi_-^{\A_{T}}\xi_T\right) .
$$
Hence $\xi\in P_1^{+-}$ and $D_A^{+-}\xi=\sum_{j=1}^n a_j \beta_j$.
Now the left hand side lies in the image of $D_A^{+-}$
and the right hand side in the orthogonal complement,
hence it is the zero vector. Since $\Bb$ is a basis all coefficients
$a_j$ are zero. Contradiction.
This proves Step~4.
\end{proof}

By Steps~2--4 the dimension of $\coker{D_A^{+-}}$ equals the one of
$\coker{\mathfrak{D}_A}$. Hence, by Step~1, the Fredholm indices are equal.
This proves Corollary~\ref{cor:Fredholm-I_T}.
\end{proof}

\subsubsection{Index and spectral content}

\begin{definition}\label{def:vhj44}
Given $A\in \Aa_{I_T}$, pick non-eigenvalues
$\lambda_{\pm T}\in\Rr (A(\pm T))$, set
$$
   \lambda:=(\lambda_{-T},\lambda_T)
   ,\quad
   \A_{-T}^{\lambda_{-T}}
   :=A(-T)-\lambda_{-T} \iota 
   ,\quad
   \A_T^{\lambda_T}
   :=A(T)-\lambda_T \iota ,
$$
then $\A_{-T}^{\lambda_{-T}}$ and $\A_T^{\lambda_T}$ lie in
$\Ll_{sym_0}^*(H_1,H_0)$.
Define moreover
$\mathfrak{D}_A^\lambda=\mathfrak{D}_A^{\lambda_{-T},\lambda_T}$ by
\begin{equation}\label{eq:frak-D-lambda}
\begin{split}
   \mathfrak{D}_A^\lambda\colon P_1(I_T)
   &\to
   P_0(I_T)\times H^+_{\frac12}(\A_{-T}^{\lambda_{-T}}) \times H^-_{\frac12}(\A_T^{\lambda_T})
   =: \Ww(I_T; \A_{-T}^{\lambda_{-T}},\A_T^{\lambda_T})
   \\
   \xi
   &\mapsto \left(D_A\xi\,,\,\pi_+^{\A_{-T}^{\lambda_{-T}}}\xi(-T)\,,\,
   \pi_-^{\A_T^{\lambda_T}}\xi(T)\right).
\end{split}
\end{equation}
\end{definition}
To compute the index difference of the Fredholm operator
$\mathfrak{D}_A^\lambda$ for different $\lambda$
we introduce the spectral content as follows.
For $a\in\spec{A}$ we denote by $E_a:=\ker{A-a\iota}$ the eigenspace
of $A$ to the eigenvalue $a$. Pick elements $\lambda\le\mu$ of $\Rr(A)$.
We define the \textbf{eigenspace interval}
$$
   E_{(\lambda,\mu)}
   :=\bigoplus_{a\in\spec{A}\atop a\in (\lambda,\mu)} E_a .
$$
The resulting decomposition defines
projections along $E_{(\lambda,\mu)}$, notation
\begin{equation}\label{eq:proj_interval_+}
   \pi^{A}_{(\lambda,\mu)}
   \colon 
   {\color{gray}\underbrace{{\color{black}H^+_\frac12(\A^{\lambda})}}
      _{H^{>\lambda}_\frac12(A)}}
   =
   E_{(\lambda,\mu)}
   \oplus 
   {\color{gray}\underbrace{{\color{black}H^+_\frac12(\A^{\mu})}}
      _{H_\frac12^{>\mu}(A)}}
   \to
   H^+_\frac12(\A^{\mu})
\end{equation}
and
\begin{equation}\label{eq:proj_interval_-}
   \pi^{A}_{(\mu,\lambda)}
   \colon 
   {\color{gray}\underbrace{{\color{black}H^-_\frac12(\A^{\mu})}}
      _{H^{<\mu}_\frac12(A)}}
   =
   {\color{gray}\underbrace{{\color{black}H^-_\frac12(\A^{\lambda})}}
      _{H_\frac12^{<\lambda}(A)}}
   \oplus E_{(\lambda,\mu)}
   \to H^-_\frac12(\A^{\lambda}) .
\end{equation}
We further define
the \textbf{spectral content} of $A$ 
between the two non-eigenvalues
$\lambda\le\mu$ as the number of eigenvalues in between, with
multiplicities, in symbols
\begin{equation}\label{eq:spectral-density}
   \rho_A(\lambda,\mu)
   :=\sum_{a\in\spec A\atop a\in(\lambda,\mu)}\dim{\ker{(A-a\iota)}}
   =\dim E_{(\lambda,\mu)} \in\N_0.
\end{equation}
Note that $\rho_A(\lambda,\lambda)=0$.
Moreover, we define
\begin{equation}\label{eq:spectral-density_-}
   \rho_A(\mu,\lambda)
   :=-\rho_A(\lambda,\mu) \in -\N_0.
\end{equation}
Note that due to the same summand $E_{(\lambda,\mu)}$
in~(\ref{eq:proj_interval_+}) and~(\ref{eq:proj_interval_-}) we have
\begin{equation}\label{eq:codim}
\begin{split}
   \codim\left(\text{$H^+_\frac12 (\A^\mu)$ in
   $H^+_\frac12 (\A^\lambda)$}\right)
   &=\rho_A(\lambda,\mu)
\\
   &=\codim\left(\text{$H^-_\frac12 (\A^\lambda)$ in
   $H^-_\frac12 (\A^\mu)$}\right) .
\end{split}
\end{equation}

\begin{lemma}\label{le:vhj44}
Given $A\in \Aa_{I_T}$, 
pick non-eigenvalues $\lambda_{-T},\mu_{-T}$ of $A(-T)$
and non-eigenvalues $\lambda_T,\mu_T$ of $A(T)$. Set
$\lambda:=(\lambda_{-T},\lambda_T)$ and $\mu:=(\mu_{-T},\mu_T)$.
Then
$$
   \INDEX{\mathfrak{D}_A^\mu}-\INDEX{\mathfrak{D}_A^\lambda}
   = \rho_{A(-T)}(\lambda_{-T},\mu_{-T})- \rho_{A(T)}(\lambda_T,\mu_T)
$$
where $\rho$ is the spectral content defined by~(\ref{eq:spectral-density}).
\end{lemma}

\begin{proof}
In the proof we distinguish four cases.

\medskip
\noindent
\textbf{Case~1.} $\lambda_{-T}\le\mu_{-T}$ and $\lambda_T\ge\mu_T$

\begin{proof}
Recall the projections defined
by~(\ref{eq:proj_interval_+}--\ref{eq:proj_interval_-}).
Since $\lambda_{-T}\le\mu_{-T}$ we have
\begin{equation*}
   \pi^{A(-T)}_{(\lambda_{-T},\mu_{-T})}
   \colon 
   H^+_\frac12(\A^{\lambda_{-T}}_{-T})
   =
   E_{(\lambda_{-T},\mu_{-T})}\oplus 
   H^+_\frac12(\A^{\mu_{-T}}_{-T})
  \to H^+_\frac12(\A^{\mu_{-T}}_{-T}) .
\end{equation*}
Since $\lambda_T\ge\mu_T$ we have
\begin{equation*}
   \pi^{A(T)}_{(\lambda_T,\mu_T)}
   \colon 
   H^-_\frac12(\A^{\lambda_T}_T)
   =
   H^-_\frac12(\A^{\mu_T}_T)
   \oplus E_{(\mu_T,\lambda_T)}
   \to H^-_\frac12(\A^{\mu_T}_T) .
\end{equation*}
Consider the projection defined by
$$
   p=\left(\1\,,\, \pi^{A(-T)}_{(\lambda_{-T},\mu_{-T})} \,,\, \pi^{A(T)}_{(\lambda_T,\mu_T)} \right)
   \colon
   \Ww(I_T; \A_{-T}^{\lambda_{-T}},\A_T^{\lambda_T})
   \to
   \Ww(I_T; \A_{-T}^{\mu_{-T}},\A_T^{\mu_T}) 
$$
and observe that $\mathfrak{D}_A^\mu=p\circ \mathfrak{D}_A^\lambda$.
By Theorem~\ref{thm:Fredholm-compose} we get identity 1 in
\begin{equation*}
\begin{split}
   &\INDEX{\mathfrak{D}_A^\mu}-\INDEX{\mathfrak{D}_A^\lambda}
\\
   &\stackrel{1}{=}
   \codim\left(\text{$\Ww(I_T; \A_{-T}^{\lambda_{-T}},\A_T^{\lambda_T})$ in
   $\Ww(I_T; \A_{-T}^{\mu_{-T}},\A_T^{\mu_T})$}\right)
\\
   &\stackrel{2}{=}
   \codim\left(\text{$H^+_\frac12 (\A_{-T}^{\lambda_{-T}})$ in
   $H^+_\frac12 (\A_{-T}^{\mu_{-T}})$}\right)
   \\
   &\quad
   +\codim\left(\text{$H^-_\frac12 (\A_T^{\lambda_T})$ in
   $H^-_\frac12 (\A_T^{\mu_T})$}\right)
\\
   &\stackrel{3}{=}
   \rho_{A(-T)}(\lambda_{-T},\mu_{-T})+\rho_{A(T)}(\mu_T,\lambda_T)
\\
   &=\rho_{A(-T)}(\lambda_{-T},\mu_{-T})-\rho_{A(T)}(\lambda_T,\mu_T) .
\end{split}
\end{equation*}
Identity 2 holds by the inclusion
$\Ww(I_T; \A_{-T}^{\mu_{-T}},\A_T^{\mu_T})\subset
\Ww(I_T; \A_{-T}^{\lambda_{-T}},\A_T^{\lambda_T})$
due to the inclusions of the second and third factors.
Identity 3 holds by~(\ref{eq:codim}).
\end{proof}

\medskip
\noindent
\textbf{Case~2.} $\lambda_{-T}\ge\mu_{-T}$ and $\lambda_T\le\mu_T$

\begin{proof}
Interchanging the roles of $\lambda$ and $\mu$ in Case~1 we obtain
$$
   \INDEX{\mathfrak{D}_A^\lambda}-\INDEX{\mathfrak{D}_A^\mu}
   = \rho_{A(-T)}(\mu_{-T},\lambda_{-T})- \rho_{A(T)}(\mu_T,\lambda_T) .
$$
Taking the negative of both sides we obtain
\begin{equation*}
\begin{split}
   \INDEX{\mathfrak{D}_A^\mu}-\INDEX{\mathfrak{D}_A^\lambda}
   &= -\rho_{A(-T)}(\mu_{-T},\lambda_{-T}) +\rho_{A(T)}(\mu_T,\lambda_T)
\\
   &=\rho_{A(-T)}(\lambda_{-T},\mu_{-T}) - \rho_{A(T)}(\lambda_T,\mu_T)
\end{split}
\end{equation*}
where the second identity is by~(\ref{eq:spectral-density_-}).
\end{proof}

\medskip
\noindent
\textbf{Case~3.} $\lambda_{-T}\le\mu_{-T}$ and $\lambda_T\le\mu_T$

\begin{proof}
Add zero to obtain
\begin{equation*}
\begin{split}
   &\INDEX{\mathfrak{D}_A^\mu}-\INDEX{\mathfrak{D}_A^\lambda}
\\
   &=\INDEX{\mathfrak{D}_A^{\mu_{-T},\mu_T}}
   -\INDEX{\mathfrak{D}_A^{\lambda_{-T},\mu_T}}
   +\INDEX{\mathfrak{D}_A^{\lambda_{-T},\mu_T}}
   -\INDEX{\mathfrak{D}_A^{\lambda_{-T},\lambda_T}}
\\
   &\stackrel{2}{=}
   \rho_{A(-T)}(\lambda_{-T},\mu_{-T})-\rho_{A(T)}(\mu_T,\mu_T) 
   +\rho_{A(-T)}(\lambda_{-T},\lambda_{-T}) - \rho_{A(T)}(\lambda_T,\mu_T)
\\
   &=\rho_{A(-T)}(\lambda_{-T},\mu_{-T}) - \rho_{A(T)}(\lambda_T,\mu_T)
\end{split}
\end{equation*}
where in identity 2 we used Case~1 for the first difference
and Case~2 for the second one.
\end{proof}

\medskip
\noindent
\textbf{Case~4.} $\lambda_{-T}\ge\mu_{-T}$ and $\lambda_T\ge\mu_T$

\begin{proof}
This follows by literally the same computation
as in Case~3. What differs is the explanation: now
in identity 2 we use Case~2 for the first difference
and Case~1 for the second one.
\end{proof}

This proves Lemma~\ref{le:vhj44}.
\end{proof}

\subsubsection{Path concatenation}\label{sec:concatenation}

Let $A\in\Aa_{I_T}^*$ be a Hessian path such that not only
$\A_{-T}:=A(-T)$ and $\A_T:=A(T)$ are invertible, but also the
Hessian operator at time zero $\A_0:=A(0)$ is.
Decomposing the time interval at time zero
$$
   I_T:=[-T,T]=I_T^-\cup I_T^+
   ,\qquad
   I_T^-:=[-T,0],\quad I_T^+:=[0,T] ,
$$
gives rise to three augmented operators, one along each of the
three intervals.
Firstly, there is the operator along $I_T$, defined by
\begin{equation*}
\begin{split}
   \mathfrak{D}_A
   \colon P_1(I_T)
   &\to  P_0(I_T)\times H^+_{\frac12}(\A_{-T}) \times H^-_{\frac12}(\A_T)
   {\color{gray}\;=: \Ww(I_T;\A_{-T},\A_T)}
\\
   \xi
   &\mapsto \left(D_A\xi,\pi_+^{\A_{-T}}\xi_{-T},\pi_-^{\A_{T}}\xi_T\right)
\end{split}
\end{equation*}
where $\xi_{\pm T}:=\xi(\pm T)$.
We define the operator along $I_T^-=[-T,0]$ by
\begin{equation*}
\begin{split}
   \mathfrak{D}_{A|_{[-T,0]}}
   \colon P_1(I_T^-)
   &\to
   P_0(I_T^-)\times H^+_{\frac12}(\A_{-T}) \times H^-_{\frac12}(\A_0)
   {\color{gray}\;=: \Ww(I_T^-; \A_{-T},\A_0)}
   \\
   \xi
   &\mapsto \left(D_A\xi\,,\,\pi_+^{\A_{-T}}\xi_{-T}\,,\,
   \pi_-^{\A_0}\xi_0\right)
\end{split}
\end{equation*}
and the operator along $I_T^+=[0,T]$ is defined by
\begin{equation*}
\begin{split}
   \mathfrak{D}_{A|_{[0,T]}}
   \colon P_1(I_T^+)
   &\to
   P_0(I_T^+)\times H^+_{\frac12}(\A_0) \times H^-_{\frac12}(\A_T)
   {\color{gray}\;=:\Ww(I_T^+; \A_0,\A_T)}
   \\
   \xi
   &\mapsto \left(D_A\xi\,,\,\pi_+^{\A_{0}}\xi_0\,,\,
   \pi_-^{\A_T}\xi_T\right).
\end{split}
\end{equation*}

\begin{theorem}[Path concatenation]\label{thm:concatenation}   
Suppose $A\in\Aa_{I_T}^*$ is such that $\A_0:=A(0)$ is invertible.
Then the Fredholm index is additive under concatenation
$$
   \INDEX{\mathfrak{D}_A}
   =\INDEX{\mathfrak{D}_{A|_{[-T,0]}}}
   +\INDEX{\mathfrak{D}_{A|_{[0,T]}}} .
$$
\end{theorem}

A main proof ingredient is a domain-homotopy of Fredholm operators.
For paths $\xi\in P_1(I_T^-)$ and $\eta\in P_1(I_T^+)$ we abbreviate
\begin{align}\label{eq:xi_=-}
   \xi_\pm(0)
   &:=\pi_\pm^{\A_0}\xi(0),
   &
   \eta_\pm(0)
   &:=\pi_\pm^{\A_0}\eta(0) .
\end{align}
For $r\in[0,1]$ we define a family of spaces by
$$
   \mathfrak{P}_r
   :=\{(\xi,\eta)\in P_1(I_T^-)\times P_1(I_T^+)
   \mid \xi_-(0)=r \eta_-(0)\;\wedge\; r\xi_+(0)=\eta_+(0)\} .
$$

\begin{proposition}\label{prop:concat-r}
Consider the family of operators defined, for $r\in[0,1]$, by
\begin{equation}\label{eq:D_A,r}
\begin{split}
   \mathfrak{D}_{A,r}
   \colon 
   X\supset
   \mathfrak{P}_r
   &\to  P_0(I_T)\times H^+_{\frac12}(\A_{-T}) \times H^-_{\frac12}(\A_T) =: Y
\\
   (\xi,\eta)
   &\mapsto \left( (D_A\xi)\#(D_A\eta),\pi_+^{\A_{-T}}\xi_{-T},\pi_-^{\A_{T}}\eta_T\right)
\end{split}
\end{equation}
where $X=P_1(I_T^-)\times P_1(I_T^+)$ and
$$
   \left(D_A\xi\right)\#\left(D_A\eta\right)\, (s)
   :=
   \begin{cases}
      D_A\xi(s) &\text{, $s\in [-T,0)$,}
      \\
      D_A\eta(s)&\text{, $s\in [0,T]$.}
   \end{cases}
$$
Then the following is true.
a)~Each member of the family is a Fredholm operator and
b)~the Fredholm index is constant along the family.
\end{proposition}

\begin{proof}[Proof of Theorem~\ref{thm:concatenation}]
Let us present right away the proof in a nutshell
\begin{align*}
   \INDEX{\mathfrak{D}_A}
   &\stackrel{{\color{gray}1.}}{=}\INDEX{\mathfrak{D}_{A,1}}
   &&\text{\color{gray}\small , equal operators (Step~1)}&
\\
   &\stackrel{{\color{gray}2.}}{=}\INDEX{\mathfrak{D}_{A,0}}
   &&\text{\color{gray}\small , homotopy, Prop.\,\ref{prop:concat-r}}&
\\
   &\stackrel{{\color{gray}3.}}{=}
   \INDEX{\left(\mathfrak{D}_{A|_{[-T,0]}}^{\bullet -}
   \oplus \mathfrak{D}_{A|_{[0,T]}}^{+ \bullet}\right)}
   &&\text{\color{gray}\small , decompose $P_0(I_T)$ (Step~3)}&
\\
   &\stackrel{{\color{gray}4.}}{=}
   \INDEX{\mathfrak{D}_{A|_{[-T,0]}}^{\bullet -}}
   +\INDEX{\mathfrak{D}_{A|_{[0,T]}}^{+ \bullet}}
   &&\text{\color{gray}\small , direct sum (obvious)}&
\\
   &\stackrel{{\color{gray}5.}}{=}
   \INDEX{\mathfrak{D}_{A|_{[-T,0]}}}+\INDEX{\mathfrak{D}_{A|_{[0,T]}}}
   &&\text{\color{gray}\small , summand-wise equality.}&
\end{align*}
Before filling in the details of steps~1-5
we need to define the $+-$ operators appearing individually
in step~3 and as direct sum in step~4. The operators\footnote{
  \textit{Bullet notation.}
  An interval has two boundary points, left and right.
  The bullet '$\bullet$' symbolizes 'no boundary condition'
  at that boundary point whose position corresponds to the position of the bullet,
  left or right.
  The sign $+/-$ tells that the positive/negative part
  of the spectrum must vanish at the other boundary point.
}
\begin{equation*}
\begin{split}
   \mathfrak{D}_{A|_{[-T,0]}}^{\bullet -}
   \colon P_1^{\bullet -}(I_T^-) 
   &\to P_0(I_T^-)\times H^+_{\frac12}(\A_{-T}) 
   {\color{gray}\; =:\Ww(I_T^-; \A_{-T})}
\\
   \xi
   &\mapsto \left(D_A\xi,\pi_+^{\A_{-T}}\xi_{-T}\right)
\end{split}
\end{equation*}
and
\begin{equation*}
\begin{split}
   \mathfrak{D}_{A|_{[0,T]}}^{+ \bullet}
   \colon P_1^{+ \bullet}(I_T^+)
   &\to P_0(I_T^+)\times H^-_{\frac12}(\A_T)
   {\color{gray}\; =:\Ww(I_T^+;\A_T)}
\\
   \eta
   &\mapsto \left(D_A\eta,\pi_-^{\A_T}\eta_T\right)
\end{split}
\end{equation*}
are defined by the usual formula $\p_s+A(s)$.
Observe that one boundary condition is imposed on the domain and the
other one on the co-domain as follows
\begin{equation*}
\begin{split}
   P_1^{\bullet -}(I_T^-) 
   &:=\{\xi\in P_1(I_T^-)\mid \xi_-(0)\stackrel{\text{(\ref{eq:xi_=-})}}{=}0\},
\\
   P_1^{+ \bullet}(I_T^+) 
   &:=\{\eta\in P_1(I_T^+)\mid \eta_+(0)\stackrel{\text{(\ref{eq:xi_=-})}}{=}0\}.
\end{split}
\end{equation*}

The proof proceeds in six steps 0-5 as enumerated in the nutshell.

\smallskip
\noindent
\textbf{Step~0.} 
The operators appearing in the nutshell are all Fredholm.

\begin{proof}
The operators $\mathfrak{D}_A$,  $\mathfrak{D}_{A|_{[-T,0]}}$, and
$\mathfrak{D}_{A|_{[0,T]}}$ are Fredholm, by Corollary~\ref{cor:Fredholm-I_T},
and $\mathfrak{D}_{A,1}$ and $\mathfrak{D}_{A,0}$
are Fredholm, by Proposition~\ref{prop:concat-r}.
The operators $\mathfrak{D}_{A|_{[-T,0]}}^{\bullet +}$ and
$\mathfrak{D}_{A|_{[0,T]}}^{- \bullet}$ are Fredholm
by the arguments in the proof of Corollary~\ref{cor:Fredholm-I_T}.
This concludes Step~0.
\end{proof}

\smallskip
\noindent
\textbf{Step~1.}
$\mathfrak{D}_A=\mathfrak{D}_{A,1}$

\begin{proof}
This follows immediately using the obvious identification of
$P_1(I_T)$ with
$$
   \mathfrak{P}_1
   =\{(\xi,\eta)\in P_1(I_T^-)\times P_1(I_T^+)
   \mid \xi_-(0)=\eta_-(0)\;\wedge\; \xi_+(0)=\eta_+(0)\} .
$$
by cutting the elements of $P_1(I_T)$ at time $0$
into two pieces. See also Appendix~\ref{sec:evaluation} on the
evaluation map.
\end{proof}

\smallskip
\noindent
\textbf{Step~2.} 
Homotopy of Fredholm operators between $\mathfrak{D}_{A,0}$ and
$\mathfrak{D}_{A,1}$.

\begin{proof}
Proposition~\ref{prop:concat-r}.
\end{proof}

\smallskip
\noindent
\textbf{Step~3.} 
$\mathfrak{D}_{A,0}$
and
$\mathfrak{D}_{A|_{[-T,0]}}^+\oplus\mathfrak{D}_{A|_{[0,T]}}^-$
correspond naturally

\begin{proof}
This follows from equal domain
$$
   \mathfrak{P}_0
   :=\{(\xi,\eta)\in P_1(I_T^-)\times P_1(I_T^+)
   \mid \xi_-(0)=0=\eta_+(0)\}
   =P_1^{\bullet -}(I_T^-) \times P_1^{+ \bullet}(I_T^+)
$$
and, in the co-domain, the identification of $P_0(I_T)$ with
$P_0(I_T^-)\times P_0(I_T^+)$.
\end{proof}

\smallskip
\noindent
\textbf{Step~4.} 
Direct sum of Fredholm operators is Fredholm of index the index sum.

\begin{proof}
Well known.
\end{proof}

\smallskip
\noindent
\textbf{Step~5.} 
$\INDEX{\mathfrak{D}_{A|_{[-T,0]}}^{\bullet -}}=\INDEX{\mathfrak{D}_{A|_{[-T,0]}}}$
  {\color{gray}and}
  $\INDEX{\mathfrak{D}_{A|_{[0,T]}}^{+ \bullet}}=\INDEX{\mathfrak{D}_{A|_{[0,T]}}}$

\begin{proof}
This follows by the arguments in the proof of
Corollary~\ref{cor:Fredholm-I_T}.
\end{proof}

This concludes the proof of Theorem~\ref{thm:concatenation}.
\end{proof}

\begin{proof}[Proof of Proposition~\ref{prop:concat-r}]
We define, for each $r\in[0,1]$, a bounded linear map
\begin{equation*}
\begin{split}
   F_r\colon X:=
   P_1(I_T^-)\times P_1(I_T^+)
   &\to
   H^-_{\frac12}(\A_0)\times H^+_{\frac12}(\A_0) =: Z
\\
   (\xi,\eta)
   &\mapsto
   \left( \xi_-(0)-r\eta_-(0) , r\xi_+(0)-\eta_+(0) \right) .
\end{split}
\end{equation*}
The kernel is the domain $\mathfrak{P}_r=\ker{F_r}$ of the restriction
$\mathfrak{D}_{A,r}\colon \mathfrak{P}_r\to Y$ of
$$
   \mathfrak{D}_A
   \colon
   X
   \to
   Y
   ,\quad
   (\xi,\eta)
   \mapsto \left( (D_A\xi)\#(D_A\eta),\pi_+^{\A_{-T}}\xi_{-T},\pi_-^{\A_{T}}\xi_T\right).
$$

\smallskip
\noindent
a)
$\mathfrak{D}_{A,r}\colon\mathfrak{P}_r\to Y$ is Fredholm $\forall r\in[0,1]$.

\smallskip
\noindent
By Proposition~\ref{prop:Sim-r} each operator
$\mathfrak{D}_{A,r}\colon\mathfrak{P}_r\to Y$ is semi-Fredholm.
Theorem~\ref{thm:Fred-var-domain}
for $D=\mathfrak{D}_A$ and $D_r=\mathfrak{D}_{A,r}$
therefore implies that its semi-Fredholm index is independent of
$r\in[0,1]$.
By Step~1 in the proof of Theorem~\ref{thm:concatenation}
the operator $\mathfrak{D}_{A,1}$ is Fredholm and therefore
has finite index. Hence, by independence of $r$, every $\mathfrak{D}_{A,r}$
has finite index and is therefore Fredholm.

\smallskip
\noindent
b)
$\INDEX{\mathfrak{D}_{A,r}}$ does not depend on $r$.

\smallskip
\noindent
Use Theorem~\ref{thm:Fred-var-domain}
for $D=\mathfrak{D}_A$ and $D_r=\mathfrak{D}_{A,r}$.
This proves Proposition~\ref{prop:concat-r}.
\end{proof}

\begin{proposition}\label{prop:Sim-r}
Let $A\in\Aa_{I_T}^*$.
Then there is a constant $c>0$ such that
\begin{equation*}
\begin{split}
   \norm{\xi}_{P_1(I_T^-)}+\norm{\eta}_{P_1(I_T^+)}
   &\le c\Bigl(
   \norm{\xi}_{P_0(I_T^-)}+\norm{\eta}_{P_0(I_T^+)}
   +\norm{D_\A\xi}_{P_0(I_T^-)}+\norm{D_\A\eta}_{P_0(I_T^+)}
   \\
   &\qquad
   +\norm{\xi_+(-T)}_{\frac12}+\norm{\eta_-(T)}_{\frac12} .
   \Bigr)
\end{split}
\end{equation*}
for every $(\xi,\eta)\in\mathfrak{P}_r$ and $r\in[0,1]$.
\end{proposition}

\begin{proof}
The proof follows the same way as the proof
of Theorem~\ref{thm:finite-interval}.
The only step which needs to be adjusted is the
estimate~(\ref{eq:Step-1-Sim}) in step 1. For this adjustment
the assumption that $r\in[0,1]$ is crucial.
Therefore it suffices to show for any constant path
$A(s)\equiv\A {\color{gray}\; =\A_T=\A_{-T}}$
consisting of an invertible operator
existence of a constant $c>0$ such that
\begin{equation}\label{eq:Step-1-Sim-r}
\begin{split}
   \norm{\xi}_{P_1(I_T^-)}+\norm{\eta}_{P_1(I_T^+)}
   &\le c\Bigl(
   \norm{D_\A\xi}_{P_0(I_T^-)}+\norm{D_\A\eta}_{P_0(I_T^+)}
   \\
   &\qquad
   +\norm{\xi_+(-T)}_{\frac12}+\norm{\eta_-(T)}_{\frac12} .
   \Bigr)
\end{split}
\end{equation}
for every $(\xi,\eta)\in\mathfrak{P}_r$.

To see this we proceed as follows.
As in Step~1 in the proof of Theorem~\ref{thm:finite-interval},
by changing the constant $C_1$ if necessary, we can assume without
loss of generality, as explained in Section~\ref{sec:int-spaces-1/2}, that
$\A\colon H_1\to H_0$ is a symmetric isometry.
We abbreviate $\pi_\pm:=\pi_\pm^\A$.
Analogous to~(\ref{eq:gh654-NEW}) it holds that
\begin{equation}\label{eq:gh654-xi-r}
\begin{split}
   \Norm{D_{\A}\xi}_{P_0(I_T^-)}^2
   &\stackrel{2}{=} \Norm{\A\xi}_{P_0(I_T^-)}^2+\Norm{\p_s\xi}_{P_0(I_T^-)}^2
   +\INNER{\xi_0}{\A\xi_0}_0-\INNER{\xi_{-T}}{\A\xi_{-T}}_0 
\end{split}
\end{equation}
for every $\xi\in P_1(I_T^-)$ and that
\begin{equation}\label{eq:gh654-eta-r}
\begin{split}
   \Norm{D_{\A}\eta}_{P_0(I_T^+)}^2
   &= \Norm{\A\eta}_{P_0(I_T^+)}^2+\Norm{\p_s\eta}_{P_0(I_T^+)}^2
   +\INNER{\eta_T}{\A\eta_T}_0-\INNER{\eta_0}{\A\eta_0}_0 
\end{split}
\end{equation}
for any $\eta\in P_1(I_T^+)$. The opposite signs are crucial.
As $\A$ is an isometry we get
{\small
\begin{equation*}
\begin{split}
   \INNER{\xi(0)}{\A\xi(0)}_0
   =\INNER{\xi_+(0)}{\A \xi_+(0)}_0
   +\INNER{\xi_-(0)}{\A \xi_-(0)}_0
   &=\norm{\xi_+(0)}_{\frac12}^2-\norm{\xi_-(0)}_{\frac12}^2,
\\
   \INNER{\eta(0)}{\A\eta(0)}_0
   =\INNER{\eta_+(0)}{\A \eta_+(0)}_0
   +\INNER{\eta_-(0)}{\A \eta_-(0)}_0
   &=\norm{\eta_+(0)}_{\frac12}^2-\norm{\eta_-(0)}_{\frac12}^2 .
\end{split}
\end{equation*}
}Taking the difference and using the relations in $\mathfrak{P}_r$ we obtain
{\small
\begin{equation}\label{eq:0er}
\begin{split}
   &\INNER{\xi(0)}{\A\xi(0)}_0-\INNER{\eta(0)}{\A\eta(0)}_0
   \\
   &=\norm{\xi_+(0)}_{\frac12}^2-\norm{\xi_-(0)}_{\frac12}^2
   -\norm{\eta_+(0)}_{\frac12}^2+\norm{\eta_-(0)}_{\frac12}^2
   \\
   &=\norm{\xi_+(0)}_{\frac12}^2-r^2\norm{\eta_-(0)}_{\frac12}^2
   -r^2\norm{\xi_+(0)}_{\frac12}^2+\norm{\eta_-(0)}_{\frac12}^2
   \\
   &=(1-r^2)\left(\norm{\xi_+(0)}_{\frac12}^2+\norm{\eta_-(0)}_{\frac12}^2\right)
   \\
   &\ge 0
\end{split}
\end{equation}
}where we used that $r\in[0,1]$.
By definition~(\ref{eq:P_1-inner}) of the $P_1$ norm we have
\begin{equation*}
\begin{split}
   &\norm{\xi}_{P_1(I_T^-)}^2+\norm{\eta}_{P_1(I_T^+)}^2\\
   &=
   \norm{\A\xi}_{P_0(I_T^-)}^2+\norm{\p_s\xi}_{P_0(I_T^-)}^2
   +\norm{\A\eta}_{P_0(I_T^+)}^2+\norm{\p_s\eta}_{P_0(I_T^+)}^2
   \\
   &\stackrel{2}{=}
   \Norm{D_{\A}\xi}_{P_0(I_T^-)}^2
   +\underline{\INNER{\xi_{-T}}{\A\xi_{-T}}_0}
   -\INNER{\xi_0}{\A\xi_0}_0
   \\
   &\quad
   +\Norm{D_{\A}\eta}_{P_0(I_T^+)}^2
   +\INNER{\eta_0}{\A\eta_0}_0
   -\underline{\INNER{\eta_T}{\A\eta_T}_0}
   \\
   &\stackrel{3}{\le}
   \Norm{D_{\A}\xi}_{P_0(I_T^-)}^2
   +\Norm{\xi_+(-T)}_{H_{\frac12}}^2
   +\Norm{D_{\A}\eta}_{P_0(I_T^+)}^2
   +\Norm{\eta_+(-T)}_{H_{\frac12}}^2 .
\end{split}
\end{equation*}
Equality 2 uses the two displayed identities after~(\ref{eq:Step-1-Sim-r}).
In inequality 3 we used~(\ref{eq:0er}) to drop the terms at time $s=0$
and we used the estimates~(\ref{eq:vghvhj95478})
on the underlined terms.
This proves~(\ref{eq:Step-1-Sim-r}) and
Proposition~\ref{prop:Sim-r}.
\end{proof}

\subsubsection{Theorem~\ref{thm:A} -- Index formula}

In order to prove the assertion
$\INDEX{\mathfrak{D}_A}=\varsigma(A)$
of Theorem~\ref{thm:A} we utilize
%
\begin{theorem}[\hspace{-.1pt}{\cite[Thm.\,D]{Frauenfelder:2024c}}]
\label{thm:proj}
For a Hilbert space pair $(H_0,H_1)$ the maps
\begin{equation}\label{eq:proj}
    \pi_\pm\colon \Ll_{sym_0}^*(H_1,H_0)\to \Ll(H_\frac{1}{2})
   ,\quad
   \A\mapsto \pi_\pm^\A
\end{equation}
are continuous.\footnote{
  $\Ll_{sym_0}^*(H_1,H_0)$ consists of the 
  invertible $A\in\Ll(H_1,H_0)$ which are
  $H_0$-symmetric~(\ref{eq:H_0-symmetric-INTRO}).
  }
\end{theorem}

\begin{proof}[Proof of Theorem~\ref{thm:A} -- Index formula]\mbox{}\\
Let $A\in \Aa^*_{I_T}:=
\{A\in \Aa_{I_T}\mid \text{$A(-T)$ and $A(T)$ are invertible} \}$.
We write $\A(-T)$ and $\A(T)$ to indicate invertibility.

\medskip
\noindent
\textbf{Step~1.} Theorem~\ref{thm:A} holds if every operator $\A(s)$
in the path $\A$ is invertible.

\begin{proof}
Homotop to constant invertible $\A(0)$.
Consider the homotopy of paths $\A_r(s):=\A(rs)$ for $r\in[0,1]$.
Then the initial path $\A_0\equiv \A(0)$ is constant and invertible and
the end path $\A_1=\A$ is the given path.
We claim the identity
\begin{equation}\label{eq:claim-index}
   \INDEX{\mathfrak{D}_\A}=\INDEX{\mathfrak{D}_{\A(0)}} .
\end{equation}
The proof uses Theorem~\ref{thm:Fredholm-index-constant}.
To homotopy member $\A_r$ we assign the operator
\begin{equation*}\label{eq:augmented-T-7}
\begin{split}
   \Dd_r\colon P_1(I_T)
   &\to P_0(I_T)\times H_{\frac12} \times H_{\frac12}
   \\
   \xi
   &\mapsto \left(D_{\A_r}\xi,\xi(-T),\xi(T)\right)
\end{split}
\end{equation*}
and the projection
\begin{equation*}\label{eq:augmented-P-7}
\begin{split}
   p_r=\left(\1,\pi_+^{\A_r(-T)},\pi_-^{\A_r(T)}\right)\colon
   P_0(I_T)\times H_{\frac12} \times H_{\frac12}
   &\to 
   P_0(I_T)\times H_{\frac12} \times H_{\frac12} .
\end{split}
\end{equation*}
Observe that $\A(\pm rT)=\A_r(\pm T)$.
Composing both operators we get
\begin{equation*}\label{eq:augmentedfrak-7}
\begin{split}
   \mathfrak{D}_{A_r}:=p_r\circ \Dd_r\colon
   P_1(I_T)
   &\to P_0(I_T)\times H^+_{\frac12}(\A(-rT)) \times H^-_{\frac12}(\A(rT)) .
   \\
   \xi
   &\mapsto \left(D_{\A_r}\xi,\pi_+^{\A(-rT)}\xi(-T),\pi_-^{\A(rT)}\xi(T)\right)
\end{split}
\end{equation*}
The projections $p_r$ depend continuously on $r$
in view of Theorem~\ref{thm:proj}.
Therefore Theorem~\ref{thm:Fredholm-index-constant} implies that
$$
   \INDEX{\mathfrak{D}_{\A_0}}=\INDEX{\mathfrak{D}_{\A_1}} .
$$
Since $\mathfrak{D}_{\A_0}=\mathfrak{D}_{\A(0)}$ and
$\mathfrak{D}_{\A_1}=\mathfrak{D}_{\A}$ this proves
the claimed identity~(\ref{eq:claim-index}).

By~(\ref{eq:augmented-T-AA}) in Step~1 of the proof of
Theorem~\ref{thm:finite-interval},
the operator $\mathfrak{D}_{\A(0)}$ is an isomorphism
and therefore a Fredholm operator of index zero.
Hence, in view of~(\ref {eq:claim-index}), we have
$\INDEX{\mathfrak{D}_{\A}}=0$.
Since $\A(s)$ is invertible for every $s$,
the spectral flow $\varsigma(\A)$ is zero.
This proves Theorem~\ref{thm:A}
in case of a family of invertible operators
along a finite interval $I_T$.
This proves Step~1.
\end{proof}

\smallskip
\noindent
\textbf{Step~2.} Let $A\in \Aa^*_{I_T}$.
There exists an integer $N\in\N_0$ and real numbers
$$
   {\color{gray} -T=\,}t_0<t_1<\dots <t_N<t_{N+1}{\color{gray}\,=T ,}
   \qquad
   {\color{gray} 0=\,}\lambda_0,\lambda_1,\dots,\lambda_{N-1},
   \lambda_N
$$
such that
$$
   \A_j(s)
   :=A(s)-\lambda_j\iota \colon H_1\to H_0
$$
is invertible for $s\in[t_j,t_{j+1}]$ whenever $j\in\{0,\dots,N\}$.
\begin{figure}[h]
  \centering
  \includegraphics
                             {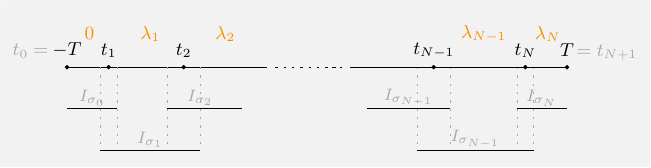}
  \caption{Step~2: Invertibility shifts $\lambda_i$ and intervals $[t_j,t_{j+1}]$}
   \label{fig:fig-I_sigma}
\end{figure} 

\begin{proof}
For each $\sigma\in[-T,T]$ choose
$\mu_\sigma\in\R\setminus \spec A(\sigma)$.
So $A(\sigma)-\mu_\sigma\iota\colon H_1\to H_0$  is invertible.
Since invertibility is an open condition, there exists $\eps_\sigma>0$
such that $A(\tau)-\mu_\sigma\iota$ is invertible for every
$$
   \tau\in
   I_\sigma:=(\sigma-\eps_\sigma,\sigma+\eps_\sigma)\cap[-T,T]
   .
$$
Since $\A(-T)$ is invertible we choose 
$$
   \mu_{-T}=0 .
$$
Since $[-T,T]$ is compact, there exists a finite subset $\mathfrak{S}$
of $[-T,T]$ such that the corresponding open intervals still cover $[-T,T]$,
in symbols
$$
   \bigcup_{\sigma\in \mathfrak{S}} I_\sigma
   =[-T,T] .]
   ,\qquad
   \Ii:=\{I_\sigma\mid\sigma\in \mathfrak{S}\}
   .
$$
We can assume without loss of generality that $-T,T\in \mathfrak{S}$,
otherwise just add two intervals.

Out of this finite covering we construct recursively
a further sub-covering $\Ii:=\{I_{\sigma_0}, I_{\sigma_1},\dots,
I_{\sigma_N}\}$ beginning at $\sigma_0:=-T$
and such that exactly nearest neighbors overlap.
If $T\in I_{\sigma_j}$, we set $N:=j$ and we are done.
If $T\notin I_{\sigma_j}$, then we choose
$\sigma_{j+1}\in \mathfrak{S}$ satisfying the two conditions
  \begin{enumerate}\setlength\itemsep{0ex} 
    \item
    $I_{\sigma_{j+1}}\cap I_{\sigma_j}\not=\emptyset$
    \hfill {$\scriptstyle\color{gray} \text{intersects predecessor $j$}$}
  \item
    $\sigma_{j+1}+\eps_{\sigma_{j+1}}\ge \sigma+\eps_\sigma$,
    $\forall \sigma\in \mathfrak{S}\colon
    I_{\sigma}\cap I_{\sigma_j}\not=\emptyset$
    \hfill {$\scriptstyle\color{gray} \text{farthest right among intersectors}$}
  \end{enumerate}
Condition 1 means that the chosen interval $I_{\sigma_{j+1}}$ intersects its predecessor.
Condition 2 means that the chosen interval $I_{\sigma_{j+1}}$ reaches farthest to the
right among all intersectors.
Furthermore, there are the following consequences
\begin{enumerate}
\item[\rm (i)]
  $\sigma_{j+1}+\eps_{\sigma_{j+1}}> \sigma_j+\eps_{\sigma_j}$;
    \hfill {$\scriptstyle\color{gray} \text{successor $j+1$ extends further right}$}
\item[\rm (ii)]
  If $I_{\sigma_i}\cap I_{\sigma_j}\not=\emptyset$ where
  $i,j\in\{1,\dots,N\}$.
  then $\abs{i-j}\le 1$.
  \\\mbox{}
    \hfill {$\scriptstyle\color{gray} \text{only next neighbors can intersect}$}
\end{enumerate}

\begin{proof}
(i) Since $T\notin I_{\sigma_{j}}$ it follows that
$\sigma_j+\eps_{\sigma_j}\le T$.
Therefore there exists $\sigma\in\mathfrak{S}$ such that
$\sigma_j+\eps_{\sigma_j}\in I_\sigma$.
Since $I_\sigma$ is open it follows that 
$I_\sigma\cap I_{\sigma_j}\not=\emptyset$ and
$\sigma+\eps_\sigma> \sigma_j+\eps_{\sigma_j}$.
Therefore, by condition 2, $\sigma_{i+1}+\eps_{\sigma_{i+1}}
\ge \sigma+\eps_\sigma$ which is strictly larger than
$\sigma_j+\eps_{\sigma_j}$.

(ii) We assume by contradiction that there exists an interval
$I_{\sigma_{i}}$ intersecting $I_{\sigma_{j}}$ where $0\le i<i+2\le j\le N$.
Applying condition 2 for $j=i$ and using that 
$I_{\sigma_i}\cap I_{\sigma_j}\not=\emptyset$,
we obtain that $\sigma_{i+1}+\eps_{\sigma_{i+1}}
\ge \sigma_j+\eps_{\sigma_j}$.
Now applying (i) we obtain that
$\sigma_{i+2}+\eps_{\sigma_{i+2}}>\sigma_{i+1}+\eps_{\sigma_{i+1}}$
which as we saw is $\ge \sigma_j+\eps_{\sigma_j}$.
Using that $j\ge i+2$ and using (i) again,
we conclude that $\sigma_j+\eps_{\sigma_j}
\ge \sigma_{i+2}+\eps_{\sigma_{i+2}}$
which as we saw is $> \sigma_j+\eps_{\sigma_j}$.
This contradiction proves (ii).
\end{proof}

The family of intervals $\Ii:=\{I_{\sigma_0}, I_{\sigma_1},\dots,
I_{\sigma_N}\}$ covers $[-T,T]$ and has the property that exactly nearest
neighbors overlap, as illustrated by Figure~\ref{fig:fig-I_sigma}.
Set $t_0:=-T$ and $t_{N+1}:=T$.
For $i=1,\dots,N$ choose $t_i \in I_{\sigma_{i-1}} \cap I_{\sigma_i}$
in the overlap interval.
The finite set of non-eigenvalues is then defined by
$
   \Lambda^\prime
   :=\{\lambda_i:=\mu_{\sigma_i}\mid  i=0,\dots,N\}
$.
Note that $\lambda_0:=\mu_{-T}=0$.
This proves Step~2.
\end{proof}

\smallskip
\noindent
\textbf{Step~3.}
We prove the theorem.

\begin{proof}
We continue the notation from Step~2.
By Step~1, for $j=0,\dots,N$, the index along each interval
$[t_j,t_{j+1}]$ vanishes
$$
   \INDEX{\mathfrak{D}_{\A|_{[t_j,t_{j+1}]}}^{\lambda_j,\lambda_j}}=0 .
$$
By Lemma~\ref{le:vhj44} and
anti-symmetry~(\ref{eq:spectral-density_-}) of $\rho$ we have
\begin{equation*}
\begin{split}
   \INDEX{\mathfrak{D}_{\A|_{[t_j,t_{j+1}]}}^{\lambda_j,\lambda_{j+1}}}
   &=-\rho_{A(t_{j+1})}(\lambda_j,\lambda_{j+1}) +\rho_{A(t_j)} (\lambda_j,\lambda_j)
\\
   &=\rho_{A(t_{j+1})}(\lambda_{j+1}, \lambda_j) .
\end{split}
\end{equation*}
By path concatenation, Theorem~\ref{thm:concatenation}, and the previous
identity we get
\begin{equation}\label{eq:index-I_T-1}
\begin{split}
   \INDEX{\mathfrak{D}_A^{\lambda_0,\lambda_{N+1}}}
   =\sum_{j=0}^N
   \INDEX{\mathfrak{D}_{\A|_{[t_j,t_{j+1}]}}^{\lambda_j,\lambda_{j+1}}}
   =\sum_{j=0}^N
   \rho_{A(t_{j+1})}(\lambda_{j+1}, \lambda_j) .
\end{split}
\end{equation}
Given, at a time $s\in[-T,T]$, a non-eigenvalue $\mu\in\Rr(A(s))$, we define
\begin{equation*}
\begin{split}
   \nu_\uparrow(s;\mu)
   &:=\max\{\ell\in\Lambda_0\mid a_\ell(s)<\mu\}
   .
\end{split}
\end{equation*}
Then
\begin{equation}\label{eq:hgjhg4433}
\begin{split}
   \nu_\uparrow(T;0)
   &=-\varsigma(A)\quad\text{\small , cf.~(\ref{eq:SF-finite}),}
\\
   \rho_{A(t_{j+1})}(\lambda_{j+1}, \lambda_j)
   &=\nu_\uparrow(t_{j+1}; \lambda_{j})
   -\nu_\uparrow(t_{j+1}; \lambda_{j+1})
   .
\end{split}
\end{equation}
Since
$
   \A_j(s)
   :=A(s)-\lambda_j\iota
$
is invertible for every $s\in[t_j,t_{j+1}]$,
no eigenvalue of $A(s)$ crosses $\lambda_j$ along $[t_j,t_{j+1}]$,
and therefore
\begin{equation}\label{eq:dfrr3r}
   \nu_\uparrow(t_{j}; \lambda_{j})
   =\nu_\uparrow(t_{j+1}; \lambda_{j})
   .
\end{equation}
By~(\ref{eq:index-I_T-1}) we obtain identity 1 in
\begin{equation}\label{eq:ffgcfd22}
\begin{split}
   \INDEX{\mathfrak{D}_A^{\lambda_0,\lambda_{N+1}}}
   &\stackrel{1}{=}\sum_{j=0}^N
   \rho_{A(t_{j+1})}(\lambda_{j+1}, \lambda_j)
\\
   &\stackrel{2}{=}\sum_{j=0}^N
   \left(\nu_\uparrow(t_{j+1}; \lambda_{j})
   -\nu_\uparrow(t_{j+1}; \lambda_{j+1})\right)
\\
   &\stackrel{3}{=}\sum_{j=0}^N
   \left(\nu_\uparrow(t_{j}; \lambda_{j})
   -\nu_\uparrow(t_{j+1}; \lambda_{j+1})\right)
\\
   &\stackrel{4}{=}\nu_\uparrow(t_{0}; \lambda_{0})-\nu_\uparrow(t_{N+1}; \lambda_{N+1})
\\
   &\stackrel{5}{=}\nu_\uparrow(-T;0)-\nu_\uparrow(T;\lambda_{N+1})
\\
   &\stackrel{6}{=}-\nu_\uparrow(T;\lambda_{N+1})
\end{split}
\end{equation}
Identity 2 is by~(\ref{eq:hgjhg4433}) and
identity 3 by~(\ref{eq:dfrr3r}). In identity 4 all terms cancel
pairwise except the first and the last one. Identity 5 holds by Step~2
and identity 6 since $\nu_\uparrow(s;0)=0$.

\smallskip
\noindent
\textsc{Case~1.}
$\lambda_{N+1}=0$

\smallskip
\noindent
In this case $\mathfrak{D}_A^{0,0}=\mathfrak{D}_A$
and $-\nu_\uparrow(T;\lambda_{N+1})=
-\nu_\uparrow(T;0)=\varsigma(A)$.
Hence~(\ref{eq:ffgcfd22}) tells that
$\INDEX{\mathfrak{D}_A}=\varsigma(A)$ and we are done.

\smallskip
\noindent
\textsc{Case~2.}
$\lambda_{N+1}\not=0$

\smallskip
\noindent
By Lemma~\ref{le:vhj44} we obtain identity 2 in
\begin{equation}\label{eq:ffgcfd22-22}
\begin{split}
   \INDEX{\mathfrak{D}_A^{\lambda_0,\lambda_{N+1}}}-\INDEX{\mathfrak{D}_A}
   &= \INDEX{\mathfrak{D}_A^{0,\lambda_{N+1}}}-\INDEX{\mathfrak{D}_A^{0,0}}
\\
   &\stackrel{2}{=} \rho_{A(-T)}(0,0)- \rho_{A(T)}(0,\lambda_{N+1})
\\
   &\stackrel{3}{=} 0
   -\nu_\uparrow(T;\lambda_{N+1})
   +\nu_\uparrow(T;0)
\\
    &\stackrel{4}{=} -\nu_\uparrow(T;\lambda_{N+1}) - \varsigma(A) .
\end{split}
\end{equation}
Identities 3 and 4 hold by~(\ref{eq:hgjhg4433}).
Now~(\ref{eq:ffgcfd22}) and~(\ref{eq:ffgcfd22-22}) imply that 
$\INDEX{\mathfrak{D}_A}=\varsigma(A)$.
Together with Corollary~\ref{cor:Fredholm-I_T}
this proves Step~3, hence Theorem~\ref{thm:A}.
\end{proof}

The proof of Theorem~\ref{thm:A} is complete.
\end{proof}

\newpage
\subsection{Half infinite forward interval}

Pick a Hessian path $A\in \Aa_{I_+}^*$
along the half infinite forward interval $I_+=[0,\infty)$;
see Definition~\ref{def:Aa_I^*}.
Then $A\colon [0,\infty)\to \Ff=\Ff(H_1,H_0)$
takes values in the symmetrizable Fredholm operators of index zero;
cf. Remarks~\ref{rem:symmetrizable} and~\ref{rem:compact-op}.
In order to eventually get to Fredholm operators,
it is not enough that the Hessian at zero and the limit at infinity
are invertible, notation
$$
   \A_0:=A(0)
   ,\qquad
   \A^+:=\lim_{s\to\infty} A(s) .
$$
In addition, one must impose a boundary condition at zero
formulated in terms of the spectral projection
$\pi_+^{\A_0}$ sitting at time zero and; see~(\ref{eq:H_1/2}).

\boldmath
\subsubsection{Estimate for $D_A$}
\unboldmath

Let $A\in \Aa_{I_+}^*$.
The Hilbert spaces $P_0(\R_+)$ and $P_1(\R_+)$ are defined by~(\ref{eq:P_1-I})
for $I=\R_+$.
In~this~section we study the linear operator  $\p_s+A$ as a map
\begin{equation}\label{eq:D_A-half-fin}
   D_A\colon P_1(\R_+)\to P_0(\R_+)
   ,\quad
   \xi\mapsto \p_s\xi+A(s)\xi .
\end{equation}
As in the case of the finite interval, Section~\ref{sec:I_T-D_A},
this operator is \emph{not} Fredholm: although it has closed image and finite
dimensional co-kernel, the kernel is infinite dimensional in the Floer
and Morse case; see Figure~\ref{fig:D_A}.

\begin{theorem}\label{thm:half-finite-interval}
Given $A\in \Aa_{I_+}^*$, there exist constants $T,c>0$ such that
\[
   \Norm{\xi}_{P_1(\R_+)}
   \le c\left(\norm{\xi}_{P_0([0,T])}+\norm{D_A\xi}_{P_0(\R_+)}
   +\norm{\pi_+^{\A_0}\xi(0)}_{\frac12}\right)
\]
for every $\xi\in P_1(\R_+)$.
\end{theorem}

This estimate becomes a semi-Fredholm estimate for $D_A$ restricted
to those $\xi\in P_1(\R_+)$ with $\pi_+^{\A_0}\xi(0)=0$
or even $\xi(0)=0$. We study this in Section~\ref{sec:forward-Fredholm}.

\begin{proof}[Proof of Theorem~\ref{thm:half-finite-interval}]
We prove the theorem in four steps.
It is sometimes convenient to abbreviate $A_s:=A(s)$.
We enumerate the constants by the step where they appear,
e.g. constant $C_1$ arises in Step~1.

\smallskip
\noindent
\textbf{Step~1} (Asymptotic estimate)\textbf{.}
There exist constants $T_1,C_1>0$ such that the following is true.
Suppose $\beta\in C^\infty(\R_+,\R)$ satisfies
$
   \supp\beta\subset (T_1,\infty)
$.
Then
\[
   \Norm{\beta\xi}_{P_1(\R_+)}
   \le C_1\left(
   \Norm{\beta D_A\xi}_{P_0(\R_+)}+\Norm{\beta^\prime\xi}_{P_0(\R_+)}
   \right)
\]
for every $\xi\in P_1(\R_+)$.

\begin{proof}
Step~3 in the proof of the Rabier Theorem~\ref{thm:Rabier}.
\end{proof}

\smallskip
\noindent
\textbf{Step~2} (Small interval at left boundary)\textbf{.}
There are constants $\eps_2>0$ and $C_2>0$ such that for every
compactly supported $\beta\in C^\infty(\R_+,\R)$ with the property
$$
   \sup_{\sigma,\tau\in\supp\beta}\norm{A_\sigma-A_\tau}_{\Ll(H_1,H_0)}\le\eps_2
$$
it holds that
\begin{equation*}
\begin{split}
   &\Norm{\beta\xi}_{P_1(\R_+)}\\
   &\le C_2\left(
   \Norm{\beta D_A\xi}_{P_0(\R_+)}+\Norm{\beta^\prime\xi}_{P_0(\R_+)}
   +\Norm{\beta\xi}_{P_0(\R_+)}
   +\norm{\pi_+^{\A_0}\beta(0)\xi(0)}_{\frac12}
    \right)
\end{split}
\end{equation*}
for every $\xi\in P_1(\R_+)$.

\begin{proof}
Step~4 in the proof of the finite interval
Theorem~\ref{thm:finite-interval}.
\end{proof}

\smallskip
\noindent
\textbf{Step~3} (Small interior interval)\textbf{.}
There is a finite subset $\Lambda^\prime\subset\R$
and constants $\eps_3,C_3>0$
such that for every $\beta\in C^\infty(\R_+,\R)$ which has compact
support in $(0,\infty)$ and has the property
$$
   \sup_{\sigma,\tau\in\supp\beta}\norm{A_\sigma-A_\tau}_{\Ll(H_1,H_0)}\le\eps_3
$$
it holds that
\[
   \Norm{\beta\xi}_{P_1(\R_+)}
   \le C_3\left(
   \Norm{\beta D_A\xi}_{P_0(\R_+)}+\Norm{\beta^\prime\xi}_{P_0(\R_+)}
   +\Norm{\beta\xi}_{P_0(\R_+)}
    \right)
\]
for every $\xi\in P_1(\R_+)$.

\begin{proof}
This is Step~6 in the proof of the Rabier Theorem~\ref {thm:Rabier}
with $\frac{1}{C_3}=\eps_3$.
\end{proof}

\smallskip
\noindent
\textbf{Step~4} (Partition of unity)\textbf{.}
We prove Theorem~\ref{thm:half-finite-interval}.

\begin{proof}
Set $\eps:=\min\{\eps_2,\eps_3\}$ and $C:=\max\{C_1,C_2,C_3\}$.
Choose $T>T_1$ and a finite partition of unity
$\{\beta_j\}_{j=0}^{M+1}$ for $[0,\infty)$ with the properties that
$\beta_0$ is compactly supported in $[0,T)$ and
$$
   \beta_0(0)=1,\quad
   \sup_{\sigma,\tau\in\supp\beta_0}\norm{A_\sigma-A_\tau}_{\Ll(H_1,H_0)}
   \le\eps,\quad
   \supp\beta_{M+1}\subset (T_1,\infty),
$$
and
 $$
   \sup_{\sigma,\tau\in\supp\beta_j}\norm{A_\sigma-A_\tau}_{\Ll(H_1,H_0)}
   \le\eps
   ,\qquad
  \supp\beta_j\subset (0,T) ,
$$
for $j=1,\dots,M$.
That such a partition exists follows from the continuity
of $s\mapsto A(s)$ and the fact that on the compact set $[0,T_1]$
continuity becomes uniform continuity.
Let $\xi\in P_1(\R_+)$. Then by Steps~{\color{gray}2,1,3} we have the estimates
\begin{equation*}
\begin{split}
   \Norm{\beta_0\xi}_{P_1(\R_+)}
   &\stackrel{{\color{gray}2}}{\le} C\left(
   \Norm{\beta_0 D_A\xi}_{P_0(\R_+)}+\Norm{\beta_0^\prime\xi}_{P_0(\R_+)}
   +\Norm{\beta_0\xi}_{P_0(\R_+)}
   +\norm{\pi_+^{\A_0}\xi(0)}_{\frac12}
    \right)
\\
   \Norm{\beta_{M+1}\xi}_{P_1(\R_+)}
   &\stackrel{{\color{gray}1}}{\le} C\left(
   \Norm{\beta_{M+1} D_A\xi}_{P_0(\R_+)}+\Norm{\beta_{M+1}^\prime\xi}_{P_0(\R_+)}
   \right)
\\
   \Norm{\beta_j\xi}_{P_1(\R_+)}
   &\stackrel{{\color{gray}3}}{\le} C\left(
   \Norm{\beta_j D_A\xi}_{P_0(\R_+)}+\Norm{\beta_j^\prime\xi}_{P_0(\R_+)}
   +\Norm{\beta_j\xi}_{P_0(\R_+)}
   \right)
\end{split}
\end{equation*}
for $j=1,\dots,M$.
We abbreviate $B:=\max\{\norm{\beta_0^\prime}_\infty,
\norm{\beta_1^\prime}_\infty,\dots, \norm{\beta_{M+1}^\prime}_\infty\}$.
Putting these estimates together we obtain
\begin{equation*}
\begin{split}
   \Norm{\xi}_{P_1(\R_+)}
   &\le\sum_{j=0}^{M+1} \Norm{\beta_j\xi}_{P_1(\R_+)}\\
   &\le C \sum_{j=0}^{M+1}\left(
   \Norm{\beta_j D_A\xi}_{P_0(\R_+)}+\Norm{\beta_j^\prime\xi}_{P_0([0,T])}
  \right)
\\
   &\quad
   +C \sum_{j=0}^{M}\Norm{\beta_j\xi}_{P_0([0,T])}
   +C \norm{\pi_+^{\A_0}\xi(0)}_{\frac12}
\\
   &\le C (M+2) \Norm{D_A\xi}_{P_0(\R_+)}
   +C \left(B (M+2)+ M+1\right) \Norm{\xi}_{P_0([0,T])}
\\
   &\quad
   +C \norm{\pi_+^{\A_0}\xi(0)}_{\frac12}
\end{split}
\end{equation*}
where in the second inequality we replaced the $P_0(\R_+)$ norm by the
$P_0([0,T])$ norm due to the supports of the $\beta_j$'s and their
derivatives.\footnote{
   Along $[T,\infty)$ we have
   $\beta_{M+1}\equiv 1$, so $\beta_{M+1}^\prime\equiv 0$.
   }
Setting
$$
   c:=\max\{C (M+2), C \left(B (M+2)+ M+1\right)\}
$$
proves Step~4.
\end{proof}
\noindent
The proof of Theorem~\ref{thm:half-finite-interval} is complete.
\end{proof}

\boldmath
\subsubsection{Estimate for the adjoint $D_A^*$}
\unboldmath

Let $A\in\Aa_{\R_+}^*$.
We call the following operator the
\textbf{adjoint of \boldmath$D_A$}, namely
$$
   D_A^*:=D_{-A^*}\colon P_1(\R_+;H_0^*,H_1^*)\to P_0(\R_+;H_1^*)
   ,\quad
   \eta\mapsto \p_s\eta-A(s)^*\eta .
$$

\begin{corollary}\label{cor:forward-interval_-A*}
For $A\in\Aa_{\R_+}^*$ there exists a constant $c>0$ such that
\begin{equation*}
\begin{split}
   \norm{\eta}_{P_1(\R_+;H_0^*,H_1^*)}
   &\le c\Bigl(\norm{\eta}_{P_0(\R_+;H_1^*)}+\norm{D_A^*\eta}_{P_0(\R_+;H_1^*)}
   +\norm{\pi_+^{-\A^*_{0}}\eta(0)}_{\frac12}\Bigr)
\end{split}
\end{equation*}
for every $\eta\in P_1(\R_+;H_0^*,H_1^*)$.
\end{corollary}

\begin{proof}
Theorem~\ref{thm:half-finite-interval} and Lemma~\ref{le:adjoint-map};
see also Remark~\ref{rem:symmetrizable}.
\end{proof}

\boldmath
\subsubsection{Fredholm under boundary conditions: $D_A^+$}
\label{sec:forward-Fredholm}
\unboldmath

Given $A\in\Aa_{\R_+}^*$,
let $\pi_\pm:=\pi_\pm^{\A_0}$ be defined by~(\ref{eq:H_1/2}).
To get from Theorem~\ref{thm:half-finite-interval}
to semi-Fredholm we endow the domain
of the operator $D_A\colon P_1(\R_+)\to P_0(\R_+)$ with the
boundary condition $\pi_+^{\A_0}\xi(0)=0$ which cuts the operator
kernel down to finite dimension and has a finite dimensional co-kernel.
Hence $\coker D_A$ is finite dimensional, too.

To this end define a subspace of the Hilbert space
$P_1 (\R_+)=P_1(\R_+;H_1,H_0)$ from~(\ref{eq:P_1-I}) as follows
\begin{equation}\label{eq:P_1^+}
\begin{split}
   P_1^+(\R_+,\A_0){\color{gray}\;=P_1^+(\R_+,\A_0 ;H_1,H_0)}
   &:=\{\xi\in P_1(\R_+)\mid \pi_+^{\A_0}\xi(0)=0\}.
\end{split}
\end{equation}
The restriction of the operator $D_A\colon P_1 (\R_+)\to P_0 (\R_+)$
in~(\ref{eq:D_A-half-fin}) we denote by
\begin{equation*}
\begin{split}
   D_A^+\colon P_1^+(\R_+,\A_0) \to P_0 (\R_+)
   ,\quad
   \xi\mapsto \p_s\xi+A(s)\xi .
\end{split}
\end{equation*}

\begin{remark}[Goal and idea of proof]
Our goal is to show that
$D_A$ has finite dimensional cokernel.

To achieve this goal we show that $D_A^+$ is a Fredholm operator
whose co-kernel is isomorphic to $\ker D_{-A^*}^+$.
The fact that $D_A^+$ has closed image is crucial
to show that $D_A$ itself has closed image (since it contains $\im D_A^+$).

For the proof that $D_A^+$ is a semi-Fredholm operator
we need the full strength of the estimate in
Theorem~\ref{thm:half-finite-interval}, in particular, that
the third term on the right is just
$\norm{\pi_+^{\A_0}\xi(0)}_{H_{1/2}}$ and not $\norm{\xi(0)}_{H_{1/2}}$.
\end{remark}

\begin{figure}[h]
\begin{center}
\begin{tabular}{lll}
      
  & $D_A\colon P_1\to P_0$         
  & $D_A^+\colon P_1^+\to P_0$ 
\\
\midrule
  $\dim \ker$
  & \small ${\color{red}\infty}$
  & \small $k<\infty$ 
\\
  $\dim \coker$
  & \small $\le\ell$ \qquad\qquad ${\color{gray}\Leftarrow}$
  & \small $\ell<\infty$
\\
\midrule
     \small 
  & \small co-semi-Fredholm
  & \small Fredholm 
\\
\midrule
  $\mathrm{image}$
  & \small closed \qquad\quad ${\color{gray}\Leftarrow}$
  & \small closed
\\
  $\coker$
  & \small 
  & \small $\coker D_A^+\simeq \ker D_{-A^*}^+$
\\
  $\ker$
  & \small {\color{red} huge}
  & \small ${\color{gray} \ker D_A^+\simeq\coker D_{-A^*}^+}$
\end{tabular}
\caption{$D_A=\p_s+A(s)$ on $P_1$ and its restriction $D_A^+$ to $P_1^+$}
\label{fig:D_A}
\end{center}
\end{figure}

\begin{theorem}[Fredholm]\label{thm:Fredholm-forward}
$D_A^+\colon P_1^+(\R_+,\A_0 ;H_1,H_0)\to P_0(\R_+;H_0)$ 
is a Fredholm operator for any Hessian path $A\in\Aa_{\R_+}^*$.
\end{theorem}

\begin{corollary}\label{cor:D_A-closed-image}
The operator $D_A\colon P_1(\R_+;H_1,H_0)\to P_0(\R_+;H_0)$
in~(\ref{eq:D_A-half-fin}) has closed image of finite co-dimension
for any Hessian path $A\in\Aa_{\R_+}^*$.
\end{corollary}

\begin{proof}
By Theorem~\ref{thm:Fredholm-forward} the image of $D_A^+$
is closed and of finite co-dimension.
Since $D_A^+$ is a restriction of $D_A$ we have inclusion
$\im D_A^+\subset \im D_A {\color{gray}\,\subset P_0(\R_+)}$.
So $\im D_A$ is of finite co-dimension.
Thus $\im D_A$ is closed by~\cite[Prop.\,11.5]{brezis:2011a}.
\end{proof}

\begin{proof}[Proof of Theorem~\ref{thm:Fredholm-forward}]
Pick $A\in\Aa_{\R_+}^*$, then $\A_0:=A(0)$ is invertible.
By Corollary~\ref{cor:forward-interval_-A*}
the operator  $D_A^+$ (and also $D_{-A^*}^+$)
has finite dimensional kernel and closed image.
By the same reasoning as in the proof of
Theorem~\ref{thm:Fredholm-I_T-D_A+-}
one shows that the co-kernel of $D_A^+$
can be identified with the kernel of $D_{-A^*}^+$,
in symbols
\begin{equation}\label{eq:coker-sim-ker}
   \coker{D_A^+}\simeq\ker{D_{-A^*}^+} .
\end{equation}
This proves Theorem~\ref{thm:Fredholm-forward}.
\end{proof}

\boldmath
\subsubsection{Paths of invertibles}
\unboldmath

\begin{proposition}[Constant path]\label{prop:const-path}
Let $\A\in\Aa_{\R_+}^*$ be a constant path, then the Fredholm operator
$D_{\A}^+\colon P_1^+(\R_+,\A_0)\to P_0(\R_+)$ is an isomorphism
and therefore its Fredholm index vanishes.
\end{proposition}

\begin{proof} The proof is in three steps.
After replacing the inner products by $\A$-adaptable
inner products, see Definition~\ref{def:A-adapted},
we can assume without loss of generality
that $\A\colon H_1\to H_0$ is a symmetric isometry.

\smallskip\noindent
\textbf{Step~1:} $\ker D_{\A}^+=\{0\}$\textbf{.}
We first show that the kernel of $D_{\A}^+$ vanishes. 
For this purpose suppose that $\xi\in P_1^+(\R_+,\A_0)$ lies in the kernel
of $D_{\A}^+$.
Then $\xi$ is a solution of the problem
$\p_s\xi(s)=-\A\xi(s)$ and $\pi_+\xi(0)=0$.
\\
Pick an orthonormal basis $\Vv(\A)=\{v_\ell\}_{\ell\in\Lambda}$ of $H_0$
consisting of eigenvectors $\A v_\ell=a_\ell v_\ell$.
We write $\xi=\sum_{\ell\in\Z^*}\xi_\ell v_\ell$.
Then each coefficient $\xi_\ell$ satisfies the ODE in one variable
$\p_s\xi_\ell(s)=-a_\ell\xi_\ell(s)$ whose solution is
$\xi_\ell(s)=e^{-a_\ell s}\xi_\ell(0)$.
Since $\pi_+\xi(0)=0$ we have $\xi_\nu(0)=0$ for every $\nu\in\N$.
Therefore $\xi_\nu=0$ for every $\nu\in\N$.
Since $a_{-\nu}<0$ is negative for $\nu\in\N$
we have that $\xi_{-\nu}(s)=e^{-a_{-\nu} s}\xi_{-\nu}(0)$ grows
exponentially unless $\xi_{-\nu}(0)=0$.
Since $\xi\in P_1^+(\R_+,\A_0)\subset L^2(\R_+,H_1)$ negative modes
$\xi_{-\nu}$ cannot grow exponentially which implies that
$\xi_{-\nu}\equiv 0$  for every $\nu\in\N$.
This shows that $\xi=0$. So $\ker D_{\A}^+=\{0\}$ is trivial.

\smallskip\noindent
\textbf{Step~2:} $\coker D_{\A}^+=\{0\}$\textbf{.}
But $\coker D_{\A}^+\simeq\ker D_{-\A^*}^+$, by~(\ref{eq:coker-sim-ker}),
and the latter is zero by Step~1.

\smallskip\noindent
Step~1 and Step~2 show that $D_{\A}^+$ is bijective and hence, by the open
mapping theorem, an isomorphism.
This proves Proposition~\ref{prop:const-path}.
\end{proof}

\begin{corollary}\label{cor:D_A^+_index=0-new}
Assume that $\A\in\Aa_{\R_+}^*$ has the property that $\A(s)$ is invertible
for every $s\in\R_+$. Then the Fredholm index
$\INDEX(\mathfrak{D}_\A)=0$ vanishes.
\end{corollary}

\begin{proof}
For constant paths this is true by Proposition~\ref{prop:const-path}.
The family of paths $\{\A_r\}_{r\in[0,1]}\subset \Aa_{\R_+}^*$
defined by
\begin{equation}\label{eq:homotopy-A_lambda}
   \A_r(s):=\A(s+\varphi(r))
   ,\qquad
   \varphi\colon[0,1]\to \R_+\cup\{\infty\}
   ,\quad
   r\mapsto\tfrac{r^2}{1-r^2},
\end{equation}
provides a homotopy between $\A$ and the constant path $\A^+$ at infinity.
Therefore, since by Theorem~\ref{thm:Fredholm-index-constant} the
Fredholm index is invariant under homotopies
$$
   r\mapsto \mathfrak{D}_{\A_r} =\Dd_r\circ p_r
   \colon P_1\to P_0\times H_\frac12^+(\A_r(0))
$$
through Fredholm operators (true by Corollary~\ref{cor:Fredholm-I_+}),
the index of $\mathfrak{D}_{\A_r}$ is constant. 
In the case at hand
the operators are the following
$$
   \Dd_r\colon P_1\to P_0\times H_\frac12
   ,\quad
   \xi\mapsto\left(D_{A_r}\xi,\xi(0)\right)
$$
and
$$
   p_r=\left(\Id,\pi_+^{\A_r(0)}\right)
   \colon
   P_0\times H_\frac12\to P_0\times H_\frac12 .
$$
The map $r\mapsto p_r$ is continuous by~\cite[Thm.\,D]{Frauenfelder:2024c},
see Theorem~\ref{thm:proj}.
It remains to show continuity of the homotopy
$[0,1]\ni r\mapsto \Dd_r$, hence of $r\mapsto D_{A_r}$.
Next we show this at $r=1$.
By continuity of the path $\sigma\mapsto A(\sigma)$,
given $\eps>0$, there exists $\sigma_0=\sigma_0 (\eps)>0$ such that
$\norm{\A^+-A(\sigma)}_{\Ll(H_1,H_0)}\le\eps$ for every $\sigma\ge \sigma_0$.
Let $r_0$ be such that $r_0^2/(1-r_0^2)=\sigma_0$.
Since the function $\varphi$ is monotone increasing,
for every $r\in[r_0,1]$ we have
$r^2/(1-r^2)\ge \sigma_0$.
Therefore for every $s\in\R_+$ we have
$\norm{\A^+-A_r(s)}_{\Ll(H_1,H_0)}^2\le \eps$.
Hence there is the estimate
\begin{equation*}
\begin{split}
   \norm{(D_{\A^+}-D_{A_r})\xi}_{P_0(\R_+)}^2
   &=\int_0^1\norm{(\A^+-A_r(s))\xi(s)}_0^2 ds\\
   &\le \int_0^1\norm{\A^+-A_r(s)}_{\Ll(H_1,H_0)}^2 \norm{\xi(s)}_1^2 ds\\
   &\le\eps^2\norm{\xi}_{P_1(\R_+)}^2 .
\end{split}
\end{equation*}
This proves that $\norm{D_{\A^+}-D_{A_r}}_{\Ll(P_1^+,P_0)}\le\eps$.
This shows continuity at $r=1$.
For $r\in[0,1]$ one compares $D_{A_r}$ and
$D_{A_{\tilde r}}$ by a similar argument where, in addition,
uniform continuity of $\sigma\mapsto A(\sigma)$ enters.
Now it follows from Theorem~\ref{thm:Fredholm-index-constant}
that $\INDEX{\mathfrak{D}_{\A_0}}=\INDEX{\mathfrak{D}_{\A_1}}$.
Since $\A_0=\A$ and $\A_1\equiv \A^+$, we get identity 1 in
$$
   \INDEX{\mathfrak{D}_\A}
   \stackrel{1}{=}\INDEX{\mathfrak{D}_{\A^+}}
   \stackrel{2}{=}\INDEX{D_{\A^+}^+}
   \stackrel{3}{=} 0
$$
where identity 2 holds by the same arguments
as in the proof of Corollary~\ref{cor:Fredholm-I_T}
and identity 3 is by Proposition~\ref{prop:const-path}.
This proves Corollary~\ref{cor:D_A^+_index=0-new}.
\end{proof}

\boldmath
\subsubsection{Theorem~\ref{thm:A} -- Fredholm property}
\unboldmath

To prove Theorem~\ref{thm:A}
in the half infinite forward interval case, namely that $\mathfrak{D}_A$
is Fredholm and $\INDEX{\mathfrak{D}_A}=\varsigma(A)$,
we use Theorem~\ref{thm:proj}
(\hspace{-.1pt}\cite[Thm.\,D]{Frauenfelder:2024c}).

\begin{corollary}[to Theorem~\ref{thm:half-finite-interval}, Fredholm]
\label{cor:Fredholm-I_+}
For any $A\in \Aa_{\R_+}^*$ the operator
\begin{equation*}
\begin{split}
   \mathfrak{D}_A=\mathfrak{D}_A^{\R_+}\colon P_1(\R_+)
   &\to P_0(\R_+)\times H^+_{\frac12}(\A_0)
   =: \Ww(\R_+;\A_0)
   \\
   \xi
   &\mapsto
   \left(D_A\xi,\pi_+^{\A_0}\xi_0\right)
\end{split}
\end{equation*}
is Fredholm, where $\A_0:= A(0$,
and of the same index as $D_A^+$.
More precisely, the kernels coincide and the co-kernels are of equal dimension.
\end{corollary}

\begin{proof}
By Theorem~\ref{thm:half-finite-interval}
the operator $\mathfrak{D}_A$ is semi-Fredholm.
So it has finite dimensional kernel and closed image.
That kernel and image of $\mathfrak{D}_A$ are equal,
respectively isomorphic, to those of the Fredholm operator $D_A^+$
from Theorem~\ref{thm:Fredholm-forward}
follows by the arguments in the proof of
Corollary~\ref{cor:Fredholm-I_T}.
\end{proof}

\subsubsection{Theorem~\ref{thm:A} -- Index formula}

Pick $A\in \Aa_{\R_+}^*$.
Choose $T>0$ sufficiently large such that
$A(s)$ is invertible for every $s\ge T$.
Analogous to Theorem~\ref{thm:concatenation} there is
concatenation formula 1
\begin{equation*}
\begin{split}
   \INDEX{\mathfrak{D}_A}
   &\stackrel{1}{=}
   \INDEX{\mathfrak{D}_A|_{[0,T]}}+\INDEX{\mathfrak{D}_A|_{[T,\infty)}}\\
   &\stackrel{2}{=}\INDEX{\mathfrak{D}_A|_{[0,T]}}\\
   &\stackrel{3}{=}\varsigma(\mathfrak{D}_A|_{[0,T]})\\
   &\stackrel{4}{=}\varsigma(\mathfrak{D}_A) 
\end{split}
\end{equation*}
Identity 2 is Corollary~\ref{cor:D_A^+_index=0-new}
and identity 3 is the spectral flow formula of Theorem~\ref{thm:A}
in the already proved finite interval case.
Identity 4 holds since $A(s)$ is invertible for every
$s\in[T,\infty]$; no eigenvalues cross zero.

\subsection{Half infinite backward interval -- Theorem~\ref{thm:A}}

Let $\R_-=(-\infty,0]$.
For $k=0,1$ we define the Hilbert space isomorphism
$\Rr_k\colon P_k(\R_-)\to P_k(\R_+)$
by $(\Rr_k\xi)(s)=\xi(-s)$ for $s\in\R_-$.
We define the map $\Aa_{\R_-}^*\to \Aa_{I_+}^*$ by
$(\Rr A)(s):=-A(-s)$ for $s\in\R_-$.
Let $A\in \Aa_{\R_-}^*$. Consider
\begin{equation*}
\begin{split}
   \mathfrak{D}_A=\mathfrak{D}_A^{\R_-}\colon P_1(\R_-)
   &\to P_0(\R_-)\times H^-_{\frac12}(\A_0)
   =: \Ww(\R_-;\A_0)
   \\
   \xi
   &\mapsto
   \left(D_A\xi,\pi_-^{\A_0}\xi_0\right)
\end{split}
\end{equation*}
where $\A_0:= A(0$.
Note that
$$
   (\Rr_0,\1)\circ \mathfrak{D}_A\circ\Rr_1
   =\mathfrak{D}_{\Rr A}
   \colon P_1(\R_+)\to P_0(\R_+) \times
   \underbrace{H^-_\frac12(\A_0)}_{H^+_\frac12((\Rr A)_0)} .
$$
Hence $\mathfrak{D}_A$ is a Fredholm operator and it holds that
\begin{equation*}
\begin{split}
   \INDEX{\mathfrak{D}_A}
   &\stackrel{1}{=}\INDEX{\mathfrak{D}_{\Rr A}}\\
   &\stackrel{2}{=}\varsigma(\Rr A)\\
   &\stackrel{3}{=}\varsigma(A) .
\end{split}
\end{equation*}
Here identity 1 is by the previous displayed conjugation,
identity 2 is by the already proven Theorem~\ref{thm:A} for $\R_+$,
and identity 3 holds since the path $\Rr A$ is
the negative of the path $A$ traversed backwards,
the two minus signs cancel.

\subsection{Real line -- Theorem~\ref{thm:A}}

Let $A\in \Aa_{\R}^*$.
Corollary~\ref{cor:Fredholm-D_A}
shows that $\mathfrak{D}_A=D_A\colon P_1(\R)\to P_0(\R)$ is Fredholm.

To show the spectral flow formula pick $T>0$ such that $A(s)$
invertible whenever $s\ge \abs{T}$.
There is the concatenation identity 1
\begin{equation*}
\begin{split}
   \INDEX{\mathfrak{D}_A}
   &\stackrel{1}{=}
   \INDEX{\mathfrak{D}_A|_{(-\infty,-T]}}
   +\INDEX{\mathfrak{D}_A|_{[-T,T]}}
   +\INDEX{\mathfrak{D}_A|_{[T,\infty)}}\\
   &\stackrel{2}{=}\INDEX{\mathfrak{D}_A|_{[-T,T]}}\\
   &\stackrel{3}{=}\varsigma(\mathfrak{D}_A|_{[-T,T]})\\
   &\stackrel{4}{=}\varsigma(\mathfrak{D}_A) 
\end{split}
\end{equation*}
Identity 2 is Corollary~\ref{cor:D_A^+_index=0-new}
and identity 3 is the spectral flow formula of Theorem~\ref{thm:A}
in the already proved finite interval case.
Identity 4 holds since $A(s)$ is invertible whenever
$\abs{s}\ge T$; no eigenvalues cross zero.

\newpage 
\appendix
\section{Hilbert space pairs}\label{sec:HS-pairs}

\boldmath
\subsection{Interpolation and extrapolation: Hilbert $\R$-scales}
\unboldmath

Let $H=(H_0,H_1)$ be a Hilbert space pair.
Then both Hilbert spaces $H_0$ and
$H_1$ are separable by~\cite[Cor.\,A.5]{Frauenfelder:2024c}.
By Riesz' theorem there is a unique bounded linear map $T\in\Ll(H_1)$,
called the \textbf{growth operator} of the pair, with
\begin{equation}\label{eq:Riesz}
   \INNER{\xi}{\eta}_0=\INNER{\xi}{T\eta}_1
\end{equation}
for all $\xi,\eta\in H_1$.
Since $\INNER{\cdot}{\cdot}_0$ and $\INNER{\cdot}{\cdot}_1$ are inner
products, the operator $T$ is positive definite and
symmetric. Moreover, in~\cite[Le.~A.7]{Frauenfelder:2024c}
we showed that compactness of the inclusion
$\iota\colon H_1\to H_0$ implies that the operator $T$ is compact.
In particular, the spectrum of $T$ consists of positive eigenvalues $\kappa$,
of finite multiplicity $m_\kappa$, whose only accumulation point is zero.
Define
$$
   \forall \kappa\in \spec T
   ,\quad
   V_\kappa:=\Eig_\kappa T:=\{v \in H_1\mid T v=\kappa v\}
   ,\quad
   m_\kappa:=\dim V_\kappa {\color{gray}\;<\infty} ,
$$
then the \textbf{eigenspace core} of the pair $(H_0,H_1)$ is the
direct sum of eigenspaces
$$
   V:=\bigoplus_{\kappa\in\spec T} V_\kappa ,
   \qquad V \subset H_1\subset H_0 .
$$
For later use, the direct sum is in decreasing eigenvalue order
$\kappa_1>\kappa_2>\dots>0$.
As a consequence of the spectral theorem for compact symmetric operators
$$
   H_1=\overline{V}^{\norm{\cdot}_1}.
$$
Since $H_1$ is a dense subset of $H_0$
we further have
$$
   H_0=\overline{V}^{\norm{\cdot}_0}.
$$

\begin{lemma}\label{le:0-1-orthogonality}
Let $\kappa_1\not = \kappa_2$ be different eigenvalues of $T$.
Then the eigenspaces $V_{\kappa_1}$ and $V_{\kappa_2}$ are
orthogonal with respect to both inner products
$\inner{\cdot}{\cdot}_0$ and~$\inner{\cdot}{\cdot}_1$.
Two vectors of $V$ are $0$-orthogonal
iff they are $1$-orthogonal, in symbols
$\perp_1\Leftrightarrow\perp_0$.
\end{lemma}

\begin{proof}
Pick $\xi_1\in V_{\kappa_1}$ and $\xi_2\in V_{\kappa_2}$.
This means that $T\xi_1=\kappa_1\xi_1$ and $T\xi_2=\kappa_2\xi_2$.
Using~(\ref{eq:Riesz}) we compute 
\begin{equation*}
\begin{split}
   \kappa_2 \INNER{\xi_1}{\xi_2}_1
   &=\INNER{\xi_1}{T\xi_2}_1\\
   &\stackrel{2}{=}\INNER{\xi_1}{\xi_2}_0\\
   &=\INNER{\xi_2}{\xi_1}_0\\
   &=\INNER{\xi_2}{T\xi_1}_1\\
   &=\kappa_1 \INNER{\xi_1}{\xi_2}_1.
\end{split}
\end{equation*}
The hypothesis $\kappa_1\not = \kappa_2$ implies $1$-orthogonality
$\INNER{\xi_1}{\xi_2}_1=0$.
So $\INNER{\xi_1}{\xi_2}_0=0$, by equality 2.
For $\xi_1,\xi_2\in V_{\kappa_2}$ equality 2
proves assertion two of the lemma.
\end{proof}

Another immediate consequence of~(\ref{eq:Riesz})
is the \textbf{length relation} in $V_\kappa$, namely
\begin{equation}\label{eq:0-1-lengths}
   \xi\in V_\kappa\qquad\Rightarrow\qquad
   \norm{\xi}_1=\frac{1}{\sqrt{\kappa}}\norm{\xi}_0 .
\end{equation}
We write $\xi\in V$ uniquely in the form $\xi=\sum_{\kappa\in\spec T}
\xi_\kappa$ where $\xi_\kappa\in V_\kappa$. Then
\begin{equation}\label{eq:0-1-squares}
   \norm{\xi}_0^2=\sum_{\kappa\in\spec T}\norm{\xi_\kappa}_0^2,
   ,\qquad
   \norm{\xi}_1^2=\sum_{\kappa\in\spec T}\frac{1}{\kappa}\norm{\xi_\kappa}_0^2.
\end{equation}
The first formula is by $0$-orthogonality in
Lemma~\ref{le:0-1-orthogonality}
and the second formula by $1$-orthogonality in
Lemma~\ref{le:0-1-orthogonality} combined with~(\ref{eq:0-1-lengths}).

\medskip
For any real $r\in\R$ we define an $r$-norm for $\xi\in V$ by
$$
   \norm{\xi}_{H_r}
   :=\biggl(\sum_{\kappa\in\spec T}
   \frac{1}{\kappa^r}\norm{\xi_\kappa}_0^2\biggr)^\frac12.
$$
Since $V$ is a direct product, for any element $\xi$ only finitely many
components $\xi_\kappa$ are non-zero, hence the number of non-zero
summands, also in~(\ref{eq:0-1-squares}), is finite.
By~(\ref{eq:0-1-squares}), the definition of the $r$-norm coincides
for $r=0,1$ with the original norms in $H_0$ and $H_1$, respectively.
Now we 
take the completion
\begin{equation}\label{eq:H_r}
   H_r:=\overline{V}^{\norm{\cdot}_{H_r}} .
\end{equation}
We endow $H_r$ with the \textbf{pair \boldmath$r$-inner product}
and the \textbf{pair \boldmath$r$-norm} defined~by
\begin{equation}\label{eq:inner_r}
    \INNER{\xi}{\eta}_{H_r}
   :=\sum_{\kappa\in\spec T} \frac{1}{\kappa^r}
   \INNER{\xi_\kappa}{\eta_\kappa}_0
   ,\quad
   \norm{\xi}_{H_r}
   :=\biggl(\sum_{\kappa\in\spec T}
   \frac{1}{\kappa^r}\norm{\xi_\kappa}_0^2\biggr)^\frac12 ,
\end{equation}
whenever $\xi,\eta\in H_r$. Here the number of non-zero summands could
be infinite, but the sum is still finite due to the completion property.

\medskip
To summarize, a Hilbert space pair $H=(H_0,H_1)$
canonically induces a \textbf{Hilbert \boldmath$\R$-scale},
roughly speaking a real family of Hilbert spaces $H_r$, notation
\begin{equation}\label{eq:(H_r)}
   H_\R:=(H_r)_{r\in\R} .
\end{equation}
The \textbf{dual Hilbert \boldmath$\R$-scale} is defined by
$
   H^*_\R=(H_r^*)_{r\in\R}
$
where
$
   H_r^*:=\Ll(H_r,\R)
$.

\boldmath
\subsubsection{The model Hilbert $\R$-scale}
\unboldmath

Let $\gf{f}\colon \N\to(0,\infty)$ be a \textbf{growth function}, i.e.
a monotone unbounded function.
Let $\ell^2_\gf{f}=\ell^2_\gf{f}(\N)$
be the space of all real sequences $x=(x_\nu)_{\nu \in \mathbb{N}}$ with
$$
   \sum_{\nu=1}^\infty \gf{f}(\nu) x_\nu^2<\infty.
$$
The space $\ell^2_\gf{f}$ is a Hilbert space with respect to the inner product
\begin{equation}\label{eq:ell2_f}
   \INNER{x}{y}_\gf{f}:=\sum_{\nu\in \N} \gf{f}(\nu) x_\nu y_\nu 
   ,\qquad
   \norm{x}_\gf{f}:=\Bigl(\sum_{\nu\in \N} \gf{f}(\nu) x_\nu^2\Bigr)^{\tfrac12} ,
\end{equation}
where $\norm{x}_\gf{f}=\sqrt{\INNER{x}{x}_\gf{f}}$ is the induced norm.
Note that $\ell^2_{\gf{f}^0}=\ell^2$.

\smallskip\noindent
\textsc{Hilbert space pair.}
The pair $(\ell^2,\ell^2_\gf{f})$ is a Hilbert space pair
by~\cite[Le.\,2.1]{Frauenfelder:2009b};
see also~\cite[Thm.\,8.1]{Frauenfelder:2021b}.
For $\nu\in\N$ let $e_\nu=(0,\dots,0,1,0,\dots)$ be the sequence
whose members are all $0$ except for member $\nu$ which is $1$.
The set of all $e_\nu$'s
\begin{equation}\label{eq:stand-basis-ell2}
   \Ee=\{e_\nu\}_{\nu\in\N}
\end{equation}
is called the \textbf{standard basis} of $\ell^2=\ell^2(\N)$.
While $\Ee$ is an orthonormal basis of~$\ell^2$,
it is still an orthogonal basis of~$\ell^2_\gf{f}$.

\smallskip\noindent
\textsc{Growth operator.}
The growth operator $T\in\Ll(\ell^2_\gf{f})$ is characterized by the identity
$
   \INNER{y}{x}_{\ell^2}
   =\INNER{y}{Tx}_{\ell^2_\gf{f}}
$
for all $x,y \in \ell^2_\gf{f}$.
Thus the growth operator $T\colon \ell^2_\gf{f}\to
\ell^2_\gf{f}$ of the pair $(\ell^2,\ell^2_\gf{f})$ is given by
\begin{equation}\label{eq:growth-T-ell2}
   T(x_\nu)=(\tfrac{x_\nu}{\gf{f}(\nu)})
   ,\qquad (x_\nu):=(x_\nu)_{\nu\in\N} .
\end{equation}
By monotonicity and unboundedness of $\gf{f}$
there exists $\nu_0$ such that for any $\nu\le \nu_0$
it holds $\tfrac{1}{\gf{f}(\nu)^2}\le \tfrac{1}{\gf{f}(1)^2}$
and for any $\nu\ge \nu_0+1$ 
it holds $\tfrac{1}{\gf{f}(\nu)}\le \gf{f}(\nu)$. Thus
$$
   \INNER{Tx}{Tx}_\gf{f}
   =\sum_{\nu=1}^\infty \tfrac{x_\nu^2}{\gf{f}(\nu)^2} \gf{f}(\nu)
   =\sum_{\nu=1}^{\nu_0} \tfrac{x_\nu^2}{\gf{f}(\nu)}
      {\color{gray}\tfrac{\gf{f}(\nu)}{\gf{f}(\nu)}}
   +\sum_{\nu=\nu_0+1}^\infty \tfrac{x_\nu^2}{\gf{f}(\nu)}
   \le\max\{\tfrac{1}{\gf{f}(1)^2},1\} \INNER{x}{x}_\gf{f} 
$$
which shows that $T$ indeed maps $\ell^2_\gf{f}$ to $\ell^2_\gf{f}$.
The elements of the standard basis~$\Ee$ are the eigenvectors $e_\nu$
of $T$ with eigenvalues $\kappa(\nu)=\tfrac{1}{\gf{f}(\nu)}$, in symbols
$$
   T e_\nu = \tfrac{1}{\gf{f}(\nu)} e_\nu
   ,\quad
   \kappa(\nu)=\tfrac{1}{\gf{f}(\nu)}
   ,\qquad
   \kappa(\nu)\ge \kappa(\nu+1)>0
   ,\quad
   \kappa(\nu)\searrow 0 .
$$
%
%
\textsc{Eigenspace core.}
We write the eigenvalues $\kappa(\nu)$ to the eigenvectors $e_\nu$
in the form of a list $\Ss(T)=(\tfrac{1}{\gf{f}(\nu)})_{\nu\in\N}$
in which can occur finite repetitions.
The eigenspace core of $T$ is then equal to
$$
   V=\bigoplus_{\nu\in\N} \R e_\nu=\R^\infty_0
   ,\qquad
   \R^\infty_0\subset \ell^2_\gf{f}\subset \ell^2.
$$
%
%
\textsc{Scale levels.}
Let $\ell^2_{\gf{f}^r}$, for $r\in\R$, consist of all sequences
$x=(x_\nu)$ such that
$$
   \norm{\xi}_{\gf{f}^r}
   =\biggl(\sum_{\nu\in\N}
   \gf{f}(\nu)^r x_\nu^2\biggr)^\frac12
   <\infty
$$
is finite.
The $r$-inner product is given by
$\INNER{x}{y}_{\gf{f}^r}=\sum_{\nu\in\N} \gf{f}(\nu)^r x_\nu y_\nu$.

\smallskip\noindent
\textsc{Real scale.}
The real Hilbert scale associated to $(\ell^2,\ell^2_\gf{f})$ is the family
$$
   \ell^{2,\gf{f}}_\R
   :=(\ell^2_{\gf{f}^r})_{r\in\R} .
$$
Because the function $\gf{f}$ is monotone increasing,
it follows that whenever ${s\le r}$ there is an inclusion 
$
   \ell^2_{\gf{f}^r} \hookrightarrow \ell^2_{\gf{f}^s}
$
of Hilbert spaces and the corresponding linear inclusion operator is bounded.
For strict inequality $s<r$, by unboundedness of~$\gf{f}$, the inclusion
operator is compact. Moreover, its image is
dense. For details we refer to~\cite[Thm.\,8.1]{Frauenfelder:2021b}
and~\cite[Sec.\,2]{Frauenfelder:2024c}.

\subsection{Scale bases}\label{sec:scale-bases}

Let $(H_0,H_1)$ be a Hilbert space pair.
Then both Hilbert spaces $H_0$ and $H_1$ are infinite dimensional,
by definition, and separable, by~\cite[Cor.\,A.5]{Frauenfelder:2024c}.

\begin{definition}\label{def:basis}
A Hilbert space is called \textbf{separable} if it contains a
countable dense subset.
An \textbf{orthonormal basis} (\textbf{ONB}) of a separable Hilbert space $H$
is a countable orthonormal subset of $H$ whose linear
span is dense in $H$.
Weakening the condition from norm $1$ to positive norm
we speak of an \textbf{orthogonal basis}.
\end{definition}

Each separable Hilbert space admits an ONB.
To see this pick a dense sequence $(v_k)_{k\in\N}$,
throw out any member $v_k$ if it is a linear combination of $v_1,\dots, v_{k-1}$,
then apply Gram-Schmidt
orthogonalization to what remains.

\begin{definition}\label{def:scale-basis}
A \textbf{scale basis} for a Hilbert space pair $(H_0,H_1)$ is an
orthonormal basis $E=\{E_\nu\}_{\nu\in\N}$ of $H_0$
that is simultaneously an orthogonal basis~of~$H_1$,
and which is ordered such that the
function
\begin{equation}\label{eq:h-E_nu}
   \gf{h}\colon \N\to (0,\infty),\quad \nu\mapsto \norm{E_\nu}_1^2
\end{equation}
is monotone increasing. Following~\cite[Thm.~A.4]{Frauenfelder:2024c}
we refer to $\gf{h}$ as the \textbf{pair growth function} of $H$.
It is automatically unbounded.
\end{definition}

\smallskip\noindent
\textsc{Existence of scale bases.}
We can construct a scale basis as follows.
We associated to $(H_0,H_1)$ an operator $T\colon H_1\to H_1$
by~(\ref{eq:Riesz}) .
For every eigenvalue $\kappa\in\spec T$
we choose an ordered $H_0$-orthonormal basis
of $V_\kappa:=\Eig_\kappa T$, notation
\begin{equation}\label{eq:E^kappa}
   E^\kappa=\{ E^\kappa_1,\dots,E^\kappa_{m_\kappa}\} .
\end{equation}
By Lemma~\ref{le:0-1-orthogonality}
the basis $E^\kappa$ of $V_\kappa$ is $H_1$-orthogonal as well and,
furthermore, all vectors have the same $H_1$-length, namely
in~(\ref{eq:0-1-lengths}) we obtained
$$
   \norm{E^\kappa_i}_1=\frac{1}{\sqrt{\kappa}}.
$$
We order the eigenvalues of $T$ decreasingly
$$
   \kappa_1>\kappa_2>\dots >0 .
$$
Now we define a function $\kappa\colon \nu\mapsto \kappa_{j(\nu)}$
that enlists the eigenvalues accounting for multiplicities.\footnote{
  E.g. if the eigenvalues $\kappa_i$ and their respective
  multiplicities $m_{\kappa_i}$ are
  $$
     \sqrt{5}>\tfrac{4}{7}>\tfrac12>\dots>0,\qquad
     2,4,m_{\kappa_3=\frac12},\dots
  $$
  the functions $\nu\mapsto j(\nu)$ and $\nu\mapsto \kappa_{j(\nu)}$
  return, respectively, the values
  $$
     \underbrace{1,1}_{m_{\kappa_1}},\underbrace{2,2,2,2}_{m_{\kappa_2}},
     3,\dots,\qquad
     \sqrt{5}, \sqrt{5}, \tfrac{4}{7}, \tfrac{4}{7}, \tfrac{4}{7}, \tfrac{4}{7}, \tfrac12,\dots
  $$
}
More precisely, 
for $\nu\in\N$ we define
$$
   j(\nu):=\min\Bigl\{j\in\N\mid \sum_{i=1}^j m_{\kappa_i}\ge \nu\Bigr\}
   ,\qquad
   \kappa(\nu):=\kappa_{j(\nu)} ,
$$
and set
$$
   E_\nu:= E^{\kappa(\nu)}_{\nu-\sum_{i=1}^{j(\nu)-1} m_{\kappa_i}}
   ,\qquad
   E:=\{E_\nu\}_{\nu\in \N} .
$$
Note that the ordered orthonormal basis $E$ of $H_0$ starts at $E_1=E_1^1$.
The pair growth function is related to the growth operator eigenvalues
$\kappa(\nu)$ by
\begin{equation}\label{eq:h-kappa}
   \gf{h}(\nu)=\frac{1}{\kappa(\nu)} .
\end{equation}

\smallskip\noindent
Next we address the question of moduli of scale bases.
For this we show the following lemma.

\begin{lemma}\label{le:base=eigenvectors}
All elements of a scale basis $E=\{E_\nu\}_{\nu\in\N}$ are $T$-eigenvectors
$$
   TE_\nu=\tfrac{1}{\Norm{E_\nu}_1^2} E_\nu
   ,\qquad
   \forall \nu\in\N.
$$
\end{lemma}

\begin{proof}
For $\nu\in\N$ write $TE_\nu=\sum_{\mu\in\N} t_{\mu\nu} E_\mu$. Then
$$
   \delta_{\rho\nu}
   \stackrel{\perp_0}{=}\INNER{E_\rho}{E_\nu}_0
   \stackrel{\text{(\ref{eq:Riesz})}}{=}\INNER{E_\rho}{TE_\nu}_1
   =\sum_{\mu\in\N} t_{\mu\nu} \INNER{E_\rho}{E_\mu}_1
   \stackrel{\perp_1}{=} t_{\rho\nu} \norm{E_\rho}_1^2
$$
where the last step uses that $\INNER{E_\rho}{E_\mu}_1$
is $0$ for $\mu\not= \rho$ and $\norm{E_\rho}_1^2$ otherwise.
\end{proof}

\medskip\noindent
\textsc{Moduli of scale bases.}

In view of Lemma~\ref{le:base=eigenvectors} all scale bases are
constructed as in~(\ref{eq:E^kappa}).
In particular, a scale basis is unique up to an action by the group
$\oplus_{\kappa\in\spec T} \O(E^\kappa)$.

\boldmath
\subsubsection{Isometry to model Hilbert $\R$-scale.}
\unboldmath

Consider a Hilbert space pair $H=(H_0,H_1)$. Let $\gf{h}$ be
a pair growth function and let $H_\R=(H_r)_{r\in\R}$ be the Hilbert
$\R$-scale associated to the pair. Any scale basis $E=\{E_\nu\}_{\nu\in\N}$
of $H$ determines, for each $r\in\R$, a Hilbert space isometry
$$
   \Psi^E_r\colon H_r\to \ell^2_{\gf{h}^r}
   ,\quad
   \xi=\sum_{\nu\in\N} \xi_\nu E_\nu\mapsto (\xi_\nu)_{\nu\in\N}
$$
by assigning to $\xi$ its coordinate sequence;
see~\cite[proof of Thm.\,A.4]{Frauenfelder:2024c}.
So
\begin{equation}\label{eq:isometric}
   \INNER{\xi}{\eta}_{H_r}=\sum_{\nu\in \N} \gf{h}(\nu) \xi_\nu \eta_\nu 
   ,\qquad
   \norm{\xi}_{H_r}=\Bigl(\sum_{\nu\in \N} \gf{h}(\nu) x_\nu^2\Bigr)^{\tfrac12} ,
\end{equation}
for all $\xi,\eta\in H_r$ and where $\gf{h}$ relates
to the growth operator eigenvalues $\kappa(\nu)$ by
$$
   \frac{1}{\kappa(\nu)}
   \stackrel{\text{(\ref{eq:h-kappa})}}{=}\gf{h}(\nu)
   \stackrel{\text{(\ref{eq:h-E_nu})}}{=}\norm{E_\nu}_1^2 .
$$

\boldmath
\subsection[Musical $\R$-scale isometry $\flat=\sharp^{-1}$ and
shift isometries]{Musical $\R$-scale isometry $\flat$ and
shift isometries}
\unboldmath

Let $H=(H_0,H_1)$ be a Hilbert space pair and $E=\{E_\nu\}_{\nu\in\N}$
a scale basis.
With~$H$ comes the growth function $\gf{h}\colon \N\to[0,\infty)$
and the Hilbert $\R$-scale $H_\R$.

\begin{definition}[Canonical $\R$-scale isometry
$\flat=\sharp^{-1}\colon H_{-r}\to  H_r^*$]
For $r\in\R$ insertion into the ${\color{red}0}$-inner product
\begin{equation}\label{eq:flat}
   \flat\colon H_{-r}\to H_r^*
   ,\quad
   \xi\mapsto \xi^\flat:=\INNER{\xi}{\cdot}_{\color{red}0}
\end{equation}
is for $\xi=\sum_{\nu\in \N} \xi_\nu E_\nu\in H_{-r}$ and
$\eta=\sum_{\nu\in\N}\eta_\nu E_\nu\in H_r$ given by the sum
$$
   (\flat \xi)\eta
   =\sum_{\nu\in \N} \xi_\nu \eta_\nu .
$$
We show that $\flat\colon H_{-r}\to H_r^*$ is an isometry.

The special case $H_0\to H_0^*$, $\xi\mapsto\INNER{\xi}{\cdot}_0$,
is the usual insertion isometry. For their common notation $\flat$
and $\sharp:=\flat^{-1}$ these are called
\textbf{musical isometries}.\footnote{
  $\flat E_\nu=E_\nu^*$ since
  $(\flat E_\nu) E_\mu=\INNER{E_\nu}{E_\mu}_0=\delta_{\nu\mu}$.
  Exactly isometries take ONB's to ONB's.
  }
\end{definition}

To see that $\flat$ is well defined,
note that $\xi\in H_{-r}$ and $\eta\in H_{r}$ implies finiteness 
$$
   \norm{\xi}_{-r}^2=\sum_{\nu\in \N} \xi_\nu^2 \gf{h}(\nu)^{-r}<\infty
   ,\qquad
   \norm{\eta}_{r}^2=\sum_{\nu\in \N} \eta_\nu^2 \gf{h}(\nu)^{r}<\infty .
$$
Thus by Cauchy-Schwarz the sum
$$
   \sum_{\nu\in \N} \xi_\nu \eta_\nu
   =\sum_{\nu\in \N} \xi_\nu \gf{h}(\nu)^{\tfrac{-r}{2}} \eta_\nu \gf{h}(\nu)^{\tfrac{r}{2}} 
   \le\Bigl(\sum_{\nu\in \N} \xi_\nu^2 \gf{h}(\nu)^{-r}\Bigr)^{\tfrac12}
   \Bigl(\sum_{\nu\in \N} \eta_\nu^2 \gf{h}(\nu)^{r}\Bigr)^{\tfrac12}
   <\infty
$$
is finite, so $\flat$ is well defined.
The fact that $\sum_{\nu\in \N} \xi_\nu \eta_\nu=0$ for every $\eta$
implies $\xi=0$ and this proves injectivity of $\flat\colon H_{-r}\to H_r^*$.
Since $H_0$ is a Hilbert space, in particular complete, $\flat\colon
H_0\to H_0^*$ is an isomorphism, as is well known.\footnote{
  Surjective: pick $\eta\in H_0^*$ non-zero, then $\ker  \eta\subset
  H_0$ is a closed subspace of co-dimension~$1$. Hence $(\ker
  \eta)^\perp=\R\hat v$ for a unit vector $\hat v\in H_0$. Now
  $\eta=\flat_0((\eta\hat v)\hat v){\color{gray}\,
  =\INNER{(\eta\hat v)\hat v}{\cdot}_0}$
  since both sides are equal on $\ker\eta=(\R\hat v)^\perp$, namely zero,
  and on $\hat v$, namely $\eta\hat v$ since $\INNER{\hat v}{\hat v}_0=1$.
  }
To see that $\flat\colon H_{-r}\to H_r^*$ is an isomorphism, in fact
an isometry, whenever $r\in\R$
consider the shift isometries introduced next,
then apply Lemma~\ref{le:comp-flat_r}.

\begin{definition}[Shift isometries]\label{def:shift-isos}
Given reals $r,s\in\R$, we define
$$
   \phi_r^s\colon H_r\to H_s
   ,\quad
   \xi\mapsto \sum_{\nu\in \N}\xi_\nu \gf{h}(\nu)^{\frac{r-s}{2}} E_\nu
$$
where $\xi=\sum_{\nu\in \N} \xi_\nu E_\nu$.
\end{definition}

The maps $\phi_r^s$ are norm preserving with inverse $(\phi_r^s)^{-1}=\phi_s^r$.
For $\xi\in H_r$ we compute
\begin{equation*}
\begin{split}
   \norm{\phi_r^s\xi}_s^2
   &=\sum_{\nu\in \N} \gf{h}(\nu)^s \left(\xi_\nu \gf{h}(\nu)^{\frac{r-s}{2}}\right)^2
   =\sum_{\nu\in \N} \xi_\nu^2 \gf{h}(\nu)^r
   =\norm{\xi}_r^2
\end{split}
\end{equation*}
which proves norm preservation. But this implies inner product preservation
by the polarization identity. Thus adjoint is inverse:
$(\phi_r^s)^*=(\phi_r^s)^{-1}=\phi_s^r$. In particular, the adjoint is an
isometry as well.

\begin{lemma}\label{le:comp-flat_r}
$\flat=(\phi_r^0)^*\flat\phi_{-r}^0\colon H_{-r}\to H_r^*$ is composed of
isometries $\forall r\in\R$.
\end{lemma}

\begin{proof}
The maps on the right hand side are isometries, as was shown above.
Given $\xi\in H_{-r}$ and $\eta\in H_r$, use the characterization of
the adjoint to~get
\begin{equation*}
\begin{split}
   ({\phi_r^0}^*\flat\phi_{-r}^0\xi)\eta
   &=(\flat\phi_{-r}^0\xi) \phi_{r}^0\eta\\
   &\stackrel{{\color{gray}2}}{=}
   \INNER{{\textstyle\sum_{\nu\in \N}}\xi_\nu\gf{h}(\nu)^{-\frac{r}{2}}E_\nu}
   {{\textstyle\sum_{\mu\in \N}}\eta_\mu\gf{h}(\mu)^{\frac{r}{2}}E_\mu}_0\\
   &=\sum_{\nu,\mu\in \N}\xi_\nu\gf{h}(\nu)^{-\frac{r}{2}}\eta_\mu\gf{h}(\mu)^{\frac{r}{2}}
   \underbrace{\INNER{E_\nu}{E_\mu}_0}_{\delta_{\nu\mu}}\\
   &=\sum_{\nu\in\N}\xi_\nu\gf{h}(\nu)^{-\frac{r}{2}}\eta_\nu\gf{h}(\nu)^{\frac{r}{2}}\\
   &=\sum_{\nu\in\N}\xi_\nu\eta_\nu\\
   &=(\flat\xi)\eta
\end{split}
\end{equation*}
where in equality two we used the definition of $\phi_{-r}^0\xi$ and
of $\phi_r^0\eta$.
\end{proof}

\begin{lemma}[$\A^*\simeq \A$]\label{le:A*-A}
Assume that $\A\colon H_1\to H_0$ is a symmetric isometry.
Then the composition of $\A^*\colon H_0^*\to H_1^*$ with the four isometries
$$
   \A\colon
   H_1
   \stackrel{\phi_1^0}{\longrightarrow}
   H_0
   \stackrel{\flat}{\longrightarrow}
   H_0^*
   \stackrel{\A^*}{\longrightarrow}
   H_1^*
   \stackrel{\flat^{-1}}{\longrightarrow}
   H_{-1}
   \stackrel{\phi_{-1}^0}{\longrightarrow}
   H_0
$$
is equal to $\A$. 
\end{lemma}

\begin{proof}
The matrix of $\A$ for an eigenvector orthonormal basis
$\Vv(\A)$ of $H_0$, (\ref{eq:ONB-A}), is diagonal.
Let $\xi\in H_1$. To show equality
$\flat\phi_0^{-1} \A\xi=\A^*\flat\phi_1^0\xi\in H_1^*$, apply both
sides to $\eta\in H_1$. By linearity, basis elements $\xi=v_\nu$
and $\eta=v_\mu$ suffice. We get
\begin{equation*}
\begin{split}
   (\flat\phi_0^{-1} A v_\nu)v_\mu
   &\stackrel{\flat}{=}\INNER{\phi_0^{-1} a_\nu v_\nu}{v_\mu}_0\\
   &=
   \INNER{a_\nu\abs{a_\nu} v_\nu}{v_\mu}_0\\
   &=a_\nu\abs{a_\nu} \delta_{\nu\mu}
\end{split}
\end{equation*}
and
\begin{equation*}
\begin{split}
   (A^*\flat\phi_1^0v_\nu) v_\mu
   &=(\flat\phi_1^0 v_\nu) A v_\mu\\
   &\stackrel{\flat}{=}\INNER{\phi_1^0 v_\nu}{A v_\mu}_0\\
   &=\INNER{\abs{a_\nu}v_\nu}{a_\mu v_\mu}_0\\
   &=\abs{a_\nu} a_\mu \delta_{\nu\mu} .
\end{split}
\end{equation*}
This proves Lemma~\ref{le:A*-A}.
\end{proof}

\section{Quantitative invertibility}

In the proof of Theorem~\ref{thm:Rabier}
we will use the following well-known lemma
which shows, also quantitatively,
that invertibility is an open condition.

\begin{lemma}[Quantitative invertibility]\label{le:invertibility-open}
Given Banach spaces $X$ and $Y$, suppose the operator $T\in\Ll(X,Y)$
is invertible and $P\in\Ll(X,Y)$ is small in the sense that
$\norm{P}< 1/\norm{T^{-1}}$.
Then $T+P$ is invertible as well with bound
$$
   \norm{(T+P)^{-1}}
   \le
   \frac{\norm{T^{-1}}}{1-\norm{T^{-1}}\cdot\norm{P}}
$$
where all norms are operator norms.
\end{lemma}

\begin{proof}
We define $S:=\Id-T^{-1}(T+P)$ and we estimate
\begin{equation}\label{eq:S77}
   \norm{S}
   =\norm{\Id -T^{-1}(T+P)}
   =\norm{T^{-1}P}
   \le \norm{T^{-1}}\norm{P}
   <1.
\end{equation}
Hence $T^{-1}(T+P)=\Id -S$ is invertible with the help of the Neumann
series
$$
   (\Id-S)^{-1}
   = \sum_{n=0}^\infty S^n
$$
whose norm we can estimate via the geometric series
$$
   \norm{(\Id-S)^{-1}}
   \le \sum_{n=0}^\infty \norm{S}^n
   =\frac{1}{1-\norm{S}}.
$$
An inverse of $T^{-1} (T+P)$ is
$
   (T^{-1}(T+P))^{-1}
   =(\Id-S)^{-1}
$
and bounded by   
$$
   \norm{(T^{-1}(T+P))^{-1}}
   =\norm{(\Id-S)^{-1}}
   \le \frac{1}{1-\norm{S}}
   \stackrel{\text{(\ref{eq:S77})}}{\le}
   \frac{1}{1-\norm{T^{-1}}\norm{P}} .
$$
Therefore $T+P=T(T^{-1}(T+P))$ is invertible and the inverse
$(T+P)^{-1}=(T^{-1}(T+P))^{-1} T^{-1}$ is bounded by
$
   \norm{(T+P)^{-1}}
   \le \norm{T^{-1}}/(1-\norm{T^{-1}}\norm{P})
$.
\end{proof}

\newpage 
\section{Evaluation map \boldmath$P_1\to H_{1/2}$}
\label{sec:evaluation}

Let $H=(H_0,H_1)$ be a Hilbert space pair.
Let $\gf{h}\ge 1$
be a growth function representing the pair growth type.
For the time interval $I=[0,1]$
we define the path space $P_1=P_1(I)$ by~(\ref{eq:P_1-I}).
Let $E=\{E_\nu\}_\nu\in\N$ be a scale basis of~$H$;
see Appendix~\ref{sec:scale-bases}.

\begin{proposition}\label{prop:tr1-NEW}
Let $x \in W^{1,2}([0,1],H_0) \cap L^2([0,1],H_1)$, then $x(0) \in H_{1/2}$.
\end{proposition}

\begin{proof}
Writing $x=\sum_\nu x_\nu E_\nu$
we estimate for $s \in\Big[0,\tfrac{1}{\sqrt{{\gf{h}}(\nu)}}\Big]$
the initial point
\begin{eqnarray*}
|x_\nu(0)| &\leq& |x_\nu(s)|+\int_0^s|\partial_t x_\nu(t)|dt\\
&\leq& |x_\nu(s)|+\int_0^{\frac{1}{\sqrt{{\gf{h}}(\nu)}}}|\partial_t x_\nu(t)|dt\\
&\leq& |x_\nu(s)|+\frac{1}{{\gf{h}}(\nu)^{1/4}}||\partial_t x_\nu||_{L^2} 
\end{eqnarray*}
where the last step is by H\"older's inequality.
Therefore
\begin{eqnarray*}
x_\nu(0)^2 &\leq& \bigg(|x_\nu(s)|+\frac{1}{{\gf{h}}(\nu)^{1/4}}||\partial_t x_\nu||_{L^2}\bigg)^2\\
&\leq& 2x_\nu(s)^2+\frac{2}{\sqrt{{\gf{h}}(\nu)}}||\partial_t x_\nu||^2_{L^2} .
\end{eqnarray*}
Taking advantage of this estimate in step four we obtain that
\begin{eqnarray*}
||x||^2_{L^2(H_1)}&=&\int_0^1||x(s)||^2_{\gf{h}} ds\\
&=&\sum_{\nu=1}^\infty \int_0^1 {\gf{h}}(\nu)x_\nu(s)^2ds\\
&\geq&\sum_{\nu=1}^\infty \int_0^{\frac{1}{\sqrt{{\gf{h}}(\nu)}}} {\gf{h}}(\nu)x_\nu(s)^2ds\\
&\geq&\frac{1}{2}\sum_{\nu=1}^\infty \int_0^{{\gf{h}}(\nu)^{\color{cyan} -1/2}} 
{\gf{h}}(\nu)^{\color{cyan} 1} x_\nu(0)^2 ds-
\sum_{\nu=1}^\infty \int_0^{\frac{1}{\sqrt{{\gf{h}}(\nu)}}} \sqrt{{\gf{h}}(\nu)}||\partial_t x_\nu||^2_{L^2}ds\\
&=&\frac{1}{2}\sum_{\nu=1}^\infty  {\gf{h}}(\nu)^{\color{cyan} 1-1/2} x_\nu(0)^2 -
\sum_{\nu=1}^\infty ||\partial_t x_\nu||^2_{L^2}\\
&\geq&\frac{1}{2}||x(0)||^2_{H_{\color{cyan} 1/2}}-||x||^2_{W^{1,2}(H_0)} .
\end{eqnarray*}
Hence
$$||x(0)||_{H_{1/2}}\leq \sqrt{2}||x||_{L^2(H_1) \cap W^{1,2}(H_0)}.$$
This completes the proof of Proposition~\ref{prop:tr1-NEW}. 
\end{proof}

\begin{definition}
By Proposition~\ref{prop:tr1-NEW} we obtain well defined
evaluation maps
$$
   \mathrm{ev} \colon P_1=W^{1,2}([0,1],H_0) \cap L^2([0,1],H_1) \to
   H_{1/2}, \quad x \mapsto x(0)
$$
and
$$
   \mathrm{Ev}\colon P_1\to H_{1/2} \times H_{1/2}
   ,\quad
   x\mapsto (x(0),x(1)) .
$$
The evaluation maps are linear continuous maps between
Hilbert spaces. 
\end{definition}

\begin{proposition}\label{prop:ev-onto}
The evaluation map $\mathrm{ev}\colon P_1\to H_{1/2}$
is surjective. 
\end{proposition}

\begin{proof}
Suppose that $x^0=(x^0_\nu)_{\nu \in \mathbb{N}} \in H_{1/2}$. 
Define $x_\nu \in C^\infty([0,1],\mathbb{R})$ by
$$x_\nu(s)=e^{-\sqrt{{\gf{h}}(\nu)}s}x_\nu^0, \quad s \in [0,1].$$
Note that
$$x_\nu(0)=x_\nu^0$$
so that if we set $x=(x_\nu)_{\nu \in \mathbb{N}}$ we have
$$\mathrm{ev}(x)=x^0.$$
Therefore in order to prove the proposition it suffices to show that
$$x \in  W^{1,2}([0,1],H_0) \cap L^2([0,1],H_1).$$
In order to achieve this we estimate
\begin{eqnarray*}
||x||^2_{W^{1,2}(H_0) \cap L^2(H_1)}&=&||x||_{L^2(H_1)}^2+||x||^2_{W^{1,2}(H_0)}\\
&=&\sum_{\nu=1}^\infty \Bigg(\int_0^1 {\gf{h}}(\nu) x_\nu^2(s)ds+\int_0^1 \partial_s x_\nu(s)^2 ds+\int_0^1 x_\nu(s)^2 ds\Bigg)\\
&=&\sum_{\nu=1}^\infty \int_0^1 \big(2{\gf{h}}(\nu)+1) e^{-2\sqrt{{\gf{h}}(\nu)}s} \big(x_\nu^0\big)^2 ds\\
&=&-\sum_{\nu=1}^\infty \frac{2 {\gf{h}}(\nu)+1}{2\sqrt{{\gf{h}}(\nu)}}\big(x_\nu^0\big)^2  e^{-2\sqrt{{\gf{h}}(\nu)}s}\bigg|_0^1\\
&=&\sum_{\nu=1}^\infty \frac{2 {\gf{h}}(\nu)+1}{2\sqrt{{\gf{h}}(\nu)}}\big(x_\nu^0\big)^2\bigg(1-e^{-2\sqrt{{\gf{h}}(\nu)}}\bigg)\\
&\leq&\sum_{\nu=1}^\infty2 \sqrt{{\gf{h}}(\nu)}\big(x_\nu^0\big)^2\\
&=&2||x||^2_{H_{1/2}} .
\end{eqnarray*}
This finishes the proof of the proposition.
\end{proof}

\begin{corollary}\label{cor:Ev-onto}
The evaluation map $\mathrm{Ev}\colon P_1\to H_{1/2}\times H_{1/2}$
is surjective. 
\end{corollary}

\begin{proof}
Given $x^0,x^1\in H_{1/2}\times H_{1/2}$,
there exist, by Proposition~\ref{prop:ev-onto}, paths
$y^0,y^1\in P_1$ such that $y^0(0)=x^0$ and $y^1(1)=x^1$.
Pick cutoff functions $\beta_0,\beta_1\in C^\infty([0,1],[0,1])$
such that $\beta_0(0)=1$ and $\beta_0\equiv 0$ on $[1/2,1]$
and $\beta_1(1)=1$ and $\beta_1\equiv 0$ on $[0,1/2]$.
Then the combination $y:=\beta_0 y^0+\beta_1 y^1$ still lies in $P_1$
and $y(0)=y^0(0)=x^0$ and $y(1)=y^1(1)=x^1$.
\end{proof}

\newpage 
\section[Invariance of Fredholm index (non-constant target space)]
{Invariance of the Fredholm index}\label{sec:Fredholm-index}

\subsection{Varying target space}

Assume that $X$ and $Y$ are Hilbert spaces
and $\Dd_r\colon X\to Y$ for $r\in[0,1]$ is a continuous family
of bounded linear maps.
Assume further that $p_r\in\Ll(Y)$ is a family of
projections, orthogonal or not, depending continuously on $r\in[0,1]$.
Since $p_r$ is a projection ($p_r p_r=p_r$) its
image is equal to its fixed point set which is closed by continuity.
Hence $\im p_r$ is a closed subspace of $Y$.
We abbreviate the composition by
$$
   \mathfrak{D}_r\colon X
   \stackrel{\Dd_r}{\longrightarrow} Y
   \stackrel{p_r}{\longrightarrow}\im p_r\subset Y .
$$

\begin{theorem}\label{thm:Fredholm-index-constant}
Assume that $\mathfrak{D}_r\colon X\to \im p_r$ is Fredholm for any
$r\in[0,1]$, then its Fredholm index is independent of $r$.
\end{theorem}

\begin{proof}
The case $p_r=\Id_Y$ is well known.
We first discuss that case as warmup.

\medskip\noindent
\textbf{Case~1:} $p_r\equiv\Id_Y$\textbf{.}
The Fredholm index of $\Dd_r\colon X\to Y$ is independent of $r$.

\begin{proof}[Proof of Case~1]
In this case $\mathfrak{D}_r=\Dd_r\colon X\to Y$ is a Fredholm
operator between fixed Hilbert spaces.
For fixed $r,s\in[0,1]$ we abbreviate
$$
   D:=\Dd_r\colon X\to Y
   ,\qquad
   Q:=\Dd_s-\Dd_r \colon X\to Y .
$$
Abbreviate $X_0:=\ker D$ and $Y_1:=\im D$ and decompose orthogonally
$$
   X=\underbrace{X_0}_{\ker D}\stackrel{\perp}{\oplus} X_1
   ,\qquad
   Y=Y_0\stackrel{\perp}{\oplus}  \underbrace{Y_1}_{\im D} .
$$
Let $D_{ij}\colon X_i\to Y_j$ denote the restriction of $D$ to $X_i$
followed by projection onto $Y_j$, and similarly for $Q$.
Note that $D$ is of the form
$$
   D=
   \begin{pmatrix}
     D_{00}&D_{01}\\D_{10}&D_{11}
   \end{pmatrix}
   =
   \begin{pmatrix}
     0&0\\0&D_{11}
   \end{pmatrix}
   \colon X_0\oplus X_1\to Y_0\oplus Y_1
$$
where $D_{11}\colon X_1\to Y_1$ is bijective, hence an isomorphism by
the open mapping theorem.
The operator $Q$ is of the form
$$
   Q=
   \begin{pmatrix}
     Q_{00}&Q_{01}\\Q_{10}&Q_{11}
   \end{pmatrix}
   \colon X_0\oplus X_1\to Y_0\oplus Y_1
$$
If $s$ is close to $r$, then $Q$ is close to the zero operator,
and so is $Q_{11}$. So by openness of invertibility
$D_{11}+Q_{11}\colon X_1\to Y_1$ is still an isomorphism.
The linear map between finite dimensional vector spaces
$$
   F:=Q_{00}-Q_{01}(D_{11}+Q_{11})^{-1} Q_{10}\colon X_0=\ker D\to
   Y_0=\coker D
$$
is Fredholm and its index is the dimension difference
of domain and target
$$
   \INDEX{F}=\dim X_0 -\dim Y_0
   =\INDEX{D} .
$$

\smallskip\noindent
\textsc{Claim 1.} $\dim \ker (D+Q)=\dim \ker F$.

\smallskip\noindent
Write $x\in \ker (D+Q)\subset X_0\oplus X_1$ uniquely in the form
$x=x_0+x_1$ where $x_i\in X_i$. Then we get two equations in the form
$$
   \begin{pmatrix}
     0\\0
   \end{pmatrix}
   =(D+Q)x
   =
   \begin{pmatrix}
     Q_{00}&Q_{01}\\Q_{10}&D_{11}+Q_{11}
   \end{pmatrix}
   \begin{pmatrix}
     x_0\\x_1
   \end{pmatrix}
   =
   \begin{pmatrix}
     Q_{00} x_0+ Q_{01}x_1\\
     Q_{10} x_0 +(D_{11}+Q_{11}) x_1
   \end{pmatrix} .
$$
The second equation tells that
\begin{equation}\label{eq:x_1}
   x_1=-(D_{11}+Q_{11})^{-1} Q_{10} x_0 .
\end{equation}
Insert this into equation one to get
$
   0=Q_{00}x_0 -Q_{01}(D_{11}+Q_{11})^{-1} Q_{10} x_0
   = Fx_0
$.
Consequently projection to the $X_0$-component is well defined as a map
$$
   \pi_0 \colon {\color{gray} X_0\oplus X_1\subset\,} \ker (D+Q)\to
   \ker F{\color{gray}\,\subset X_0}
   ,\quad
   x=x_0+x_1\mapsto x_0 .
$$
We show that $\pi_0$ is an isomorphism by constructing
an inverse, the candidate~is
$$
   \tau\colon \ker F\to \ker (D+Q),\quad 
   x_0\mapsto (x_0, -(D_{11}+Q_{11})^{-1} Q_{10} x_0).
$$
The image of $\tau$ lies in the kernel of $D+Q$, indeed
$$
   \begin{pmatrix}
     Q_{00}&Q_{01}\\Q_{10}&D_{11}+Q_{11}
   \end{pmatrix}
   \begin{pmatrix}
     x_0\\-(D_{11}+Q_{11})^{-1} Q_{10} x_0
   \end{pmatrix}
   =
   \begin{pmatrix}
     F x_0\\
     0
   \end{pmatrix} 
   =
   \begin{pmatrix}
     0\\
     0
   \end{pmatrix}  .
$$
Clearly $\pi_0\tau=\Id$.
Vice versa $\tau\pi_0=\Id$ holds by~(\ref{eq:x_1}).
This proves Claim~1.

\smallskip\noindent
\textsc{Claim 2.} $\dim \coker (D+Q)=\dim \coker F$.

\smallskip\noindent
This amounts to prove that the dimensions of the orthogonal
complements $(\im D+Q)^\perp$ and $(\im F)^\perp$ coincide.
\\
Suppose that $y=y_0+y_1\in Y_0\oplus Y_1$ is element of $(\im (D+Q))^\perp$,
equivalently
\begin{equation}\label{eq:hghg77}
\begin{split}
   0
   =
   \INNER{Q_{00} x_0+ Q_{01}x_1}{y_0}
   +\INNER{Q_{10} x_0 +(D_{11}+Q_{11}) x_1}{y_1}
\end{split}
\end{equation}
for every $x=x_0+x_1\in X_0\oplus X_1$. We take two particular choices.

Firstly, for the choice $x_0=0$ condition~(\ref{eq:hghg77}) reduces to
$$
   0=\inner{Q_{01}x_1}{y_0}+\inner{(D_{11}+Q_{11})x_1}{y_1}
   =\inner{x_1}{Q_{01}^*y_0+(D_{11}+Q_{11})^* y_1}
$$
for every $x_1\in X_1$.
By non-degeneracy of the inner product this means that
\begin{equation}\label{eq:y_1}
   y_1=-(D_{11}+Q_{11})^{*^{-1}} Q_{01}^* y_0
\end{equation}
whenever $y_0+y_1\in Y_0\oplus Y_1$ is element of $(\im (D+Q))^\perp$.

Secondly, in~(\ref{eq:hghg77}) choose $x_1$ according to~(\ref{eq:x_1}).
Then the first factor in the first inner product is $Fx_0$ and in the
second inner product the first factor is $0$, thus what remains is
$0=\INNER{F x_0}{y_0}_Y$ for every $x_0\in X_0$.
Hence $y_0\perp \im F$ and therefore projection to the $Y_0$-component
is well defined as a map
$$
   \Pi_0\colon {\color{gray}Y_0\oplus Y_1\supset\,} (\im (D+Q))^\perp
   \to (\im F)^\perp{\color{gray}\,\subset Y_0}
   ,\quad
   y_0+y_1\mapsto y_0 .
$$
We show that $\Pi_0$ is an isomorphism by constructing
an inverse, the candidate~is
$$
   \Tt\colon (\im F)^\perp\to (\im (D+Q))^\perp,\quad 
   y_0\mapsto y_0+y_1
$$
where $y_1$ is given by~(\ref{eq:y_1}).
To see that the image of $\Tt$ lies in $(\im (D+Q))^\perp$, insert $\Tt y_0=y_0+y_1$
into the right hand side of condition~(\ref{eq:hghg77}) and note that
\begin{equation*}
\begin{split}
   &\INNER{Q_{00} x_0}{y_0}_{Y}
   +\underline{\INNER{Q_{01}x_1}{y_0}_{Y}}
   +\inner{Q_{10} x_0}{-(D_{11}+Q_{11})^{*^{-1}} Q_{01}^* y_0}_{Y}\\
   &\quad
   +\underline{\inner{(D_{11}+Q_{11}) x_1}{-(D_{11}+Q_{11})^{*^{-1}} Q_{01}^* y_0}_{Y}}
\\
   &=\inner{Q_{00} x_0-Q_{01}(D_{11}+Q_{11})^{-1} Q_{10} x_0}{y_0}_Y\\
   &=\inner{F x_0}{y_0}_Y\\
   &=0
\end{split}
\end{equation*}
indeed vanishes for every $x=x_0+x_1\in X_0\oplus X_1$.
This proves that $\Tt y_0\in (\im (D+Q))^\perp$.
In the calculation the two underlined terms canceled each other
and the last equality is due to the domain of $\Tt$,
namely $y_0\in (\im F)^\perp$.
\\
Clearly $\Pi_0\Tt=\Id$.
Vice versa $\Tt\Pi_0=\Id$ holds by~(\ref{eq:y_1}).
This proves Claim~2.

\medskip\noindent
We prove Claim~1.
By definition of $D$ and $Q$ the above discussion shows that
\begin{equation*}
\begin{split}
   \INDEX{\Dd_s}
   &=\INDEX(D+Q)\\
   &=\dim\ker(D+Q)-\dim\coker(D+Q)\\
   &=\dim\ker{F}-\dim\coker{F}\\
   &=\INDEX{F}\\
   &=\INDEX{D}\\
   &=\INDEX{\Dd_r}
\end{split}
\end{equation*}
for all $s,r\in[0,1]$ sufficiently close.
This proves the well known Case~1.
\end{proof}

\noindent
\textbf{Case~2:} General\textbf{.}
The Fredholm index of the composed operator
$\mathfrak{D}_r:=p_r\circ \Dd_r\colon X\to Y\to \im p_r$
is independent of $r\in[0,1]$.

\begin{proof}[Proof of Case~2]
We reduce the proof of Case~2 to Case~1 via Step~1:

\smallskip\noindent
\textsc{Step~1.} For any $r\in[0,1]$ there is $\eps>0$
such that $p_r|_{\im{p_s}}\colon\im{p_s}\to\im{p_r}$
is an isomorphism for every $s\in (r-\eps,r+\eps)\cap [0,1]$.

\smallskip\noindent
To see this, given $r\in[0,1]$, by continuity of projections
we choose $\eps>0$ sufficiently small such that
$\norm{p_r-p_s}_{\Ll(Y)}\le\min\{1/4 \norm{p_r}_{\Ll(Y)},\tfrac12\}$
for every $s\in (r-\eps,r+\eps)\cap [0,1]$.
Now, for any such $s$, we estimate
\begin{equation*}
\begin{split}
   \norm{p_r\circ p_s|_{\im{p_r}}-\1_{\im{p_r}}}_{\Ll(\im{p_r})}
   &=\norm{p_r\circ p_s|_{\im{p_r}}-p_r\circ p_r|_{\im{p_r}}}_{\Ll(\im{p_r})}\\
   &=\norm{p_r\left(p_s|_{\im{p_r}}-p_r|_{\im{p_r}}\right)}_{\Ll(\im{p_r})}\\
   &\le\norm{p_r}_{\Ll(Y)}\cdot\norm{p_s-p_r}_{\Ll(Y)} \\
   &\le\tfrac{1}{4} .
\end{split}
\end{equation*}
Analogously we get the estimate
\begin{equation*}
\begin{split}
   \norm{p_s\circ p_r|_{\im{p_s}}-\1|_{\im{p_s}}}_{\Ll(\im{p_s})}
   &=\norm{p_s\circ p_r|_{\im{p_s}}-p_s\circ p_s|_{\im{p_s}}}_{\Ll(\im{p_s})}\\
   &=\norm{p_s\left(p_r|_{\im{p_s}}-p_s|_{\im{p_s}}\right)}_{\Ll(\im{p_s})}\\
   &\le\norm{p_s-p_r+p_r}_{\Ll(Y)}\cdot\norm{p_r-p_s}_{\Ll(Y)}\\\
   &\le\norm{p_s-p_r}_{\Ll(Y)}^2
   +\norm{p_r}_{\Ll(Y)}\cdot\norm{p_r-p_s}_{\Ll(Y)} \\
   &\le\tfrac{1}{4} + \tfrac{1}{4} .
\end{split}
\end{equation*}
This proves that both compositions
$$
   p_r\circ p_s|_{\im{p_r}}\in\Ll(\im{p_r})
   ,\qquad
   p_s\circ p_r|_{\im{p_s}}\in\Ll(\im{p_s})
   ,
$$
are invertible. Hence $p_r|_{\im{p_s}}\colon\im{p_s}\to\im{p_r}$
is surjective by the first composition and injective by the second,
thus an isomorphism by the open mapping theorem.
This proves Step~1.

\smallskip\noindent
\textsc{Step~2.} 
We prove Case~2.

\smallskip\noindent
Fix $r\in[0,1]$.
We consider the family of operators, continuous in $s\in[0,1]$,
between fixed Hilbert spaces
$$
   p_r\circ\mathfrak{D}_s
   \colon X \to \im{p_s} \to \im{p_r} .
$$
Let $\eps>0$ be as in Step~1.
Because for $s\in(s-\eps,s+\eps)\cap [0,1]$ the projection
$p_r|_{\im{p_s}}\colon\im{p_s}\to\im{p_r}$ is an isomorphism,
we conclude that $p_r\circ\mathfrak{D}_s$
is a Fredholm operator\footnote{
  same kernel, isomorphism preserves closedness of image and
  dimension of cokernel
  }
satisfying
$$
   \INDEX{(p_r\circ\mathfrak{D}_s)}=\INDEX{\mathfrak{D}_s} .
$$
By Case~1 we further have
$$
   \INDEX{(p_r\circ\mathfrak{D}_s)}
   =\INDEX{(p_r\circ\mathfrak{D}_r)}
$$
for every $s\in (r-\eps,r+\eps)\cap [0,1]$.
Since $p_r\circ\mathfrak{D}_r=\mathfrak{D}_r$,
we combine the two index equalities to obtain
$\INDEX{\mathfrak{D}_s}=\INDEX{\mathfrak{D}_r}$
for every $s\in (r-\eps,r+\eps)\cap [0,1]$.
This proves that the index is locally constant
and, since $[0,1]$ is connected, we obtain the the index is globally
constant on $[0,1]$.

This proves the Case~2.
\end{proof}
This concludes the proof of Theorem~\ref{thm:Fredholm-index-constant}.
\end{proof}

\subsection{Composition}

\begin{theorem}[Composition]\label{thm:Fredholm-compose-two}
Let $X,Y,Z$ be Banach spaces.
\\
i)~Let $S\colon X\to Y$ and $T\colon Y\to Z$ be Fredholm
operators between Banach spaces, then the composition $R\circ S$ is
Fredholm and
$$
   \INDEX{R\circ S}=\INDEX{R}+\INDEX{S} .
$$
ii) If both $S\colon X\to Y$ and $T\colon Y\to Z$ are bounded linear
maps with finite dimensional kernel and closed range, then
the above index formula is still valid, although with values in
$\Z\cup\{-\infty\}$.
\end{theorem}

\begin{proof}
See e.g.~\cite[\S 16 Thm.\,5 and Thm.\,12]{Muller:2007a}.
\end{proof}

\begin{theorem}\label{thm:Fredholm-compose}
Let $D\colon X\to Y$ be a bounded linear operator between Hilbert
spaces. Let $p\colon Y\to Y$ be a projection
whose image $Z:=\im p$ is of finite co-dimension.
Then the following is true.
The operator $D\colon X\to Y$ is Fredholm iff
$D_p:=p\circ D\colon X\to Z$ is Fredholm as a map to $Z$
and in this case the indices are related by
$\INDEX{D}=\INDEX{D_p}-\codim{Z}$.
\end{theorem}

\begin{proof}
As a map $p\colon Y\to Z$ is Fredholm
and $\INDEX{p}=\dim\ker{p}=\codim Z$.

\medskip
\noindent
\textbf{Case~1.} $D\colon X\to Y$ is Fredholm.

\begin{proof}
The composition of Fredholm operators $D_p=p \circ D\colon X\to Y\to Z$
is Fredholm, by Theorem~\ref{thm:Fredholm-compose-two}, and
$
   \INDEX{D_p}
   =\codim Z +\INDEX{D} 
$.
\end{proof}

\medskip
\noindent
\textbf{Case~2.} $p\circ D\colon X\to Y\to Z$ is Fredholm.

\begin{proof}
a) The kernel of $D$ is finite dimensional:
True since $\ker{D}\subset \ker{(p\circ D)}$.
b) The image of $D$ is closed: It is the pre-image
under the continuous map $p$ of the, by assumption closed,
image of $p\circ D$, in symbols
$
   \im{D}=p^{-1}\left(\im{(p\circ D)}\right)
$.
c) The co-kernel of $D$ is finite dimensional:
By a) and b) part ii) of Theorem~\ref{thm:Fredholm-compose-two}
applies and its index formula yields that
$$
   \dim{\coker{D}}
   =\codim{Z}+\dim{\ker{D}} +\dim{\coker{(p\circ D)}}
   -\dim{\ker{(p\circ D)}}.
$$
But the right hand side is finite by a) and assumption.
\end{proof}
This concludes the proof of Theorem~\ref{thm:Fredholm-compose}.
\end{proof}

\subsection{Varying domain}

\begin{theorem}\label{thm:Fred-var-domain}
Let $X,Y,Z$ be Hilbert spaces and $D\in\Ll(X,Y)$.
Suppose that $[0,1]\ni r\mapsto F_r\in\Ll(X,Z)$
is a continuous family of linear surjections.
Then the following is true.
If, for each $r\in[0,1]$, the restriction of $D$ to $\ker F_r$, notation
$$
   D_r:=D|\colon X\supset V_r\to Y
   ,\qquad
   V_r:=\ker{F_r} ,
$$
is a semi-Fredholm operator, then the semi-Fredholm index\footnote{
  The semi-Fredholm index
   $\INDEX{D_r}:=\dim\ker{D_r}-\dim\coker{D_r}$
   takes values in $\{-\infty\}\cup\Z$.
  }
of $D_r$ does not depend on $r$.
\end{theorem}

\begin{proof}
The proof is in two Steps.

\smallskip
\noindent
\textbf{Step~1.} (Kernel of $F_r$ as a graph)\textbf{.}
For $r$ near zero $V_r:=\ker{F_r}$ is the graph of
$$
   T_r:=(F_r|_{V_0^\perp})^{-1}(F_0-F_r) \colon V_0 \to V_0^\perp
$$
and $T_r\to 0$ in $\Ll(V_0,V_0^\perp)$, as $r\to 0$.

\begin{proof}
Given $x\in V_0$, we shall determine $y=y(x,r)$
such that a) $y\in V_0^\perp$ and b) $F_r(x+y)=0$.
\\
By b) and since $x\in \ker{F_0}$ we get
$
   0=F_r(x+y)=F_r x+F_r y=(F_r-F_0) x+F_r y
$.
Hence $F_r y=(F_0-F_r) x$.
Since $F_0$ is onto, it holds that the restriction to a complement of
the kernel $F_0|_{V_0^\perp} \colon V_0^\perp \to Z$ is an isomorphism.
Since the map $r\mapsto F_r\in\Ll(X,Z)$ is continuous, so is in
particular $r\mapsto F_r|_{V_0^\perp}\in\Ll(V_0^\perp,Z)$,
Since the condition to be an isomorphism is an open property, each 
$$
   F_r|_{V_0^\perp}\colon V_0^\perp\stackrel{\simeq}{\longrightarrow} Z
   ,\quad
   \text{$r\ge 0$ small},
$$
is still an isomorphism.
\\           
Consequently $y$ is given in the form
$y=(F_r|_{V_0^\perp})^{-1}(F_0-F_r) x$.
We abbreviate
$$
   T_r:=(F_r|_{V_0^\perp})^{-1}(F_0-F_r) \colon V_0 \to V_0^\perp .
$$
Then $V_r=\graph{\,T_r}$.
The linear map $(F_r|_{V_0^\perp})^{-1}\colon Z \to V_0^\perp$ is bounded,
uniformly in $r\ge 0$ small.
Hence, since $r\mapsto F_r$ is continuous,
it holds that $T_r$ converges to the zero operator
in $\Ll(V_0,V_0^\perp)$, as $r\to 0$.
\end{proof}

\smallskip
\noindent
\textbf{Step~2.} 
We prove the theorem.

\begin{proof}
We show that the index is locally constant.
Since the interval $[0,1]$ is connected this implies that the index is
constant on the whole interval. To simplify notation
we discuss local constancy at $r=0$.

By Step~1 we can write for small $r\ge 0$ the subspace
$V_r$ of $X$ as the graph of $T_r$.
The graph map is the isomorphism $\Gamma_r\colon V_0\to V_r$
defined by $x\mapsto (x,T_r x)$.
We further set $D_r^0:=D_r\circ \Gamma_r\colon V_0\to V_r\to Y$.
Since $D_r$ is a semi-Fredholm operator by hypothesis and $\Gamma_r$
is an isomorphism it follows that $D_r^0$ is a semi-Fredholm operator of
the same index, namely $\INDEX{D_r^0}=\INDEX{D_r}$.
\\
Note that $\Gamma_0=\Id_{V_0}$, hence $D_0^0=D_0$.
Since $T_r\to 0$ in $\Ll(V_0,V_0^\perp)$, as $r\searrow 0$, 
The map $r\mapsto D_r^0$ is continuous: indeed
$D_r^0 x= D(x+T_r x)$ and $T_r$ depends continuously on $r$ by Step~1.
Hence $r\mapsto D_r^0\in \Ll(V_0,Y)$ is a continuous
family of semi-Fredholm operators between fixed Hilbert spaces
and hence its semi-Fredholm index is constant as explained in Case~1
in the proof of Theorem~\ref{thm:Fredholm-index-constant}
for Fredholm operators; for semi-Fredholm operators
we refer to~\cite[\S 18 Thm.\,4]{Muller:2007a}.
\end{proof}

The proof of Theorem~\ref{thm:Fred-var-domain} is complete.
\end{proof}

\newpage
\section{Self-adjoint Hilbert space pair operators}
\label{sec:Hs.adj.ops}

\begin{theorem}\label{thm:ONB-basis}
Let $(H_0,H_1)$ be a Hilbert space pair.
Suppose the bounded linear map $A\colon H_1\to H_0$
is $H$-self-adjoint.\footnote{
  Fredholm of index $0$ and symmetric as unbounded operator
  on $H_0$ with dense domain~$H_1$.
  }
Then the following is true.
As unbounded operator on $H_0$ with dense domain $H_1$
the operator $A=A^*$ is selfadjoint.
The spectrum of $A$ consists of infinitely many discrete
real eigenvalues~$a_\ell$, of finite multiplicity each,\footnote{
  The \textbf{multiplicity} of an eigenvalue $a$ is the dimension of its
  eigenspace $\ker(A-a\,\Id)$.
  }
which accumulate either at $+\infty$,
or at $-\infty$, or at both.
Moreover, there exists a countable orthonormal basis
$\Vv(A)=\{v_\ell\}$ of $H_0$
composed of eigenvectors $v_\ell\in H_1$ of $A$.
\end{theorem}

In a Hilbert space pair both Hilbert spaces
are separable by~\cite[Cor.\,A.5]{Frauenfelder:2024c}.

\begin{proof}[Proof of Theorem~\ref{thm:ONB-basis}]
There are two cases for $A$, injective and not injective.

\smallskip\noindent
\textbf{Case~1:} $A$ is injective.

\smallskip\noindent
By the Fredholm property the image of $A$ is closed,
hence $(\im{A})^\perp=\coker{A}$.
Since the Fredholm index is zero and $A$ is injective we conclude
$\dim\coker{A}=\dim\ker{A}=0$.
Thus the operator $A\colon H_1\to H_0$ is surjective, hence bijective.
Since $A$ is also bounded the inverse $A^{-1}\colon H_0\to H_1$ is 
bounded, too, by the open mapping theorem. Composed with the compact
inclusion $\iota\colon H_1\to H_0$, the inverse as an operator on
$H_0$ is not only bounded, but even a compact operator with dense image
$$
   A^{-1}\colon H_0 \stackrel{\text{cp.}}{\longrightarrow} H_0
   ,\qquad
   \im{A^{-1}}=H_1\stackrel{\text{compact}\atop\text{dense}}{\INTO} H_0 .
$$
Now, by $H_0$-symmetry of $A$, the inverse $A^{-1}\in\Ll(H_0)$
is symmetric
$$
   \INNER{A^{-1}x}{y}
   =\INNER{A^{-1}x}{AA^{-1}y}
   =\INNER{AA^{-1}x}{A^{-1}y}
   =\INNER{x}{A^{-1}y}
   ,\quad
   \forall x,y\in H_0 ,
$$
which, by boundedness, is equivalent to self-adjointness
$(A^{-1})^*=A^{-1}\in\Ll(H_0)$.

To summarize, the inverse is a self-adjoint compact operator
$A^{-1}\colon H_0\to H_0$.
These are exactly the hypotheses of the Hilbert-Schmidt theorem,
see e.g.~\cite[thm.\,VI.16]{reed:1980a}, which asserts that
there is an orthonormal basis $\{v_k\}_{k\in\N}$ of $H_0$
such that $A^{-1}v_k=b_kv_k$ for non-zero real numbers $b_k\to 0$, as $k\to\infty$.
Moreover, the multiplicity of each eigenvalue $b_k$,
namely the dimension of its eigenspace
$\Eig_{b_k}(A^{-1}):=\ker(A^{-1}-b_k\,\Id)$, is finite.

Note that, while the list $(b_k)_{k\in\N}$ may contain finite
repetitions, there are still infinitely many different members.
Note further that, since $\im{A^{-1}}=H_1$, the eigenvectors $v_k\in H_0$
lie simultaneously in $H_1$: indeed $b_kv_k=A^{-1} v_k\in H_1$.
Hence we may apply $A$ to $A^{-1}v_k=b_kv_k$ and divide by $b_k$
to obtain
$$
   A v_k=a_k v_k
   ,\quad a_k:=\tfrac{1}{b_k}\in\R\setminus\{0\}
   ,\quad k\in\N
   ,\qquad \abs{a_k}\stackrel{k\to\infty}{\longrightarrow} \infty.
$$
Set $\Vv(A):=\{v_k\}_{k\in\N}\subset H_1$ to get an ONB of $H_0$
consisting of $A$-eigenvectors.

\smallskip
Self-adjointness $A=A^*$:
The operator $A^{-1}\in\Ll(H_0)$ satisfies
the hypothesis of~\cite[Thm.\,13.11 part~(b)]{rudin:1991a},
namely to be self-adjoint and injective.
The conclusion is that the operator inverse
$(A^{-1})^{-1}\colon H_0\supset \im{A^{-1}}\to H_0$ is self-adjoint.
This proves Theorem~\ref{thm:ONB-basis} for injective $A$ (Case~1).

\smallskip\noindent
\textbf{Case~2:} $A$ is not injective.

\smallskip\noindent
The linear map $A\colon H_0\supset H_1\to H_0$ decomposes as follows
\begin{equation}\label{eq:map-A}
\begin{tikzcd} 
H_0=\ker{A}\stackrel{\perp_0}{\oplus} X_0
  &(\im{A})^{\perp_0}\stackrel{\perp_0}{\oplus} X_0=H_0
\\
H_1=\ker{A}\oplus X_1
\arrow[u, hook,"\iota"', "{\text{dense}\atop\text{compact}}"]
\arrow[ur, "{A=0\oplus A|}"']
\end{tikzcd} 
\end{equation}
where
\begin{equation*}
\begin{split}
   X_0:=(\ker{A})^{\perp_0}{\color{gray}\,\subset H_0,}
   \quad
   X_1:=\iota^{-1}(X_0)=X_0\cap H_1 {\color{gray}\,\subset H_1,}
   \quad
   X_0=\im{A} .
\end{split}
\end{equation*}
We used that, by the Fredholm property, 
the kernel of $A$ is finite dimensional, so
a closed subspace of $H_0$, as well as of $H_1$.
Let $X_0$ be the orthogonal complement of $\ker{A}$ in $H_0$.
Orthogonal complements are closed subspaces.
Since $X_0$ is closed and $\iota$ is continuous, the pre-image
$X_0\cap H_1$ is a closed subspace of $H_1$.

Again by the Fredholm property,
the image of $A$ is closed, hence it too admits an orthogonal
complement which, by Fredholm index zero, is of the same
finite dimension as $\ker{A}$.
We show that $\im{A}=X_0$. 
'$\subset$' Given $y=Ax\in \im{A}$ and $z\in\ker{A}$,
by symmetry of $A$ we get
$\INNER{y}{z}_0=\INNER{Ax}{z}_0=\INNER{x}{Az}_0=\INNER{x}{0}_0=0$.
'$=$' Since the orthogonal complements $\ker{A}$ of $X_0$ and
$(\im{A})^{\perp_0}$ of $\im{A}$ are of the same finite dimension,
inclusion $\im{A}\subset X_0$ can only be true in case of equality
(otherwise the co-dimensions would be different).

We show that $H_1$ is the direct sum $\ker{A}\oplus X_1$. Note that
$\ker{A}\cap X_1=\ker{A}\cap X_0\cap H_1=\{0\}\cap H_1=\{0\}$
and $\ker{A}+ X_1=H_1$:
'$\subset$' obvious.
'$\supset$' Pick $x\in H_1$. Since $H_1\subset H_0=\ker{A}\oplus X_0$
write $x=x_*+x_0$ for unique elements $x_*\in\ker{A}$ and $x_0\in X_0$.
Then $x_0=x-x_*\in H_1\cap X_0=X_1$.

\smallskip\noindent
\textsc{Step~1:} The restriction $A|\colon X_0\supset X_1\to X_0$
meets the hypothesis of Case~1:
\begin{itemize}\setlength\itemsep{0ex} 
\item[\rm (a)]
  inclusion $\iota|\colon X_1\INTO X_0$ is compact
  and $X_1$ is a dense subset of $X_0$;
\item[\rm (b)]
  $A|$ is $X_0$-symmetric;
\item[\rm (c)]
  $A|\colon X_1\to X_0$ is a bounded bijection (hence Fredholm of index zero).
\end{itemize}

\begin{proof}[Proof of Step~1]
(a) Compactness: Let $B$ be a bounded subset of $X_1$.
Then $B$ is also subset of $H_1$, $H_0$, and $X_0$.
The closure of $B$ in $H_0$ is compact since $X_1\to H_1\to H_0$
is the composition of a bounded and a compact inclusion map,
hence itself compact.
But $X_0$ is a closed subspace of $H_0$ which contains $B$.
Thus the closure of $B$ is contained in $X_0$ as well.
\\
Density: The proof relies on $\ker{A}$ serving as finite dimensional
complement in \emph{both} $H_0$ and $H_1$.
Fix $x\in X_0{\color{gray}\;\subset H_0}$.
Since $H_1$ is dense in $H_0$, there exists a $H_0$-convergent sequence
$H_1\ni x_\nu\to x$. We use the orthogonal sum $H_0=\ker{A}\oplus X_0$
to write $x_\nu=c_\nu +z_\nu$ for unique $c_\nu\in\ker{A}\subset H_1$
and $z_\nu\in X_0$.
Now $z_\nu-x+c_\nu=x_\nu-x\to 0$ in $H_0$ and
$z_\nu=x_\nu-c_\nu\in H_1$. Thus $z_\nu\in X_0\cap H_1= X_1$.
Since $x_\nu-x=c_\nu +(z_\nu-x)$
with $c_\nu\in\ker{A}$ and $z_\nu-x\in X_0$ being $H_0$-orthogonal
Pythagoras provides the equality
$$
   \norm{c_\nu}_0^2+\norm{z_\nu-x}_0^2
   =\norm{x_\nu-x}_0^2 \stackrel{\nu\to\infty}{\longrightarrow} 0.
$$
This proves $H_0$-convergence $H_1\ni z_\nu \to x\in X_0$
and concludes the proof of~(a).

\smallskip\noindent
(b) Since $X_1\subset H_1$ and $X_0\subset H_0$,
part~(b) is true by $H_0$-symmetry of $A$.

\smallskip\noindent
(c) Injective and surjective are obvious.
The restriction of a bounded linear map to a closed
subspace is bounded.
\end{proof}

\noindent
\textsc{Step~2:} We prove Theorem~\ref{thm:ONB-basis}.

\begin{proof}[Proof of Step~2]
We decompose $A=0\oplus A|$ into two summands as in~(\ref{eq:map-A}).

\smallskip\noindent
\textsc{Summand $A|\colon X_0\supset X_1\to X_0$.}
By Step~1 the restriction $A|$
meets the hypothesis of Case~1. Thus $A|$ is self-adjoint
as an unbounded operator and its spectrum $\spec{A|}$
consists of infinitely many discrete real eigenvalues $a\not=0$
of finite multiplicity each, which accumulate either at $+\infty$,
or at $-\infty$, or at both.
Moreover, there is an ONB
$\Vv(A|)=\{v_k\}_{k\in\N}\subset X_1$ of $X_0$
consisting of eigenvectors of $A|$.

\smallskip\noindent
\textsc{Summand $0\colon \ker{A} \to (\im{A})^{\perp_0}$.}
The spectrum consists of the eigenvalue $0$.
The dimension of the eigenspace $\ker{A}$
is at least $1$ (Case~2) and finite (Fredholm assumption).
Choose an $H_0$-ONB of $\ker{A}$, notation $\Vv(\ker{A})$.

\smallskip
To see that $A\colon H_0\supset H_1\to H_0$ is self-adjoint,
unpack the definition of the domain of an adjoint operator
to get the first identity
$$
   \dom{A^*}
   =\ker{A}\oplus D(A|^*)
   =\ker{A}\oplus D(A|)
   =\ker{A}\oplus X_1
   =\dom{A}
$$
whereas the second identity holds since $A|$ is self-adjoint by Case~1.

The spectrum of $A$ is the union of the spectrum of $A|$ and $\{0\}$.
The union
$$
   \Vv(A)
   :=\Vv(\ker{A})\cup
   \Vv(A|)
$$
consists of eigenvectors of $A$. It is
an ONB of $H_0$ (eigenvectors to different
eigenvalues are orthogonal since $A=A^*$).
This proves Step~2 and Case~2.
\end{proof}

This concludes the proof of Theorem~\ref{thm:ONB-basis}.
\end{proof}

\newpage
\bibliographystyle{alpha}
\addcontentsline{toc}{section}{References}
\bibliography{$HOME/Dropbox/0-Libraries+app-data/Bibdesk-BibFiles/library_math,$HOME/Dropbox/0-Libraries+app-data/Bibdesk-BibFiles/library_math_2020,$HOME/Dropbox/0-Libraries+app-data/Bibdesk-BibFiles/library_physics}{}

\end{document}

%% file: joa-environments-english.tex
\newtheorem{theoremABC}{Theorem}

\newtheorem{theorem}{Theorem}[section]

\newtheorem{corollary}[theorem]{Corollary}

\newtheorem{lemma}[theorem]{Lemma}
\newtheorem{proposition}[theorem]{Proposition}

\theoremstyle{definition}
\newtheorem{definition}[theorem]{Definition}

\newtheorem{remark}[theorem]{Remark}

\theoremstyle{remark}
\newtheorem*{notation}{Notation}


%% file: joa-latex-shortcuts.tex
\newcommand{\A}{{\mathbb{A}}}

\newcommand{\D}{{\mathbb{D}}}

\newcommand{\N}{{\mathbb{N}}}

\newcommand{\Q}{{\mathbb{Q}}}
\newcommand{\R}{{\mathbb{R}}}

\newcommand{\Z}{{\mathbb{Z}}}
\newcommand{\Aa}{{\mathcal{A}}}   
\newcommand{\Bb}{{\mathcal{B}}}
\newcommand{\Dd}{{\mathcal{D}}}
\newcommand{\Ee}{{\mathcal{E}}}
\newcommand{\Ff}{{\mathcal{F}}}

\newcommand{\Ii}{{\mathcal{I}}}

\newcommand{\Kk}{{\mathcal{K}}}
\newcommand{\Ll}{{\mathcal{L}}}   

\newcommand{\Rr}{{\mathcal{R}}}
\newcommand{\Ss}{{\mathcal{S}}}
\newcommand{\Tt}{{\mathcal{T}}}

\newcommand{\Vv}{{\mathcal{V}}}
\newcommand{\Ww}{{\mathcal{W}}}

\newcommand{\Ann}{{\rm Ann }}            
\newcommand{\coker}{{\rm coker\, }}  
\newcommand{\im}{{\rm im\, }}             
\newcommand{\dom}{{\rm dom\, }}           
\newcommand{\Id}{{\rm Id}}
\newcommand{\codim}{{\rm codim\, }}       
\newcommand{\supp}{{\rm supp\, }}         
\newcommand{\INDEX}{\mathop{\mathrm{index}}}     

\newcommand{\cgraph}[1]{\Gamma_{\kern-.5ex{}#1}}     
\newcommand{\spec}{\mathrm{spec}\,}        
\newcommand{\Eig}{\mathrm{Eig}}            
\newcommand{\norm}{{\rm norm}}

\newcommand{\eps}{{\varepsilon}}

\renewcommand{\O}{{\rm O}}     

\newcommand{\inner}[2]{\langle #1, #2\rangle}   
\newcommand{\INNER}[2]{\left\langle #1, #2\right\rangle}

\def\NABLA#1{{\mathop{\nabla\kern-.5ex\lower1ex\hbox{$#1$}}}}
\def\Nabla#1{\nabla\kern-.5ex{}_{#1}}
\def\Tabla#1{\Tilde\nabla\kern-.5ex{}_{#1}}
\def\Babla#1{\widebar\nabla\kern-.5ex{}_{#1}}
\def\abs#1{\mathopen|#1\mathclose|}   
\def\Abs#1{\left|#1\right|}            
\def\norm#1{\mathopen\|#1\mathclose\|}
\def\Norm#1{\left\|#1\right\|}

\renewcommand{\Tilde}{\widetilde}

\newcommand{\p}{{\partial}}

\newcommand{\INTO}{\hookrightarrow}              

\renewcommand{\1}{{{\mathchoice {\rm 1\mskip-4mu l} {\rm 1\mskip-4mu l}
{\rm 1\mskip-4.5mu l} {\rm 1\mskip-5mu l}}}}
\renewcommand{\graph}{{\rm graph }}  

\newlength\eqshift
\renewcommand\theequation{\thesection.\arabic{equation}}
\let\savetheequation\theequation

\makeatletter
\renewcommand*\env@matrix[1][\arraystretch]{%
  \edef\arraystretch{#1}%
  \hskip -\arraycolsep
  \let\@ifnextchar\new@ifnextchar
  \array{*\c@MaxMatrixCols c}}
\makeatother

\makeatletter
\let\save@mathaccent\mathaccent
\newcommand*\if@single[3]{%
  \setbox0\hbox{${\mathaccent"0362{#1}}^H$}%
  \setbox2\hbox{${\mathaccent"0362{\kern0pt#1}}^H$}%
  \ifdim\ht0=\ht2 #3\else #2\fi
  }
\newcommand*\rel@kern[1]{\kern#1\dimexpr\macc@kerna}
\newcommand*\widebar[1]{\@ifnextchar^{{\wide@bar{#1}{0}}}{\wide@bar{#1}{1}}}
\newcommand*\wide@bar[2]{\if@single{#1}{\wide@bar@{#1}{#2}{1}}{\wide@bar@{#1}{#2}{2}}}
\newcommand*\wide@bar@[3]{%
  \begingroup
  \def\mathaccent##1##2{%
    \let\mathaccent\save@mathaccent
    \if#32 \let\macc@nucleus\first@char \fi
    \setbox\z@\hbox{$\macc@style{\macc@nucleus}_{}$}%
    \setbox\tw@\hbox{$\macc@style{\macc@nucleus}{}_{}$}%
    \dimen@\wd\tw@
    \advance\dimen@-\wd\z@
    \divide\dimen@ 3
    \@tempdima\wd\tw@
    \advance\@tempdima-\scriptspace
    \divide\@tempdima 10
    \advance\dimen@-\@tempdima
    \ifdim\dimen@>\z@ \dimen@0pt\fi
    \rel@kern{0.6}\kern-\dimen@
    \if#31
      \overline{\rel@kern{-0.6}\kern\dimen@\macc@nucleus\rel@kern{0.4}\kern\dimen@}%
      \advance\dimen@0.4\dimexpr\macc@kerna
      \let\final@kern#2%
      \ifdim\dimen@<\z@ \let\final@kern1\fi
      \if\final@kern1 \kern-\dimen@\fi
    \else
      \overline{\rel@kern{-0.6}\kern\dimen@#1}%
    \fi
  }%
  \macc@depth\@ne
  \let\math@bgroup\@empty \let\math@egroup\macc@set@skewchar
  \mathsurround\z@ \frozen@everymath{\mathgroup\macc@group\relax}%
  \macc@set@skewchar\relax
  \let\mathaccentV\macc@nested@a
  \if#31
    \macc@nested@a\relax111{#1}%
  \else
    \def\gobble@till@marker##1\endmarker{}%
    \futurelet\first@char\gobble@till@marker#1\endmarker
    \ifcat\noexpand\first@char A\else
      \def\first@char{}%
    \fi
    \macc@nested@a\relax111{\first@char}%
  \fi
  \endgroup
}
\makeatother

\long\def\symbolfootnote[#1]#2{\begingroup%
\def\thefootnote{\fnsymbol{footnote}}\footnote[#1]{#2}\endgroup}

\tikzset{
  symbol/.style={
    draw=none,
    every to/.append style={
      edge node={node [sloped, allow upside down, auto=false]{$#1$}}}
  }
}

\newcommand{\gf}[1]{#1}                       

